\pdfoutput=1
\RequirePackage{fix-cm}
\documentclass[smallextended]{svjour3}       
\smartqed  
\usepackage{graphicx}
 \usepackage{mathptmx}      
\usepackage{amsmath}
 \usepackage{natbib}
 \usepackage{algorithm}
\usepackage{algorithmic}
\usepackage{amssymb}
\usepackage{float}
\usepackage{amsthm}
\usepackage{subcaption}
\usepackage{hyperref}
\newcommand{\comment}[1]{}
\usepackage{xcolor}

\usepackage{bmpsize}
\usepackage{adjustbox}

\newtheorem{observation}{Observation}

\begin{document}

\title{Instance-specific linear relaxations of semidefinite optimization problems \thanks{This material is based upon work supported by the U. S. Office of Naval Research under award number N00014-21-1-2243.}
}




\author{Daniel de Roux
\and
Robert Carr 
\and
R.Ravi
}

\authorrunning{D. de Roux et al.}

\institute{Daniel de Roux  \at
           	{Tepper School of Business, Carnegie Mellon University} \\
              \email{dderoux@andrew.cmu.edu}\\
             \and
            Robert Carr \at 
            {UNM}\\
              \email{bobcarr@swcp.com}    \\
             \and 
           R. Ravi  \at
            	{Tepper School of Business, Carnegie Mellon University} \\
              \email{ravi@andrew.cmu.edu}\\
}

\date{Received: date / Accepted: date}

\maketitle

\begin{abstract}
We introduce a generic technique to obtain linear relaxations of semidefinite programs with provable guarantees based on the commutativity of the constraint and the objective matrices. 
We study conditions under which the optimal value of the SDP and the proposed linear relaxation match, which we then relax to provide a flexible methodology to derive effective linear relaxations.
We specialize these results to provide linear programs that approximate well-known semidefinite programs for the max cut problem proposed by Poljak and Rendl, and the Lov\'asz theta number; we prove that the linear program proposed for max cut certifies a known eigenvalue bound for the maximum cut value and is in fact stronger. 
Our ideas can be used to warm-start algorithms that solve semidefinite programs by iterative polyhedral approximation of the feasible region. We verify this capability through multiple experiments on the max cut semidefinite program, the Lov\'asz theta number and on three families of semidefinite programs obtained as convex relaxations of certain quadratically constrained quadratic problems.

\end{abstract}

\section{Introduction}
\label{revised_sec:1}

Semidefinite optimization, i.e., the optimization of a linear function over the set of positive semidefinite matrices intersected with an affine subspace \cite{vandenberghe1996semidefinite}, is one of the most active research areas in convex mathematical optimization. The generic formulation for a semidefinite optimization problem (SDP) is

\begin{equation}\label{SDP}\tag{SDP}
\begin{aligned}
\min_{X\in \mathbb{S}^{ n}} \  &\langle  C,X\rangle \\
\text{ s.t: } \langle A_i,X \rangle  = & \  b_i,\  \forall i \in [r],\\ 
X  \succeq & \ 0
\end{aligned}
\end{equation}

\noindent where $\mathbb{S}^n$ denotes the set $n\times n$ symmetric matrices, $C \in \mathbb{S}^{n}$ is a symmetric (without loss of generality) cost matrix, $\left \langle \cdot, \cdot \right \rangle$ denotes the Frobenius inner product, $r \in \mathbb{N}$, $[r]$ denotes the set of integers $\{1,\dots,r\}$ and a symmetric $n \times n$ matrix $X$ is \textit{positive semidefinite}, denoted $X \succeq 0 $, if and only if $v^\top X v \geq 0$ for all $v \in \mathbb{R}^n$.

SDPs arise naturally in combinatorial optimization \cite{alizadeh1995interior,goemans1995improved,Lovasz1979shannon,nesterov1997semidefinite}, control theory \cite{ahmadi2012algebraic,henrion2005positive,parrilo2000structured}, polynomial optimization \cite{lasserre2001global,parrilo2000structured,parrilo2003semidefinite}, machine learning \cite{d2004direct,lanckriet2004learning} and are solvable in polynomial time up to an arbitrary accuracy via the theory of interior-point methods \cite{nemirovski2004interior}. Nonetheless,
it is well known that SDPs are challenging to solve in practice. Typical off-the-shelf solvers use interior-point methods, which require computation of large Hessian matrices (and their inverses) and are often intractable due to memory limitations. For an illustration, see  \cite[chapter 6.7]{ben2001lectures}, \cite{majumdar2020recent}, and \cite{bertsimas2020polyhedral} where it is mentioned that state-of-the-art solver such as MOSEK \cite{mosek} cannot solve semidefinite problems with a symmetric matrix $X$ on more than $250$ rows. Inspired by these practical limitations, researchers have proposed a several ideas to solve large-scale semidefinite programs. Among them, we have $i)$ exploiting structure of the problem (such as sparsity and symmetry), $ii)$ producing low rank solutions, $iii) $ algorithms based on augmented Lagrangians and the alternating direction method of multipliers, and $iv$) approaches that trade scalability with conservatism  by iteratively finding inner and outer polyhedral approximations of the semidefinite problem. See \cite{majumdar2019survey} for a survey of all of these methods. \\

Although all of these techniques enjoy a rich literature, the algorithms of $iv)$ are of special theoretical interest. The cornerstone of these methods is constructing inner and outer polyhedral approximations of the semidefinite cone in order to find a sequence of improving feasible solutions together with tighter bounds on the objective of the SDP, allowing one to trade off between scalability and conservatism. 

Research on this class of algorithms is relevant due to its intimate connection with a fundamental question in convex geometry: can the positive semidefinite cone be approximated by polyhedra? Taking the perspective of the field of optimization, this question can be framed by asking if linear programs are strong enough to approximate semidefinite ones.
These twin questions, relevant in the fields of optimization and convex geometry respectively, have given rise to a thriving body of research ~\cite{ahmadi2017sum,bertsimas2020polyhedral,charikar2009integrality, fawzi2021polyhedral}.

Outer approximations have been the focus of substantial effort since the hardness of \ref{SDP} comes from the semidefinite constraint and so one may drop it and and add linear constraints on $X$ implied by $X\succeq 0$. In this case, (\ref{SDP}) is relaxed to a linear program. A typical example is to add the constraints $X_{i,i}\geq 0 ,\ \forall i \in [n]$
and $X_{ii}+X_{jj}\pm 2X_{i,j}\geq 0, \ \forall i\in [n], \ \forall j \in [n]$ which are valid for any $X\succeq 0$. These relaxations tend to be weak and seldom used in practice \cite{braun2015approximation,chan2016approximate}. A well studied example of this phenomenon is the maximum cut problem and the theoretical hardness of approximating it with linear programs 
which we will discuss in depth in Section \ref{revised_sec:3}.

The previous approach can be improved using ideas of Kelley \cite{kelley1960cutting}. The strategy is to sequentially refine the linear relaxations by aggregation of cutting planes. More concretely, consider the linear relaxation of \ref{SDP} given by

\begin{equation}\label{LSDP}\tag{$L_\mathcal{S}$}
\begin{aligned}
\min_{X\in \mathbb{S}^{n}} \  &\langle  C,X\rangle \\
\text{ s.t: } \langle A_i,X \rangle  = & \  b_i,\  \forall i \in [r],\\ 
v^\top Xv  \geq & \ 0 \ \forall \ v \in \mathcal{S}
\end{aligned}
\end{equation}

where $\mathcal{S}$ is a finite subset of $\mathbb{R}^n$. 
 Here, we simply insist that $v^\top X v \geq 0 $ only for the elements $v$ of the set $\mathcal{S}$. If a solution to this program is not positive semidefinite, we may update $\mathcal{S}$ iteratively. This results in the following algorithm:

\begin{algorithm}
\caption{}
\label{alg1}
\begin{algorithmic}[1] 
\STATE Fix a finite set $\mathcal{S} \subseteq \mathbb{R}^n$. Drop the semidefinite constraint $X\succeq 0$ of program $\ref{SDP}$ and solve the resulting linear program \ref{LSDP} finding a minimizer $X^*$. 
\WHILE{$X^*$ has a negative eigenvalue}
\STATE Find a eigenvector $v$ corresponding to the most negative eigenvalue of $X^*$. Add $v$ to $\mathcal{S}$.
\STATE Solve the updated linear program to find a new minimizer $X^*$.
\ENDWHILE
\RETURN $X^*$.
\end{algorithmic}
\end{algorithm}

The specific implementations of this algorithm mainly differ in how one updates the set $\mathcal{S}$.
In \cite{ahmadi2017sum,ahmadi2019dsos},the authors use the extreme rays of the set of diagonally dominant matrices, which are then rotated by matrices obtained from a Cholesky decomposition of an optimal solution to the dual of (\ref{LSDP}). They also propose an inner approximation of the positive semidefinite cone based on the so-called $DSOS_n$ and $SDSOS_{n,d}$ cones. In a different line of work \cite{baltean2019scoring,dey2021cutting,qualizza2012linear,sherali2002enhancing,wang2021polyhedral} chose the elements $v$ of $\mathcal{S}$ favoring sparsity, with the idea that the resulting linear programs will be easier to solve. 
Bundle methods, such as the \textit{spectral bundle method} of Helmberg and Rendl \cite{helmberg2000spectral}
work with the dual of (\ref{SDP}), under the further restriction that $X$ has a constant trace. \cite{krishnan2006unifying} presents a unifying framework for the latter and similar methods. In \cite{bertsimas2020polyhedral}, the constraint $X\succeq 0$ is replaced for infinitely many constraints of the form $f(X,Y)\leq 0$ which must hold for every $Y$ in some convex set $\mathcal{Y}$ and where $f$ is a Lipschitz continuous function. The authors further argue that one should instead solve a second-order cone relaxation, adding the constraints 
\[
\left\| \begin{pmatrix}2X_{i,j} \\ X_{i,i}-X_{j,j} \end{pmatrix} \right\|_2 \leq X_{i,i}+X_{j,j}, \ \forall i \in [n], \  \forall j \in [n],
\]
 which are valid for (\ref{SDP}). \\

It is noteworthy that mostly all of these works discuss how to update $\mathcal{S}$, but seldom consider how to initialize it. Typically $\mathcal{S}$ is set to the standard basis of $\mathbb{R}^n$, resulting in the linear constraints $X_{ii} \geq 0, \  i \in \{1,\dots,n \},$ which are implied by the constraint $X\succeq 0$. Interestingly, under mild conditions, there exists a finite set $\mathcal{S}$ that ensures that the optimal values of the SDP and the linear relaxation \ref{LSDP} match, supporting the approach of using Algorithm\ref{alg1}. 

\begin{observation}\label{OptimalS}
Suppose that both (\ref{SDP}) and its dual, given by the following semidefinite optimization program 

\begin{equation}\label{DSDP}\tag{DSDP}
\begin{aligned}
&\max_{y\in \mathbb{R}^{r}} \  b^\top y \\
\text{ s.t: }& C-  \sum_{i=1}^r y_iA_i  \succeq  \ 0
\end{aligned}
\end{equation}

\noindent are strictly feasible. Let $y^*$ be an optimal solution to (\ref{DSDP}).  Let $v_1,\dots,v_n$ be an orthonormal basis of $\mathbb{R}^n$ of eigenvectors of $C-  \sum_{i=1}^my^*_iA_i = S^*$ with $S^* = \sum_{i=1}^n \beta_ivv^\top$, and $\beta_i$ the eigenvalues of $S^*$. Let $\mathcal{S}^*=\{v_1,\dots,v_n\}$. Then, $L_{\mathcal{S}^*}$ is solvable, and its optimal value matches the optimal value of  (\ref{SDP}).
\end{observation}
The proof of this observation is deferred to the Appendix \ref{AppendixA}.
Similar versions of Observation \ref{OptimalS} can be found in \cite{krishnan2006unifying} and \cite{sivaramakrishnan2002linear}.
In fact, \cite{sivaramakrishnan2002linear} proves that if $\mathcal{S} = \{u_1,\dots,u_l\}$ are the vectors generated by the spectral bundle method of \cite{helmberg2000spectral} of Rendl et al. to solve \ref{DSDP}, the objective value of $(L_{\mathcal{S}})$ matches that of $(\ref{SDP})$, but this is hardly surprising: if we knew in advance the set of vectors $\mathcal{S}^*$ given by Observation \ref{OptimalS}, we could set $\mathcal{S} = \mathcal{S}^*$ and solve (\ref{SDP}) as a linear program. More importantly, we emphasize that finding the sets $\mathcal{S}^*$ and $\{u_1,\dots,u_l\}$ requires solving another comparable SDP, namely \ref{DSDP}.

In this paper, we tackle the task of finding a better set $\mathcal{S}$ to initialize Algorithm \ref{alg1} under certain computational restrictions by drawing inspiration from the question of when - if ever- one can avoid the iterative procedure suggested by Kelley and exactly solve the semidefinite program with a linear program. 
By ``exactly solving" we mean finding a linear relaxation of the SDP whose optimal value equals that of the SDP.

Technically, Observation \ref{OptimalS} indicates that the question of \textit{exactly} solving an SDP with a linear problem is ill-posed if one does not restrict the set of \textit{algorithms} one is allowed to use to process the instance. We can consider at least three possible approaches to amend this issue. First,  restricting the access one has to the given instance. For example, say we are not shown a full SDP instance, but one is allowed to sample a small subset of the entries of the objective and constraints matrices. Second, to only have access to algorithms with at most a certain computational complexity, say matrix multiplication complexity. However, this would require fixing a concrete computational model and proving lower bounds for the complexity of the algorithms to be used, which are typically very hard to obtain. A third approach, which we take in this paper, is to fix an oracle $\mathcal{O}$, that we can query at most a constant number of times. Concretely, we will assume that we have at our disposal an oracle that can compute a eigenvector decomposition of a symmetric matrix, and that can solve linear programs of polynomial size.
If the SDP can in fact be solved with such an oracle, we say it is  \textit{ solvable under } $\mathcal{O}$.

\subsection{Hardness of approximation of the max cut problem}

The question of finding a good set $\mathcal{S}$ to initialize Algorithm \ref{alg1} amounts to finding a linear approximations to a semidefinite programs together with a guarantee that the approximation is good. This line of research is motivated by the question of whether the maximum cut (max cut henceforth) problem can be approximated using a linear program by a factor strictly better than $2$. This problem consists in 
finding a bipartition of the nodes of a given graph that maximizes the number of edges with one end in both parts. The results of Poljak, Rendl, Goemans and Williamson \cite{goemans1995improved,poljak1995nonpolyhedral} show that max cut can be approximated to within a factor of $~1.13$ by an SDP relaxation. Therefore, a linear approximation of factor at most $~1.769$ to that SDP would result in a linear approximation the the max cut problem with an approximation better than $2$ \footnote{ Since  $1.77 \cdot 1.13 =2$.}. Such a result would be striking as the common belief is that max cut cannot be approximated within a factor better than $2$ with a linear program in the restricted case that the feasible region of the program is independent of the graph and solely depends on the number of vertices \cite{braun2015approximation,chan2016approximate,charikar2009integrality,de2007linear,kothari2021approximating}. In Section \ref{revised_sec:3}, we explore in detail the hardness of approximation results for max cut.

Drawing inspiration from the study of exact solvability of an SDP with an LP, we make the case that we can obtain ``good starting" linear approximations for semidefinite programs if one is allowed to let $\mathcal{S}$ depend on the dual of the semidefinite program. The heart of the argument is that the obstructions mentioned for max cut emerge specifically when the polytopes being optimized are determined solely by the number of variables (node pairs for max cut) in a given instance. Hence, we propose to let $\mathcal{S}$ depend on the matrices $C$ and $A_1,\dots,A_r$ which determine the objective and the constraints of $(\ref{SDP})$, and consequently on the feasible region of \ref{DSDP}. Crucially, such formulations trivially avoid the results in \cite{braun2015approximation} and \cite{kothari2021approximating}. We call linear approximations with such dependence ``instance-specific". Notice that making some assumption on the algorithms that we can use to interact with the instance is essential here. To illustrate this point, imagine we wish to write a linear program to find the max cut value $mc(G)$ of a graph $G$. To do so, we can compute a max cut of the graph using brute force and then write an LP with a linear constraint insisting that the objective equals $mc(G)$.

\subsection{Exact linear relaxations under $\mathcal{O}$}

To find candidate sets $\mathcal{S}$ that guarantee that the linear program \ref{LSDP} is a strong relaxation of \ref{SDP} we first explore sufficient conditions under which the SDP is solvable under the oracle $\mathcal{O}$. 
Although Observation \ref{OptimalS} suggests an answer, such a set of vectors cannot, as far as we are aware, be obtained with the oracles we are considering. 
In Section $2$, we present Theorems \ref{theorem:poly_exact} and \ref{teo:diag_exactness} which will provide solvability under $\mathcal{O}$ without requiring the solution of a semidefinite program. Our results are tied to the geometry of the dual feasible region of \ref{SDP}, and a relevant case is when the dual feasible region is a polyhedron. If such is the case and an explicit description of it is available, then program \ref{DSDP} can be solved as a linear program. Theorem \ref{theorem:poly_exact} shows that under the same condition the \textit{primal} \ref{SDP} can be solved with a linear program as well. Unfortunately, this theorem is not very useful as it requires enumerating the vertices of the feasible region, which may grow exponentially. 
The polyhedral assumption has received attention from the literature in the context of \textit{ quadratically constrained quadratic problems} (QCQPs) \cite{wang2022tightness}, and perhaps more so a weakening of it: simultaneous diagonalizability.

\begin{definition}
    A set of matrices $\{A_i\}_{i \in I} \subseteq \mathbb{R}^{n \times n} $ where $I$ is some set of indices which may be infinite, is said to be simultaneously diagonalizable (SD) if there exists an invertible, orthogonal matrix $U \in \mathbb{R}^n$ such that every element of the set $\{U^\top A_i U\}_{i \in I}$ is a diagonal matrix. Note that $U^\top U = UU^\top = I_n$ as $U$ is orthogonal.
\end{definition}

It turns out that if the set of matrices
defining the dual feasible region $\Gamma$ of \ref{SDP} is simultaneously diagonalizable, then $\Gamma$ is a polyhedron \cite{wang2022tightness}.

\begin{observation}\label{obs:simulDiag}
Let $\Gamma$ be a spectrahedron given by the representation  $\Gamma = \{y\in \mathbb{R}^n: C-\sum_{i}^r A_iy_i )\succeq 0  \}$. If the set of matrices $\{C, \{A_i\}_{ i \in [r]} \} $ is simultaneously diagonalizable, then $\Gamma$ is polyhedral.
\end{observation}

We prove this fact in Section \ref{revised_sec:2}, and point out that the given condition is sufficient but not necessary. Under this more stringent condition, we  prove in Theorem \ref{teo:diag_exactness} that $\mathcal{O}$ can be used to solve \ref{SDP}.

It will typically not be the case that the dual feasible set $\Gamma$ is polyhedral, and much less that the matrices $C, \{A_i\}_{i \in  [r]}$ are simultaneously diagonalizable. In Section \ref{revised_sec:2} we prove that this condition is equivalent 
to the simultaneous diagonalizability of matrices $C- \sum_{i}A_ip_i$ and
$C-\sum_i A_iq_i$ for all $p$ and $q$ in $\mathbb{R}^r$. This characterization suggests that we only insist of the commutativity of the matrices $C-\sum_i A_ip_i$ and $C-\sum_i A_iq_i$ for some $p$ and $q$. It turns out that this is the key idea to initialize the set $\mathcal{S}$ in Algorithm \ref{alg1}. 
In Section \ref{revised_sec:2} we set the theoretical background of these considerations, and in the following sections we explore their applications to three families of semidefinite optimization problems: the max cut problem, The Lov\'asz theta number and the more generic Shor SDP relaxation of quadratically constrained quadratic problems.

We stress that the intention of the presented approach is to further explore when an SDP can be solved with a linear program, and to improve on existing cutting plane approaches to solve SDPs (such as the conservative methods described in \cite{majumdar2019survey}). This family of methods is \textit{not} the de-facto choice to solve large scale semidefinite programs, and very strong methods exist which can scale substantially such as \cite{o2016conic,yang2015sdpnal+,zhao2010newton,wen2013feasible,wang2023decomposition}. Nevertheless, We point out that these methods might come with their own limitations and in settings where SDPs appear naturally, such as in the sum-of-squares hierarchy for polynomial optimization \cite{yang2023inexact}, or whenever optimal solutions to the SDPs are not low rank. In these regimes, polyhedral approximations might be a good alternative. In addition, developing stronger polyhedral approximations to SDPs has consequences in approaches to integer semidefinite programs, which has received attention recently \cite{coey2020outer,gally2016computing,gally2018framework,   hojny2023handling,yonekura2010global} and in spatial branch-and-bound algorithms for non-convex quadratic problems.

\subsection{Overview and outline of this paper}

 \begin{itemize}
     \item[(a)] In Sect. \ref{revised_sec:2} we derive two sufficient conditions for solvability of an SDP under $\mathcal{O}$. These conditions are then weakened to produce a strategy to provide candidate starting sets $\mathcal{S}_0$ for outer polyhedral approximation algorithms to solve SDPs.
 
     \item[(b)] In Sect. \ref{revised_sec:3}, we study the setting of finding a maximum cut of a graph $G$ using the semidefinite relaxation of Poljak, Rendel, Goemans and Williamson \cite{goemans1995improved,poljak1995nonpolyhedral}. Even though the conditions for exact solvability are not met, we use the relaxed version to provide a linear program that certifies a spectral bound in contrast to previous linear relaxations for the maximum cut problem. We then derive a solvability result under $\mathcal{O}$, recovering and generalizing a theorem of Alon and Sudakov \cite{alon2000bipartite}.

     \item[(c)] In Sect. \ref{revised_sec:4} we introduce linear relaxations of the Lov\'asz theta number SDP and Shor's semidefinite relaxation for quadratically constrained quadratic programs. We recall as well our linear relaxation of max cut, and introduce a linear strengthening of the max cut SDP.
     
     \item[(d)] In Sect. \ref{revised_sec:5} we extensively test our methods empirically on random instances of the problems introduced in Section \ref{revised_sec:4}.  
     We discuss solving times of the proposed programs.
     
     \item[(e)] In Appendix \ref{AppendixA}, we prove Observation \ref{OptimalS} and provide an alternative proof of Lemma \ref{eigenBound}. In Appendix B we show the performance of our linear program in the case where the original SDP is itself a relaxation of an underlying optimization problem. We study the case of the max cut problem and the \textit{sparse PCA } problem, where both the SDPs and our linear relaxations can be used to recover a solution to the underlying problem. We show that the quality of our linear programs is competitive with that of the SDPs.
     For max cut, we compare with results obtained by Mirka and Williamson in \cite{mirka2022experimental}. 
     
 \end{itemize}

\subsection{Notation}
We denote the set of square, real,  $n\times n$ symmetric matrices by $\mathbb{S}^n$. We denote the cardinality of a set $I$ by $|I|$.
We denote by $e_1,\dots, e_n$ the standard basis of $\mathbb{R}^n$ and the $n\times n$ identity matrix by $I_n$. For a symmetric matrix $W$ we let $\lambda_1(W) \geq \lambda_2(W) \geq \dots \geq \lambda_n$(W) be its eigenvalues. When the matrix is clear from the context, we drop the terms in parentheses and simply write $\lambda_1\geq\dots\geq\lambda_n$. For $A \in \mathbb{S}^n$ we write $tr(A)$ for the trace of $A$: $tr(A) = \sum_{i=1}^nA_{ii}$ and write $\|A\|_F$ to denote the Frobenius norm of $A$: $\|A\|_F = \sqrt{\sum_{i=1}^n\sum_{j=1}^n A_{ij}^2}$. The $\ell_1$ norm of $A$ is given by $\|A\|_1 = \sum_{i,j}|A_{ij}|$. We denote by $\langle\cdot,\cdot\rangle$ the usual Frobenius inner product of two matrices in $\mathbb{S}^n$, recalling that for two matrices $A,B \in \mathbb{S}^n, \ \langle A,B\rangle = tr(A^TB)=tr(AB)$.
We denote by $\Vec{1}$ the vector of all ones in $\mathbb{R}^n$ and by $J$ the matrix of all ones. If $A$ is a matrix, we denote by $diag(A)$ the vector given be the diagonal of $A$. If $u$ is a vector, $diag(u)$ denotes the matrix with $u$ on its diagonal. We
denote by $\mathcal{E}(A)$ an arbitrary orthonormal basis consisting of eigenvectors of $A$. In particular, if $A \in \mathbb{S}^n$ and $\mathcal{E}(A) = \{v_1,\dots,v_n\}$ then we have $A = \sum_{i=1}^n \lambda_iv_iv_i^\top$ \cite{horn2012matrix}. Finally, given a weighted graph $G$ we denote respectively the value of the max cut of $G$, the adjacency matrix, the number of edges and the laplacian matrix by $mc(G), W(G)$, $m$ and $\mathcal{L}(G)$. If $G$ is clear from the context, we drop the dependency on $G$ and simply write $mc,W,m$ and $\mathcal{L}$.

\section{Instance-specific linear relaxations of semidefinite optimization problems}
\label{revised_sec:2}

In this section we explore the question of exact solvability of semidefinite programs given access to an oracle $\mathcal{O}$, with the following properties:
\begin{itemize}
    \item Given a set of simultaneously diagonalizable matrices $\{A_1,\dots,A_r\}$, $\mathcal{O}$ can be called once to compute an orthogonal matrix $U$ such that $U^\top A_i U$ are diagonal matrices for $i=1,\dots r$. For an implementation of such an oracle see \cite{golub2013matrix}.
    \item $\mathcal{O}$ can be called a constant number of times to find an optimal solution to a linear program of polynomial size in the bit representation of the information of the SDP, namely the objective and constraint matrices.
\end{itemize}

In case we can find the optimal value of program \ref{SDP} by querying   $\mathcal{O}$ at most a constant number of times, we say that the SDP is solvable under $\mathcal{O}$, and our intention is to derive sufficient conditions that guarantee solvability of the SDP. It is to be expected that such conditions are not applicable except in some rare cases. We posit that we can derive weakenings of them to provide a starting set $\mathcal{S}$ for Algorithm \ref{alg1}.
Recall that a generic SDP is given by 

\begin{equation}\tag{SDP}
\begin{aligned}
\min_{X\in \mathbb{S}^{n}} \  &\langle  C,X\rangle \\
\text{ s.t: } \langle A_i,X \rangle  = & \  b_i,\  \forall i \in [r],\\ 
X  \succeq & \ 0.
\end{aligned}
\end{equation}

The dual of this program is:

\begin{equation}\tag{DSDP}
\begin{aligned}
&\max_{y\in \mathbb{R}^{n}} \  b^\top y \\
\text{ s.t: }& C-  \sum_{i=1}^my_iA_i  \succeq  \ 0.
\end{aligned}
\end{equation}

Throughout this paper, we will assume  ''generic SDPs" and their duals are strictly feasible, and therefore strong duality holds.
A spectrahedron $\Gamma$ is the intersection of the cone of positive semidefinite matrices and an affine subspace. If we identify the affine subspace with $\mathbb{R}^r$ then we can write $\Gamma$ as:

\[
\Gamma = \{y \in \mathbb{R}^r: y_1A_1+\dots+y_rA_r + A_{r+1} \succeq 0 \}
\]
where $A_1,\dots,A_r, \ A_{r+1}$ are symmetric $n\times n $ matrices. In general, the map $\mathcal{A}: \mathbb{R}^{r} \rightarrow \mathbb{S}^n$ given by $\mathcal{A}(y) = y_1A_1+\dots+y_rA_r + A_{r+1} $ is called an affine symmetric matrix map. Through duality, one can see that spectrahedrons are to semidefinite programs what polyhedra are to linear programs \cite{vinzant2014spectrahedron}. 
It is clear that whenever $\Gamma$ is polytope and we have an explicit representation of it given by a system of linear equations $Ax\leq d$, then program \ref{DSDP} reduces to a linear program. More interestingly perhaps is that the primal problem \ref{SDP} can also be solved as a linear program, albeit on potentially an exponential number of constraints.

\begin{theorem}\label{theorem:poly_exact}
Consider a generic semidefinite optimization problem \ref{SDP}, with dual given by \ref{DSDP}. Suppose that the set 
\[\Gamma = C-  \sum_{i=1}^ry_iA_i  \succeq  \ 0\]
is a polytope with extreme points $p_1\dots,p_k$, and define  $\mathcal{S}:=\bigcup_{i=1}^k \mathcal{E}(C- \mathcal{A}(p_k)) $. Then, $L_{\mathcal{S}}$ is a linear program and solves
\ref{SDP}.
\end{theorem}
\begin{proof}

The maximum value of the function $b^\top y$ over $\Gamma$ is achieved at some vertex $p$ of $\Gamma$. By strong duality and the solvability of \ref{DSDP}, there exists some $X^* \succeq 0 $ which solves program \ref{SDP}. In particular, program
\ref{LSDP} with $\mathcal{S}:=\bigcup_{i=1}^k \mathcal{E}(C- \mathcal{A}(p_k)) $ where $p_1,\dots,p_k$ are the vertices of $\Gamma$ is feasible. Let $\Hat{X}$ be an optimal solution to this program. Let $\{v_1\dots v_n\} \subseteq \mathcal{S}$ be an orthonormal eigenbasis for the matrix $C-\mathcal{A}(p)$. Since this matrix is positive semidefinite, we can write $C-\mathcal{A}(p) = \sum_{i=1}^n \beta_i v_iv_i^\top$ where the $\beta_i, i\in [n]$ are the (non-negative) eigenvalues of $C-\mathcal{A}(p)$. By feasibility of $\Hat{X}$,
$v^\top_i\Hat{X}v_i \geq 0$ for all $i \in [n]$. Multiplying each term by $\beta_i\geq 0$ we derive

\[
\sum_{i=1}^n \beta_i \left \langle \Hat{X}, v_iv_i^T \right \rangle =  \left \langle \Hat{X}, \sum_{i=1}^n \beta_i v_iv_i^T \right \rangle = \left \langle \Hat{X},C - \mathcal{A}(p)\right \rangle \geq 0.
\]

This implies that $\left \langle \Hat{X}, C \right \rangle \geq \left \langle \Hat{X}, \mathcal{A}(p) \right \rangle$.

To conclude, recall that for $j \in [r]$,
$\left \langle X, A_j\right \rangle = b_j$ giving the inequality $\left \langle \Hat{X}, C \right \rangle \geq b^\top p $. Again by strong duality and since the LP is a relaxation of the SDP, we have $b^\top p = \left \langle C, X^* \right \rangle \geq \left \langle C, \Hat{X} \right \rangle$  yielding the desired equality $  b^\top p = \left \langle C, X^* \right \rangle $.

\end{proof}

In \cite{ramana1997polyhedra} is its shown that deciding if a spectrahedron is a polyhedron is in co-NP, and an algorithm for deciding polyhedrality is given.  \cite{bhardwaj2015deciding} generalizes and improves the previous results. The algorithm presented in the latter paper runs in exponential time, as it requires enumerating the vertices of a certain polyhedron. Even if we knew that $\Gamma$ is polyhedral, we do not have exact solvability under $\mathcal{O}$, as the previous problem has an exponential number of constraints. A particular case in which $\Gamma$ is polyhedral and that has received attention in the literature is whenever the matrices $C$ and $A_i, \  i\in [r]$ are simultaneously diagonalizable. This is the content of observation \ref{obs:simulDiag}, which we now prove.

\begin{proof}[Proof of Observation \ref{obs:simulDiag}. Also see \cite{wang2022tightness}, Lemma $9$ ]

Let $U$ be a matrix that simultaneously diagonalizes matrices $C$ and $A_i, \  i \in [r]$ i.e. the matrices 
$C' = U^\top C U$ and $A'_i = U^\top A_i U $ are all diagonal. By Silvester's law of inertia \cite{horn2012matrix}, we have that 
$C -\mathcal{A}(y) \succeq 0 $  if and only if
$U^\top \left [ C-\mathcal{A}(y) \right ]U \succeq 0 $ if and only if 
$C'-\sum_{i=1}^r y_iA'_i \succeq 0$. Hence, we have

\[
\Gamma = \{y \in \mathbb{R}^r : C' -\sum_{i=1}^r y_iA'_i \succeq 0 \}
\]
which is a polyhedral set since all matrices involved are diagonal.
\end{proof}

For a clear exposition of the implications of this observation to QCQPs see \cite{wang2022tightness} and the references therein. In addition, the authors show that the region $\Gamma$ might be polyhedral even if the matrices $C$ and $\{A_i\}_{\ i\in [r]}$ are not simultaneously diagonalizable. Although the latter condition is much more stringent, it allows us to avoid the need to have the vertices of $\Gamma$ given to us explicitly, as Theorem \ref{theorem:poly_exact} requires.

\begin{theorem}\label{teo:diag_exactness}
Let \ref{SDP} be a semidefinite program with dual \ref{DSDP}. Suppose that the set of matrices $\{C, A_1,\dots,A_r\}$ is simultaneously diagonalizable. Then, \ref{SDP} is solvable under $\mathcal{O}$. 
\end{theorem}
\begin{proof}

Let $U$ be a orthogonal matrix that simultaneously diagonalizes $C$ and $A_i$ for each $i \in [r]$. Let $v_1,\dots. v_n$ denote the columns of $U$ and set $\mathcal{S} = \{v_1,\dots, v_n\}$. 
Let $p^*$ be a dual optimal solution with $C-\mathcal{A}(p^*) = S^* $ where $S^*$ is positive semidefinite. Since $U$ diagonalizes
each $A_i,\ i\in [r]$ and $C$, it is clear that the matrix $U^\top \left [ C-\mathcal{A}(p) \right ]U$ is diagonal. In other words, the matrix
$U^\top S^* U = D$ for some diagonal matrix $D$ with non-negative entries. This means that we can express $S^*$ as 

\[
S^* = \sum_{i=1}^n \beta^*_iv_iv_i^\top, \beta_i^* \in \mathbb{R}_+ \ \forall \ i \in [n].
\]

We turn our attention the linear relaxation of \ref{SDP} defined by $\mathcal{S}$, defined in Section \ref{revised_sec:1}, which we recall is given by

\begin{equation}\tag{$L_\mathcal{S}$}
\begin{aligned}
\min_{X\in \mathbb{S}^{n}} \  &\langle  C,X\rangle \\
\text{ s.t: } \langle A_i,X \rangle  = & \  b_i,\  \forall i \in [r],\\ 
v^\top Xv  \geq & \ 0 \ \forall \ v \in \mathcal{S}.
\end{aligned}
\end{equation}

This program is linear and is a relaxation of \ref{SDP} as any feasible solution to it is feasible for \ref{LSDP}.
Its dual is given by

\begin{equation}\label{dSrel}\tag{$DL_\mathcal{S}$}
\begin{aligned}
&\max_{y,\in \mathbb{R}^{n},\beta \in \mathbb{R}^{n}_+} \  b^\top y \\
\text{ s.t: }& C-  \sum_{i=1}^ry_iA_i  = \sum_{i=1}^n \beta_i v_iv_i^\top.
\end{aligned}
\end{equation}

Observe that this program is a strenghening of program \ref{DSDP}, and that $S^*$ is feasible for this program. Therefore, their optimal values must match, and in particular the optimal value of \ref{dSrel} is finite. By strong duality of linear programs, \ref{LSDP} is solvable and its optimal value equals the optimal value of both \ref{DSDP} and \ref{dSrel}. Again by our strong duality assumption of programs \ref{SDP} and \ref{DSDP}, program \ref{LSDP} solves \ref{SDP}.

    
\end{proof}


A class of problems that has been extensively studied in the literature and where the hypothesis of our previous theorem applies are simultaneously-diagonalizable QCQPs. Recall that a QCQP is a problem of the form 

\begin{equation}\label{QCQP}\tag{QCQP}
 \inf_{x \in \mathbb{R}^n} q_0(x) : q_i(x) \leq \  0 \  \forall  \ i \in [r].
\end{equation}

where $q_i(x) = x^\top A_i x + 2b_i^\top x + c_i $ with $A_i \in \mathbb{S}^n$, $b \in \mathbb{R}^n$ and $c_i \in \mathbb{R}$ for all $i \in \{0,\dots,r\}$. QCQPs are NP-hard to solve in general but admit tractable convex relaxations. The SDP relaxation of a QCQP is given by the following semidefinite program \cite{bao2011semidefinite,shor1990dual}:

\begin{equation}\label{shor}
\begin{aligned}
\inf_{x \in \mathbb{R}^n, X\in \mathbb{S}^n} & \left \langle A_0,X \right \rangle + 2b_0^\top x + c_0 \\
 \text{s.t }:   \left \langle A_i,X \right \rangle + & 2b_i^\top x + c_i   \leq  \ 0 \ \forall i \in [r] \\
&  
\begin{bmatrix}
X & x \\
x^\top & 1 
\end{bmatrix} \succeq 0.
\end{aligned}     
\end{equation}

Whenever the $A_i$ are simultaneously diagonalizable and we have access to a matrix $U$ such that $A_i = U^\top D_i U$ for $i \in \{0,\dots,r\}$, we can perform the change of variables $y = Ux$ and $ \Tilde{b_i} = Ub_i, \ i \in  \{0,\dots,r\} $ to obtain the a diagonalized version of the problem 

\begin{equation}\label{DQCQP} 
\inf_{y \in \mathbb{R}^n} q_0(y) : q_i(y) \leq \  0 \  \forall \  i \in [r]
\end{equation}

However, we have $q_i(y) = a_i^\top y^2 + 2\Tilde{b}_i^\top y + c_i  $, $d_i \in \mathbb{R}^n$, $\Tilde{b}_i \in \mathbb{R}^n$ and $c_i \in \mathbb{R}$ for each $i \in [0,\dots,r]$. Here, $y^2 \in \mathbb{R}^n$ is the vector whose entries are the squared entries of the vector $y \in \mathbb{R}^n$.
Ben-Tal and den Hertog \cite{ben2014hidden} and Locatelli \cite{locatelli2016exactness} study a certain second order cone relaxation of this problem, and show that the optimal value of that relaxation and that of the SDP relaxation match. Our results imply that in fact, given access to a matrix $U$ that simultaneously diagonalizes the $A_i, \ i \in [0,\dots, r]$ we can solve the SDP relaxation (\ref{shor}) with the linear program \ref{LSDP} where $\mathcal{S}$ is the set of columns of $U$.

\begin{corollary}
Consider a quadratically constrained quadratic problem given as in \ref{QCQP} and such that the matrices $\{A_i\}_{i \in \{0,\dots,r\}}$ are simultaneously diagonalizable by an orthogonal matrix $U$. Let $opt$ be the optimal value of relaxation (\ref{shor}) of \ref{QCQP}. Let $\mathcal{S}$ be the set of columns of $U$. Then, the objective value $z$ of the linear relaxation \ref{Srel} of (\ref{shor}) equals $opt$.   
\end{corollary}
\begin{proof}
The proof is immediate from Theorem \ref{teo:diag_exactness}.
\end{proof}

\subsection{Finding initial sets}

As we have seen in Theorem \ref{teo:diag_exactness}, we know some vectors whose inclusion in $\mathcal{S}$ guarantees solvability under $\mathcal{O}$. The reason this worked was that we were able to produce a feasible solution to \ref{dSrel} which matches the objective of an \textit{optimal} solution to \ref{DSDP}. Nonetheless, the previous argument still holds for a generic feasible solution to \ref{DSDP}: any dual feasible solution will generate sets $\mathcal{S}$ that satisfy the corresponding dual bound.

\begin{lemma}

Consider a generic SDP problem and let $\Hat{y}$ be a feasible solution to the dual of the SDP with objective value $b^\top \Hat{y}$. Let $\mathcal{S} = \mathcal{E}(C-\mathcal{A}(\Hat{y}))$. Then, the objective value $z^*$ of 
program \ref{LSDP} satisfies 

\[
z^* \geq b^\top \Hat{y}.
\]
    
\end{lemma}
The proof of this lemma is very similar to that of Theorem \ref{teo:diag_exactness}. This result indicates that finding a good set $\mathcal{S}$ amounts to finding feasible solutions to the dual of \ref{SDP} whose objective value is close to optimal. This task is akin to finding good feasible solutions to \ref{SDP}, or at worse to solve a semidefinite feasibility problem, which in principle may be as hard as solving the original problem. However, the results of the previous subsection suggest a way to get around this issue by exploiting simultaneous diagonalizability. Under a weakening of this assumption, we will be able to construct solutions, which will be automatically feasible for the the \ref{DSDP}.
To begin, we give in Proposition \ref{prop:commut_generalized} a characterization of simultaneous diagonalizability which we will then relax.

\begin{lemma}\label{lemm:commute}
    Let $\{A_i\}_{i \in I} \subseteq \mathbb{S}^n $ be a set of symmetric matrices. Then, there exists a basis of orthonormal vectors $\{u_1,\dots,u_n\}$ that simultaneously diagonalizes $\{A_i\}_{i \in I}$ if and only if $A_i$ and $A_j$ commute for every $i$ and $j \in I$, i.e, $A_iA_j = A_jA_i \ \forall i,j \in I$.
\end{lemma}
See \cite{conradsimultaneous} for a proof.

\begin{proposition}\label{prop:commut_generalized}
The set of matrices $\{A_1,\dots A_r\} \subseteq \mathbb{S}^n$ is simultaneously diagonalizable if and only if for every $p$ and $q \in \mathbb{R}^r$ the matrices 
$\mathcal{A}(p) = \sum_{i=1}^r p_iA_i $ and $\mathcal{A}(q) = \sum_{i=1}^r q_iA_i $ commute, and hence are simultaneously diagonalizable.
\end{proposition}
\begin{proof}
    Necessity is trivial by having  $p$ and $q$ range over the standard basis of $\mathbb{R}^r$ and Lemma \ref{lemm:commute}. For sufficiency, Let $U$ be an orthonormal matrix such that the matrices $U^\top A_i U = D_i $ are diagonal $ \forall i \in [r]$. Given $p$ and $q \in \mathbb{R}^r$ we have:
\[
U^\top \mathcal{A}(p) U = U^\top \left ( \sum_{i=1}^r p_iA_i \right ) U = \sum_{i=1}^r p_i D_i.
\]
Similarly we have $U^\top \mathcal{A}(q) U = \sum_{i=1}^m q_i D_i $. 
Since diagonal matrices commute we have 

\[
\begin{aligned}
\left (\sum_{i=1}^r p_iD_i\right) \left(\sum_{i=1}^r q_iD_i \right ) = \left (\sum_{i=1}^r q_iD_i\right) \left(\sum_{i=1}^r p_iD_i \right ).
\end{aligned}
\]
Given that $U^\top U = I$, pre-and post-multiplying by $U$ and $ U^\top$ respectively gives:

\[
\begin{aligned}
U\left (\sum_{i=1}^r p_iD_i\right)U^\top U \left(\sum_{i=1}^r q_iD_i \right ) U^\top = U \left (\sum_{i=1}^r q_iD_i\right) U^\top U \left(\sum_{i=1}^r p_iD_i \right ) U^\top
\end{aligned}
\]
and finally

\[
\mathcal{A}(p)\mathcal{A}(q) = \mathcal{A}(q)\mathcal{A}(p). 
\]

Since these matrices commute, they are simultaneously diagonalizable.

\end{proof}

Given that commutativity of the set $\{C,A_1,\dots,A_r\}$ will typically not hold, we relax the equivalent condition given by the previous lemma to require that commutativity holds only for special class of $p$'s and $q$'s. In particular we will set $p=e_{r+1}$ and
$q$ such that for some subset $J\subseteq [r]$ we have 
$\sum_{j \in J} q_jA_j = I_n$. The idea is that if we have a point $y \in \mathbb{R}^n$, not necessarily dual feasible for which the matrices $C$ and $\mathcal{A}(y)$ commute, then taking $\mathcal{S}$ to be the columns of a matrix that diagonalizes them will yield a linear program with objective value as good as the best dual feasible solution that lies on the set 

\[\{A\in \mathbb{S}^n : \exists \ x,t \in \mathbb{R}: A = tI_n + x\sum_{j \in [r]\setminus J}q_jA_j  \}.\]

\begin{theorem}\label{teo:dual_bound}
Consider a generic semidefinite optimization problem of the form \ref{SDP}, with dual \ref{DSDP}. Suppose that there exists vectors $q^1,q^2 \in \mathbb{R}^{r}$ whose support is disjoint such that $\sum_{j =1}^rq^1_jA_j = I_n$ and
such that the matrices $C$ and $\sum_{j=1 }^r q^2_jA_j $ commute and therefore are simultaneously diagonalizable by some orthogonal matrix $U$. Let $\mathcal{S}$ to be the set of columns of such an $U$. Then, the optimal value $z$ of program $L_{\mathcal{S}}$
satisfies the bound

\[
  \left (\sum_{j =1}^r b_jq^2_j \right )x + \left(\sum_{j=1}^r b_jq^1_j \right )t \leq z  
\]

for any $x$ and $t$ such that the matrix $C-x \left (\sum_{j=1}^r q_jA_j  \right ) +tI_n$ is positive semidefinite.
    
\end{theorem}

\begin{proof}

Let $U$ be a matrix that simultaneously diagonalizes $C$ and $\mathcal{A}(q^2) =\sum_{j=1}^r q^2_jA_j$. Let $z$ be the optimal value
of program $L_\mathcal{S}$ where $\mathcal{S}$ is the set of columns $v_1,\dots,v_n$ of $U$. Recall that the dual of this program is given by

\begin{equation} \tag{$DL_{\mathcal{S}}$}
\begin{aligned}
&\max_{y,\in \mathbb{R}^{n},\beta \in \mathbb{R}^{n}_+} \  b^\top y \\
\text{ s.t: }& C-  \sum_{i=1}^ry_iA_i  = \sum_{i=1}^n \beta_i v_iv_i^\top.
\end{aligned}
\end{equation}

Since $U$ diagonalizes $C$, any column $v$ of $U$ is an eigenvector of $C$ with some corresponding eigenvalue $\lambda$, and the same holds for $\mathcal{A}(q)$ with some eigenvalue $\gamma$.
Hence, $v$ is a eigenvector of $C-x\mathcal{A}(q)+tI_n$ with corresponding eigenvalue $\lambda- x\gamma + t $. Since we are looking for $x$ and $t$ values such that $C-x\mathcal{A}(q)+tI_n$ is psd, this gives rise to the equation $\lambda- x\gamma + t \geq 0$, and we have such one equation for every column of $U$. This system is always feasible as the $t$ variable is free. Hence, there exists $x^*, t^*$ for which the matrix $C-x^*\mathcal{A}(q)+t^*I_n$ is positive semidefinite. As $U$ diagonalizes $C$, $\mathcal{A}(q)$ and $I_n$ as $U^\top I_n U = U^\top U = I_n $, $C-x^*\mathcal{A}(q)+t^*I_n$ is diagonalizable by $U$ and thus can be written as $\sum_i\eta_iv_iv_i^\top$ with $\eta_i \geq 0$
for $i \in [n]$. Thus, setting $y_j = x^*q_j^2 $ if $j$ belongs to the support of $q^2$ and $y_j = t^*q_j^1$ if $j$ belongs to the support of $q^1$ (here recall that $q^1$ and $q^2$ have disjoint support) gives a feasible solution to program $DL_\mathcal{S}$ by setting $\eta_i = \beta_i$ for $i \in [n]$.
The objective value of this solution is
\begin{equation}\label{eq:teo3}
\left (\sum_{j=1}^r b_jq^2_j \right )x^* + \left( \sum_{j=1}^rb_jq^1_j \right ) t^*.
\end{equation}
\end{proof}

We make a few observations about this theorem. First and foremost, we didn't require  that 
the matrix $I_n + \sum_{j=1}^r q_jA_j $ is feasible for program \ref{DSDP}. Second, notice that we have required that we can aggregate some of the $A_j$ to form the identity matrix. Although this seems quite constraining, it is always the case that such a combination exists by our assumption that $\ref{DSDP}$ is strictly feasible, i.e if there exists $ q  \in \mathbb{R}^r$ such that $C-\mathcal{A}(q) \succ 0$. In principle, finding such $q$ would require finding finding a point in the interior of the dual feasible region, which might be non-trivial. This suggests that our theorem is easier to apply in regimes where it is more directly ``obvious " which combination of the $A_j$ forms the identity. This is the case in the max cut problem, the Lov\'asz theta number, the sparse PCA problem, the extended trust region SDP relaxation and many others. Finally, we observe that even though the bound given in Equation \ref{eq:teo3} is the best bound we can \textit{prove}, there might be other ``hidden " dual feasible solutions that certify a better bound for $L_\mathcal{S}$.

\begin{observation}[Hidden basis property]\label{hidden}
    Let $\Hat{y}$ be a dual feasible solution for program \ref{DSDP} with objective value $b^\top\Hat{y}$. Suppose that $y \in \mathbb{R}^r$ is a point such that the matrices $C-\mathcal{A}(\Hat{y})$ and $C-\mathcal{A}(y)$ (which is not necessarily PSD) share a basis of orthonormal eigenvectors. Let $\mathcal{S} = \mathcal{E}(C-\mathcal{A}(y))$ then, the optimal value $z$ of program \ref{LSDP} satisfies

    \[
     b^\top \Hat{y} \leq z.
    \]
    
\end{observation}

The proof of this observation is straightforward, but note that we have required $\mathcal{S}$ to be some eigenbasis of $C-\mathcal{A}(y)$ rather than the set of columns of some orthogonal matrix that simultaneously diagonalizes $C$ and $\mathcal{A}(y)$. Clearly, if $U$ diagonalizes both of those matrices it diagonalizes any linear combination of them. As we will see in the max cut experiments, Theorem \ref{teo:dual_bound} will certify a spectral bound, but the LP relaxation will actually have a better objective than the bound of Theorem \ref{teo:dual_bound} guarantees in practice.

\subsection{Finding commuting matrices}\label{subsec:finding_conmuting}

To apply Theorem \ref{teo:dual_bound}, we need first to find a combination of the constraints matrices which commutes with the
objective matrix $C$ of $SDP$. This can be accomplished using a linear program. Picking $L$ to be an arbitrary linear function
on $y$ gives the program

\begin{equation}
\begin{aligned}
&\min_{y\in \mathbb{R}^{n}} \  L(y) \\
\text{ s.t: }& C \mathcal{A}(y) =  \mathcal{A}(y) C.  
\end{aligned}
\end{equation}

To select $L$, we propose a function that trades off between the $\ell_1$ norm of the matrix $C-\mathcal{A}(y)$ and the dual objective function
$b^\top y$. The intention of the $\ell_1$ term is to promote solutions where $C-\mathcal{A}(y)$ is sparse, rendering the computation of an eigenbasis easier. The term $-b^\top y$ encourages having solutions with good dual objective value. This yields the program

\begin{equation}\label{conmute_generator}\tag{CG}
\begin{aligned}
&\min_{y\in \mathbb{R}^{n}} \ \sum_{i,j} \left|[C-\mathcal{A}(y)]_{ij}\right| - b^\top y   \\
\text{ s.t: }& C \mathcal{A}(y) =  \mathcal{A}(y) C.  
\end{aligned}
\end{equation}
 
Note that the null vector is always a feasible solution to this program.
In Section \ref{revised_sec:5} we experimentally test this idea.

\section{Linear Relaxations of the max cut semidefinite program}
\label{revised_sec:3}

The question of exactly - or approximately - solving an SDP with a linear program finds one of its historical roots in the max cut problem, where in a given undirected graph, we seek a bipartition of the nodes to maximize the number of edges with one end in both parts. Since linear programming has been one of the main paradigms to tackle NP-hard combinatorial optimization problems through the relax-and-round paradigm, substantial efforts were dedicated to find a linear programming relaxation of the max cut problem. A graph with $m$ edges has always a cut of size at least $\frac{1}{2}m$ and any cut can cut at most $m$ edges, so it is trivial to provide an algorithm with integrality gap \footnote{In this paper, we employ the convention that the integrality gap is a number that is at least 1 and hence is the ratio of the value of the relaxation to the optimal value of the max cut.} $2$. For example, a randomized algorithm picking vertices at random or a greedy algorithm will have this guarantee. The question was then if there exists a linear program that could have an approximation ratio better than $2$.

The starting point of this line of research was perhaps the linear relaxation for max cut given by \cite{barahona1986cut,poljak1994expected}. Let $G=(V,E)$ be an undirected, simple graph on $m$ edges and $W$ its adjacency matrix. We define
\[
\alpha(G):= \max \langle W ,X\rangle
\]
\begin{equation}\label{usualLP}
\begin{aligned}
X_{ij}+X_{ik}+X_{kj} \leq 2 \  \forall i,j,k \in V \\
X_{ij} - X_{ik} - X_{jk} \leq 0  \  \forall i,j,k \in V \\
0 \leq X_{ij} \leq 1 \ \forall i,j \in V.
\end{aligned}
\end{equation}

Here we use a binary variable $X_{ij}$ for each pair of vertices $\{i,j\}$ to denote if the edge between them is cut. The first set of `triangle' constraints specify that at most two edges can be picked in a cut from any triangle, while the second set rules out exactly one edge from any triangle from being selected in a cut. 
In ~\cite{poljak1994expected}, Poljak and Tuza prove that for sparse and dense versions of Erd\H{o}s-R\'enyi random graphs, the integrality gaps of this LP tend to $2-o(1)$ and $\frac43-o(1)$ respectively. Here, $G_{n,p}$ denotes the class of random graphs on $n$ nodes where every edge is included independently of others with probability $p$.

\begin{theorem}\label{theorem:hardnessLP} (Poljak, Tuza)~\cite{poljak1994expected}
Let $mc(G)$ denote the size of the max cut of $G$.
\begin{itemize}
    \item (Sparse graphs). Let $p(n)$ be a function such that $0<p<1$, $p(n)\cdot n \rightarrow \infty $ and $p\cdot n^{1-a}\rightarrow 0$ for every 
$a>0$, then the expected relative error $\frac{\alpha(G_{n,p})-mc(G_{n,p})}{mc(G_{n,p})}$  tends to $1$ as $n\rightarrow \infty$ with probability $1-o(1)$.
\item   (Dense graphs). Let $p(n)$ be a function such that $0<p<1$, $p(n) = \Omega\left(\sqrt{\frac{\log(n)}{n}}\right)  $. Then the expected relative error $\frac{\alpha(G_{n,p})-mc(G_{n,p})}{mc(G_{n,p})}$, tends to $\frac{1}{3}$ as $n\rightarrow \infty$ with probability $1-o(1)$.
\end{itemize}
\end{theorem}

Such integrality gap lower bounds for the basic LP encouraged two distinct approaches to solve the problem. The first one focused on adding valid constraints to formulation (\ref{usualLP}), such as ''hypermetric", and ''gap" constraints. See \cite{deza2009geometry,o2018sherali} for more details. Nonetheless, a long line of research culminated in showing that such direct strengthenings will fail to provide an approximation factor better than $2$ \cite{chan2016approximate,charikar2009integrality,de2007linear}.
In particular,  Kothari et al. \cite{kothari2021approximating} prove that this problem - and more generally Constraint Satisfaction Problems - is resistant to this strategy by showing that extended linear formulations are as powerful as the Sherali-Adams hierarchy, which in turn requires an exponential number of rounds (in $\varepsilon$) to certify an integrality gap better than $2-\varepsilon$. 
The second approach, perhaps much more influential, considered stronger optimization relaxations, such as the vector optimization relaxation of 
Poljak and Rendel \cite{poljak1995nonpolyhedral}, shown to be SDP-representable and providing an approximation ratio of 
$\sim 1.13$ in the seminal work of Goemans and Williamson \cite{goemans1995improved}. Naturally, this leads to the question if linear programs can well approximate semidefinite ones. In  \cite{braun2015approximation} Braun et. al. show that in principle one needs an exponential number of a constraints in an LP to correctly approximate an SDP. 
These two combined results extinguish the hope that linear programming may be used to approximate max cut. Since the question of finding a good set $\mathcal{S}$ to initialize Algorithm \ref{alg1} amounts to finding a linear approximation to a semidefinite program, these results suggest that no systematic procedure can generate a good set $\mathcal{S}$ as in particular they would provide an approach to obtain a low-gap linear programming approximation to the max cut problem. In this sense, we propose to use instance-specific information to avoid the hardness of approximation results, in particular by exploiting bounds relating the spectrum of the graph to the value of the max cut, resulting in linear relaxations with better approximation ratios.

For a graph $G=(V,E)$ we set $m=|E|$ and denote by $W$ its adjacency matrix. Recall that the semidefinite relaxation for max cut due to of Poljak, Rendl, Goemans and Williamson \cite{goemans1995improved,poljak1995nonpolyhedral} is given by 

\begin{equation}\label{GW}\tag{GW}
\begin{aligned}
 &\frac{1}{2}m +  \frac{1}{4} \max_X\ \left \langle -W,X\right \rangle \\
\text{ s.t: } & X\succeq 0,\ \ X_{ii}=1, \ \forall \ i \in [n].
\end{aligned}
\end{equation}

with dual

\begin{equation}\label{DGW}\tag{DGW}
\begin{aligned}
 &\frac{1}{2}m + \frac{1}{4} \min_{\gamma \in \mathbb{R}^n} \ \sum_{i=1}^n \gamma_i \\
&\text{ s.t: }  W+diag(\gamma) \succeq 0. 
\end{aligned}
\end{equation}

It is known that strong duality holds for this pair of programs: both (\ref{GW}) and (\ref{DGW}) are solvable and their objectives coincide. Delorme and Poljak show \cite{delorme1993laplacian} that the max cut value of $G$ on $n$ nodes is upper bounded by the quantity 
\[\min_{u \in \mathbb{R}^n: \sum_i u_i=0 } \ \ \frac{n}{4}\lambda_{1}(\mathcal{L}(G)+diag(u)).\]

It turns out that this program is equivalent to program \ref{DGW} \cite{goemans1995improved}. In their seminal work, Goemans and Williamson show that this program achieves an approximation ratio of roughly $\frac{1}{0.878}\sim 1.138$. Through this equivalence, one can show that the semidefinite program \ref{GW} satisfies a series of eigenvalue bounds. For instance, one may take $u$ such that $\sum_{i=1}^n u_i=0$ and all of the diagonal entries of the matrix $\mathcal{L}(G)+diag(u)$ equal $\frac{2m}{n}$. This results in what is usually known as \textit{the} eigenvalue bound for max cut due to Mohar and Poljak \cite{mohar1990eigenvalues} 

\begin{equation}\label{eigenBound}
    mc(G) \leq \frac{1}{2}m+\frac{n}{4}\lambda_{1}(-W) = \frac{1}{2}m-\frac{n}{4} \lambda_n(W)
\end{equation}

To see the second equality, recall that for any matrix $A$, $\lambda_n = \lambda_1(-A)$). See~\cite{alon2000bipartite,mohar1990eigenvalues} for an elementary proofs of this inequality. As mentioned in \cite{o2018sherali}, conventional wisdom is that LPs cannot certify even the eigenvalue bound, and we are not aware of a polynomially sized linear program that certifies this bound.
 \comment{

This bound is powerful enough to approximate the max cut value of the graphs of Theorem \ref{thm:pt}, as the next observation shows.

\begin{observation}\label{teoErdosRenyi}
Let $G$ be generated from the Erd\H{o}s-R\'enyi random graph model
$\mathcal{G}_{n,p}$ where $p \equiv p(n)$ is a function that satisfies $ 0< p < 1$ and $W$ its adjacency matrix. Let $\chi(G) \in \mathbb{R}^+$ be the following quantity:
\[\chi(G) = \frac{m}{2}-\frac{n}{4}\lambda_{max}(-W).\]

\begin{itemize}
    \item (Sparse graphs). Let $p(n)$ be a function such that $0<p<1$, $p(n)\cdot n \rightarrow \infty $ and $p\cdot n^{1-a}\rightarrow 0$ for every 
$a>0$, then the relative error $\frac{\chi(G_{n,p})-mc(G_{n,p})}{mc(G_{n,p})}$  tends to $0$ as $n\rightarrow \infty$ with probability $1-o(1)$.
\item   (Dense graphs). Let $p(n)$ be a function such that $0<p<1$, $p(n) = \Omega\left(\sqrt{\frac{\log(n)}{n}}\right)  $. Then the expected relative error $\frac{\chi(G_{n,p})-mc(G_{n,p})}{mc(G_{n,p})}$, tends to $0$ as $n\rightarrow \infty$ with probability $1-o(1)$.
\end{itemize}
\end{observation}
}

\subsection{Instance-specific linear relaxations.} 

The specialization of program \ref{LSDP} to the max cut problem results in a 
polynomially sized linear program that explicitly depends of the adjacency matrix $W$ on $G$, allowing us to circumvent the theoretical limitations of linear relaxations described in the introduction of this section. Using Theorem \ref{teo:dual_bound} this LP will be shown to satisfy the eigenvalue bound (\ref{eigenBound}) whenever $\mathcal{S}$ is chosen appropriately. Fixing $\mathcal{S} = \{v_1,\dots,v_k\}$, program \ref{LSDP} specializes to a linear program which we denote by program \ref{Srel}.

\begin{equation}\label{Srel}\tag{$SP_{\mathcal{S}}$}
\begin{aligned}
&\max_{X \in \mathbb{S}^n} \ \frac{1}{2}m + \frac14 \langle -W,X\rangle   \\
\text{ s.t: } v^\top Xv  \geq  0 \ & \forall \ v \in \mathcal{S},  X_{ii}=1, \ \forall \ i \in [n],
  \| X \|_\infty \leq 1.
\end{aligned}
\end{equation}

In this program we have included the constraint 
$  \| X \|_\infty \leq 1$. As the following observation shows, this is a valid constraint for \ref{GW}. Adding it is useful because it will guarantee that the dual of \ref{Srel} is always feasible, regardless of $G$.

\begin{observation}
Let $X$ be feasible for program (\ref{GW}). Then, it is feasible for program $SP_{\mathcal{S}}$ for any set $\mathcal{S}\subseteq \mathbb{R}^n$.
\end{observation}
\begin{proof}
Let $X$ be feasible for (\ref{GW}). This means $X$ is positive semidefinite, and that there exists a set of vectors $x_1,\dots,x_n$ such that $X_{ij}= x_i^\top x_j $ for all $i,j\in [n]$. For each $i\in [n] $ we have $X_{ii}=1$ and thus we see that $\|x_i\|_2=1$. It follows that each entry of the vectors $x_i$ is bounded by $1$ and therefore that $X_{ij}$ is bounded by $1$ for all $i$ and $j$. The other two constraints of the linear program are clearly satisfied by $X$.  

\end{proof}

It will be also be useful to consider the following strenghening of program \ref{GW} depending of $\mathcal{S} = \{v_1,\dots,v_k\}$. 

\begin{equation}\label{DSDGW}\tag{$SD_{\mathcal{S}}$}
\begin{aligned}
\frac{1}{2}m + \frac{1}{4}& \max_{\eta \in \mathbb{R}^k} \left \langle -W,\sum_{i=1}^k \eta_iv_iv_i^\top \right \rangle \\
 \text{s.t: } diag \left(\sum_{i=1}^k \eta_iv_iv_i^T \right) \leq  1, & \ \eta_i \geq 0, \ v_i \in \mathcal{S} \ \forall i \in [k], \ k=|\mathcal{S}|.
    \end{aligned}
\end{equation}

Here, and for the rest of the paper, we denote by $z_{SP_{\mathcal{S}}},z_{GW}, z_{DGW},z_{SD_{\mathcal{S}}}$ the optimal values of \ref{Srel}, \ref{GW}, \ref{DGW} and \ref{DSDGW} ignoring the additive constant $\frac{1}{2}m$ and the multiplicative constant $\frac{1}{4}$, respectively.
For illustration, we have:

\[
\begin{aligned}
 &z_{GW} = \max \ \left \langle -W,X\right \rangle \\
\text{ s.t: } & X\succeq 0,\ \ X_{ii}=1, \ \forall \ i \in [n].
\end{aligned}
\]

By duality, we get the following relationships between these optimal values 
\[
z_{SD_{\mathcal{S}}} \leq  z_{GW} = z_{DGW} \leq z_{SP_{\mathcal{S}}}.
\]

Observe that we may employ different sets $\mathcal{S}$ to define $SP$ and $SD$ and the above relations will continue to hold. As a sanity check, we first observe that program \ref{Srel} satisfies the trivial bound for max cut.

\begin{lemma}
Let $\mathcal{S}$ be an arbitrary subset of $\mathbb{R}^n$. Then $z_{SP_{\mathcal{S}}}$ satisfies:

\[
z_{SP_{\mathcal{S}}} \leq 2m
\]
and therefore $\frac{1}{2}m+\frac{1}{4}z_{SP_{\mathcal{S}}} \leq m$.
\end{lemma}
\begin{proof}
Let $\mathcal{S} = \{v_1,\dots,v_k\}$. The dual of program \ref{Srel} is given by

\begin{equation}\label{MCdSrel}\tag{$DSP_{\mathcal{S}}$}
\begin{aligned}
\min_{\lambda, \alpha,\delta, \beta,\Lambda} \ &\frac{1}{2}m - \frac14 \left  [ tr(\Lambda) - \sum_{i\neq j}\delta_{ij} - \sum_{i\neq j}\alpha_{ij}   \right]    \\
\text{ s.t: } & W - \Lambda = \sum_{i=1}^k \beta_iv_iv_i^\top, \\
 & \delta_{ij} \geq 0 \ \forall i\neq j \in [n], \\
  & \alpha_{ij} \geq 0 \ \forall i\neq j \in [n], \\
  & \lambda_{i} \in \mathbb{R} \ \forall i \in [n], \\
   & \beta_{i} \geq 0 \ \forall i \in [n], \\
    \Lambda \in \mathbb{S}^n, \Lambda_{ij} = & \delta_{ij}-\alpha_{ij}\ \forall i\neq j \in [n],
   \Lambda_{ii} = \lambda_{i} \  \forall i \in [n].
\end{aligned}
\end{equation}

The proof of this fact is deferred to Appendix \ref{AppendixA}. Letting $\beta_i = 0 \ \forall i \in [n]$, $\Lambda = W$ where $\delta_{ij} = 1, \alpha_{ij}=0$ whenever $W_{ij}=1$ and $0$ otherwise, we obtain a feasible solution for the previous program with 
$tr(\Lambda)-\sum_{i\neq j}\delta_{ij}-\sum_{i\neq j}\alpha_{ij}=-2m$. 

\end{proof}

It can be checked that for an arbitrary graph $G$, the feasible region of program \ref{DGW}, namely $\Gamma = \{\gamma \in \mathbb{R}^n: W+ diag(\gamma) \succeq 0 \}$
is not necessarily polyhedral. However, we can exploit Theorem \ref{teo:dual_bound} to derive a set $\mathcal{S}$ for the relaxation \ref{Srel} that has a good objective value. Although this statement can be proven directly by simply giving a judicious choice of $\mathcal{S}$, we derive the result in a way that explicitly uses the theorem.

\begin{theorem}\label{teo:eigenbound}
Let $G$ be a graph on $n$ vertices and $W$ its adjacency matrix. Let $\lambda_n$ denote the smallest eigenvalue of $W$. Set $\mathcal{S} = \mathcal{E}(W) $. Then
\[
z_{SP_{\mathcal{S}}} \leq  -n\lambda_n := \chi(G).
\]
\end{theorem}
\begin{proof}

For $i=1,\dots,n$ let the matrix $A_i$ denote the matrix of all zeros but with a single $1$ in its $i$-th diagonal entry. Hence, $\Gamma$ can be expressed as:

\[
\Gamma = \{\gamma \in \mathbb{R}^n: W+ \sum_{i}^n \gamma_iA_i \succeq 0 \}.
\]

To apply Theorem \ref{teo:dual_bound}, we express the identity as some combination of the $A_i$. concretely, we 
let $\Hat{\gamma} = \Vec{1} $ be the vector of all ones in $\mathbb{R}^n$ so that we have that 
$\sum_{i=1}^n \Hat{\gamma}_iA_i = I_n$. By Theorem \ref{teo:dual_bound},
it follows that if $\mathcal{S} = \mathcal{E}(W)$ then the optimal value $z_{SP_{\mathcal{S}}}$ of program \ref{Srel} satisfies

\[
z_{SP_{\mathcal{S}}} \leq  t\cdot n 
\]
for any $t$ such that $W +tI$ is positive semidefinite. Observe that  $ W -\lambda_nI $ is positive semidefinite. In particular, we obtain 

\[
z_{SP_{\mathcal{S}}} \leq  -n\lambda_n.
\]

\end{proof}

We provide an alternate direct proof of this result in Appendix \ref{AppendixA} by directly setting $\mathcal{S}=\mathcal{E}(W)$ and using the dual of program \ref{Srel}. Interestingly, this result allows us to show that the linear relaxation \ref{Srel} is strictly stronger than the linear formulation for max cut given in program (\ref{usualLP}), in the sense that it gives -in contrast to the previous LP- the correct value of max cut for the graphs considered in Theorem \ref{theorem:hardnessLP}. Perhaps more interestingly, we show that for random $d-$regular graphs the linear program \ref{Srel} with $\mathcal{S}=\mathcal{E}(W)$ approximates max cut with an approximation factor of $1 + O(\frac{1}{\sqrt{d}})$. This result is quite striking as it is precisely for random $d-$ regular graphs (with $d\in O(1)$) that the hardness of approximation for max cut using the Sherali-Adams hierarchy was shown \cite{chan2016approximate,charikar2009integrality,de2007linear,kothari2021approximating}. These two claims are the content of the next two corollaries.

\begin{corollary}\label{cor:1}
 Let $G=G(n,p)$ be sampled according to the Erd\H{o}s-R\'enyi model \cite{frieze2016introduction} where $p$ is a function of $n$. Let $g(n)$ be a non-decreasing function of $n$. Then, the ratio $\frac{\frac{1}{2}m+\frac{1}{4}\chi(G)}{\frac{1}{2}m}$ is at most $1+\sqrt{\frac{2}{g(n)}}$ as long as $np$ is at least
 $\frac{g(n)}{2}\log(n)$, with high probability. In particular, for all dense graphs of Theorem \ref{theorem:hardnessLP}, $np\geq g(n)\geq\sqrt{n}$ and the quotient converges to $1$. For sparse graphs of Theorem \ref{theorem:hardnessLP} where $np=c\log(n)$, the quotient is at most $1+\sqrt{\frac{2}{c}}$.
\end{corollary}
\begin{proof}
 Let $p = p(n)$ and $G=G(n,p)$ be sampled according to the Erd\H{o}s-R\'enyi model with $np \in \Omega(\frac{g(n)}{2} \log(n))$. Letting $\varepsilon = \frac{1}{n}$ and applying Theorem $1$ of \cite{chung2011spectra} we have that with probability at least $1-\frac{1}{n}$
 
 \[
 -\lambda_n \leq \sqrt{4np\ln\left(\frac{2n}{\varepsilon}\right)} +p.  
 \]

Recalling that the number of edges $m$ of $G$ is $\theta(n^2p)$ with high probability, a direct computation of the quantity  $\frac{\frac{1}{2}m-\frac{1}{4}n\lambda_n}{\frac{1}{2}m}$ gives the result.
\end{proof}

\begin{corollary}
Suppose that $G$ is a $d$-regular graph with $-\lambda_n \leq c\cdot\sqrt{d}$ for some constant $c$ and  $\mathcal{S} = \mathcal{E}(W)$. Then, the following inequality holds: 
\[
\frac{z_{SP_{\mathcal{S}}} }{z_{SD_{\mathcal{S}}} } \leq 1 + \frac{c}{\sqrt{d}}.
\]
\end{corollary}

\begin{proof}
Recall that a $d$-regular graph has $m= \frac{nd}{2}$ edges. This gives $n = \frac{2m}{d}$. 
Suppose $-\lambda_n \leq c \cdot\sqrt{d}$. Then, by Theorem \ref{teo:eigenbound} and that $z_{GW}\geq 0$ for any graph $G$ we get
    \[
    \frac{ \frac{1}{2}m+ \frac{1}{4}z_{SP_{\mathcal{S}}}}{ \frac{1}{2}m + \frac{1}{4} z_{GW}} \leq \frac{\frac{1}{2}m-\frac{1}{4}n\lambda_n}{\frac{1}{2}m} = \frac{\frac{1}{2}m-\frac{1}{4}\lambda_n \frac{2m}{d}}{\frac{1}{2}m} = 1-\frac{\lambda_n}{d} \leq 1+\frac{c}{\sqrt{d}}.
    \]
\end{proof}

It is known that random $d-$regular graphs satisfy the hypothesis of the theorem \cite{feige2005spectral,friedman1989second,tikhomirov2019spectral}, justifying our previous claim on the guarantees of our linear relaxation on random $d-$regular graphs. Another class of graphs which satisfies the hypothesis of the theorem are the Ramanujan expander graphs \cite{lubotzky1988ramanujan}, where $c = 2$. We contrast this result with the fact that the relative error of $\alpha(G)$- defined above in LP(\ref{usualLP})- relative to the max cut of $G$ tends to $1$ for Ramanujan graphs \cite{poljak1994expected}.

To the best of our knowledge, this is the first linear relaxation of max cut with these two guarantees.

\subsection{Hidden basis property and stronger guarantees}

In the previous subsection, we considered a bound given by Theorem \ref{teo:dual_bound} using the fact that $-\lambda_n\Vec{1}$ is feasible for program \ref{DGW}. However, it might very well be the case that there are other ``hidden" dual feasible solutions. Although we are not aware of any such solutions, it is illustrative to check whether or not program \ref{Srel} gives better solutions than the eigenvalue bound. This raises the question of the quality of our linear relaxation in the setup where the eigenvalue bounds fails to give an approximation factor better than $2$ for the maximum cut value of a graph.
Indeed, the eigenvalue bound is not powerful enough to provide an approximation factor better than $2-\varepsilon>0$ for any given $\varepsilon>0$ in general. As a matter of fact, for any given $\varepsilon>0$, there exist a family of graphs whose maximum cut is bounded above by $\frac{1}{2}m+\varepsilon m$, but the eigenvalue bound cannot certify a bound better than $2-\varepsilon$. We give an example of such as class in our next definition, which is inspired by a remark in \cite{o2018sherali}.

\begin{definition}
We say that a graph $G$ is sampled from the class of random graphs $\mathcal{G}(n,d,l)$ if $G$ has $n$ vertices and two disjoint components $G_1$ and $G_2$ where $G_1$ is a random $d$-regular graph, $G_2$ is complete bipartite graph where each side of the bipartition has $\sqrt{n}$ nodes, and $l$ random edges connect the $G_1$ and $G_2$. 
\end{definition}
Observe that the absolute value of most negative eigenvalue of the adjacency matrix of a graph sampled from $\mathcal{G}(n,d,l)$ is $\Omega(\sqrt{n})$ due to the bipartite component. If $d \in O(1)$ then the number of edges in $G$ is linear in $n$ and so is the maxcut of $G$. However, the eigenvalue bound is weak: it certifies that the maxcut size is at most $O(n^{1.5})$ (notice that this is worse even than the trivial upper bound of $m$). This class of graphs is suggested as an example in \cite{o2018sherali} as a class of graphs where the eigenvalue bound behaves poorly. However, our LP certifies a much better value, when $l=0$, as the next observation shows:

\begin{lemma}\label{removeEdges}
Let $G$ be a graph with two disconnected components $G_1$ and $G_2$, where
$|V(G_1)|= n_1$, $|V(G_2)|= n_2$, $\lambda^1$ is the smallest eigenvalue of the adjacency matrix of graph $G_1$ and $\lambda^2$ is the smallest eigenvalue of the adjacency matrix of $G_2$. Let $\mathcal{S} = \mathcal{E}(W)$. Then, $SP_\mathcal{S}$ certifies:

\[
z_{SP_{\mathcal{S}}} \leq  n_1\lambda^1 + n_2\lambda^2.
\]
\end{lemma}
\begin{proof}
    The proof is basically the same as the proof of Theorem  \ref{teo:eigenbound} by observing that the support of eigenvectors corresponding to disjoint components of a graph are disjoint.
\end{proof}

\comment{
Intuitively, this result demonstrates that our LP ``realizes" that $G$ consists of two components and and scales the quantities $n_1$ and $n_2$ by the "correct" value -$\lambda^1$ and $\lambda^2$ respectively- while the eigenvalue bound is oblivious to this.
}

This result may seem artificial in the sense that $G$ is a disconnected graph. However, we show through extensive experiments in Tables \ref{tab3} and \ref{tab4} in Section \ref{revised_sec:5} that the quotient of the optimal value \ref{Srel} to the GW relaxation is significantly better than the quotient of $\chi(G)$ to the GW relaxation, even when edges are added between the components in these difficult examples for the eigenvalue bound.\\

\subsection{Solvability of maxcut under $\mathcal{O}$}

In the previous subsection, we have seen that we can derive a good starting set $\mathcal{S}$. In general, program \ref{LSDP} does not solve the max cut SDP.
In this subsection we will show that whenever $G$ is a \textit{distance regular graph} then we have solvability of the max cut SDP under $\mathcal{O}$.
The class of distance-regular graphs contains strongly regular graphs, which have been extensively studied for their algebraic, combinatorial and spectral properties \cite{distRegGraph,spielman2012spectral}. Famous graphs such as the Petersen graph belongs to this class. In what follows, we give a sufficient condition that ensures that the value of \ref{Srel} equals the optimal value of the GW relaxation, provided that $\mathcal{S}$ includes an orthonormal eigenbasis of $W$.

\begin{definition}[Distance-regular graphs]
For a graph $G$ and $u,v$ vertices in $V(G)$ define $G_{j}(u)$ to be the set of vertices of $G$ at distance exactly $j$ of $u$, i.e., the vertices $v \in V(G)$ such that the shortest path joining $u$ and $v$ has length $j$. We say $G$ is distance regular if it is connected, $d$-regular for some $d$ and there exists integers $c_i, b_i, \  i\in \mathbb{N}$ such that for any two vertices $u,v$ at distance $i=d(u,v)$ there are precisely $c_i$ neighbours of $v$ in $G_{i+1}(u)$ and $b_i$ neighbours of $v$ in $G_{i-1}(u)$.
\end{definition}

Examples of such graphs are all strongly regular graphs, Hamming graphs, complete graphs, cycles, and odd graphs (such as the Petersen graph) \cite{distRegGraph}. The next theorem will allow us to prove that our linear relaxations are tight for this class of graphs.

\begin{theorem}\label{teoTight}
Let $G$ be a graph and $W$ its adjacency matrix. Let $\mathcal{S}= \mathcal{E}(W)$ and $W_n$ be the eigenspace of $W$ corresponding to $\lambda_n$. Suppose the dimension of $W_n$ is $k$ with $n > k\geq 1$. Suppose there exists an orthonormal basis  $\mathcal{U}=\{u_1,\dots,u_k\}$ of $W_n$ such that the matrix $A$ with rows $u_1,\dots,u_k$ has columns with constant $2$- norm, i.e. there exists some $c \in \mathbb{R}^{+}$ such that 
$\|A_{j}\|_2 = c \ \forall j\in [n] $ where $A_j$ denotes the $j$-th column of $A$. Then, $z_{SD_{\mathcal{S}}}$ equals $-n\lambda_n$ and in particular 
\[z_{SP_{\mathcal{S}}}=z_{GW} =z_{SD_{\mathcal{S}}}.\]

\end{theorem}
\begin{proof}
The proof requires two steps. First, we show that if such basis $\mathcal{U}$ exists and we let $\mathcal{S} = \mathcal{U}$ then the theorem holds. Second, we show that 
we may set $\mathcal{S}$ to be an arbitrary orthonormal basis of $W_n$. This is necessary since the dimension of $W_n\geq 2$ and hence orthonormal bases are not unique. This might break the theorem if we choose any other orthonormal basis for $\mathcal{S}$ instead of $\mathcal{U}$. 
We begin with the first step. Notice that $c = \sqrt{\frac{k}{n}}$. 
Indeed, since the $u_i$ are unitary vectors we have that for all $ \ i \in [k]$ $\sum_{j=1}^n A^2_{i,j}=1 $. Summing over $i$ gives $\sum_{i=1}^k\sum_{j=1}^n A^2_{i,j}=k$.
By our assumption of constant sum of the column vectors, we get   
$\sum_{i=1}^k A^2_{i,j}=c^2 \ \forall j\in [n]$. Summing over $j$ gives 
$\sum_{j=1}^n\sum_{i=1}^k A^2_{i,j}=nc^2$ and we get $k = nc^2$.
Let $B= \sqrt{ \frac{n}{k}}A^\top$ and $Y = BB^\top$. Let $v_i$ denote the $i$th row of $B$, and recall that $v_i$ has norm $\sqrt{\frac{k}{n}}$. This implies that $Y_{ii} = v_i \cdot v_i = \frac{n}{k} \cdot \frac{k}{n} =1 $.
Finally, observe that $Y = \frac{n}{k}\sum_{i}^k u_i(u_i)^\top$.
It follows that $Y$ is feasible for \ref{DSDGW} with $\mathcal{S}=\mathcal{U}$. This solution has an objective value

 \[
Z_{SD_{\mathcal{S}}} \geq   \left \langle -W, B^TB\right \rangle \geq \left\langle -W, \frac{n}{k}\sum_{i=1}^k u_iu_i^T\right \rangle   = -n\lambda_n.
\]
For the second part, we show that we can take $\mathcal{S}$ to be any arbitrary orthonormal basis of $W_n$. Notice that the only fact that we used from $\mathcal{U}$ is that the matrix $A$ formed by stacking the vectors $u_i$ as rows has constant column norm. Therefore, it suffices to show that any matrix $A'$ formed in the same way from an arbitrary basis $\mathcal{U}'$ has this same property. Hence, let $\mathcal{U}' = \{w_1,\dots,w_k\}$ be an arbitrary basis of $W_n$ and suppose that the basis $\mathcal{U}$ exists.

Since the vectors$\{u_1,\dots,u_k\}$ are an orthonormal basis of $W_n$ which is a lineal subspace of $\mathbb{R}^n$, we can extend this set of vectors to a full orthonormal basis $\{u_1,
\dots,u_k,u_{k+1},\dots u_n\}$ of $\mathbb{R}^n$. 
Further, observe that $\sum_{i=1}^n u_i(u_i)^\top = I_n$ where $I_n$ is the $n\times n$ identity matrix. To see this, let $v = r_1u_1 +\dots+r_nu_n\in \mathbb{R}^n$ be an arbitrary vector expressed in the $u_i$, $i \in [n]$ basis. We have 

\begin{equation}\label{eq:idmatrix}
\left(\sum_{i=1}^nu_iu_i^\top\right)v = \sum_{i=1}^n \langle u_i,v \rangle u_i = \sum_{i=1}^n r_iu_i  =  v. 
\end{equation}
We derive that $\sum_{i=1}^nu_iu_i^\top$ equals the identity matrix.
Notice that this equation remains true if we replace the $u_i$ for any arbitrary orthonormal basis of $\mathbb{R}^n$. Since the diagonal entries of $A^\top A $ equal $\frac{k}{n} = c^2$ we see that the diagonal entries of $\sum_{i=k+1}^n u_iu_i^\top$ equal $1-c^2$. Finally it follows that 
$\{w_1,\dots,w_k,u_{k+1},\dots,u_{n}\}$ is as well a basis for $\mathbb{R}^n$ and thus by Equation (\ref{eq:idmatrix}) we have $\sum_{i=1}^kw_iw_i^\top + \sum_{i=k+1}^n u_iu_i^\top = I_n$. This shows that every diagonal entry of the matrix $\sum_{i=1}^kw_iw_i^\top$ must equal $c$, and hence the matrix $A'$ formed by stacking the vectors $w_i$ as rows has constant column norm.
The conclusion of the theorem follows from the inequality
$Z_{SD_{\mathcal{S}}} \leq Z_{GW}\leq Z_{SP_{\mathcal{S}}} \leq  -n\lambda_n $.
\end{proof}

Alon and Sudakov proved something similar to the first part of our proof in  \cite{alon2000bipartite}. In the paper, the authors prove that $z_{GW} = \frac{1}{2}m -\frac{1}{4}n\lambda_n$ under the hypothesis that there exists a feasible solution $Y= B^\top B$ for the ($GW$) relaxation such that the columns of $B$ are unitary vectors $v_1,\dots, v_n$ and its rows $u_1,\dots u_k$, $1\leq k\leq n$ are eigenvectors of $W$ corresponding to $\lambda_n$. We conclude this section with the corollary for distance-regular graphs.

\begin{corollary}

Let $G$ be a distance-regular graph. Let $\mathcal{S} = \mathcal{E}(W)$. Then
\[z_{SP_{\mathcal{S}}}=z_{GW} =z_{SD_{\mathcal{S}}}.\]
\end{corollary}

\begin{proof}
 
The results follows from the following theorem. It states that the eigenspaces of distance regular graphs satisfy the hypothesis of Theorem \ref{teoTight}.
 
\begin{theorem}[\cite{distRegGraph}, Theorem $4.1.4$]
Let $G$ be a distance regular graph and $\lambda$ an eigenvalue of $G$. Then, there exists a symmetric matrix whose columns span the eigenspace corresponding to $\lambda$ and that have a constant norm.
\end{theorem}
\end{proof}

\section{Applications to semidefinite programs} 
\label{revised_sec:4}

To verify the applicability of the ideas presented, we consider three families of semidefinite optimization problems, each illustrating an aspect of our work. The first problem considered is the semidefinite relaxation of the maxcut problem which we presented in Section \ref{revised_sec:3}. We present our experimental results in Section \ref{revised_sec:5}.

\subsection{Maximum Cut}

The max cut  problem is a prime example of how our methodology can be applied as it is a hard combinatorial problem that linear programs fail to approximate. We will test our ideas using two linear programs, already introduced in Section \ref{revised_sec:3}.

\begin{equation}\tag{$SP_{\mathcal{S}}$}
\begin{aligned}
&\max_{X \in \mathbb{S}^n} \ \frac{1}{2}m + \frac14 \langle -W,X\rangle   \\
\text{ s.t: } v^\top Xv  \geq  0 \ & \forall \ v \in \mathcal{S},  X_{ii}=1, \ \forall \ i \in [n],
  \| X \|_\infty \leq 1.
\end{aligned}
\end{equation}

By the results of Section \ref{revised_sec:3}, we know that as $n \rightarrow +\infty$, the optimal value of this program will converge to the optimal value of max cut for Erd\H{o}s-R\'enyi graphs and random $d$-regular graphs whenever $\mathcal{S}$ contains a basis of eigenvector of the matrix $W$. We test the quality of the linear relaxation on such graphs, as well as on graphs of the family $\mathcal{G}(n,l,k)$ which was introduced in Section \ref{revised_sec:3}. This family was designed to have a trivial eigenvalue bound. In Appendix \ref{AppendixA} we include as well experiments on the quality of relaxations on
$16$ graphs taken from TSPLIB \cite{reinelt1991tsplib} and $14$ graphs from the network repository \cite{nr-aaai15}. Furthermore, we consider program
\begin{equation} \tag{$SD_{\mathcal{S}}$}
\begin{aligned}
\frac{1}{2}m + \frac{1}{4}& \max_{\eta \in \mathbb{R}^k} \left \langle -W,\sum_{i=1}^k \eta_ix_ix_i^\top \right \rangle \\
 \text{s.t: } diag \left(\sum_{i=1}^k \eta_ix_ix_i^T \right) \leq  1, & \ \eta_i \geq 0, \ x_i \in \mathcal{S} \ \forall i \in [k], \ k=|\mathcal{S}|.
    \end{aligned}
\end{equation}
This program is useful as we can obtain graph cuts from its solution using the rounding procedure of Goemans and Williamson~\cite{goemans1995improved}. Since the focus of this paper is comparing the optimal value of the different linear relaxations versus the optimal value of the SDPs, we defer results on rounded solutions to Appendix \ref{cutGeneration}.

\subsection{Lov\'asz theta number}

The second problem we consider is the Lov\'asz theta number $\vartheta(G)$ introduced by Lov\'asz in the seminal paper \cite{Lovasz1979shannon} as a convex relaxation for the stability number of a graph $G$. $\vartheta$ can be computed in polynomial time using a semidefinite program. Since $\vartheta(\Bar{G})$ -where $\Bar{G}$ is the complement of $G$- is lower and upper bounded resp. by the clique number and the chromatic number of $G$, it allows one to compute those numbers in polynomial time for graphs for which these two quantities coincide e.g., perfect graphs. 

$\vartheta(G)$ can be computed by the following semidefinite optimization program:
\begin{equation}\label{theta}\tag{$Tn$}
\begin{aligned}
\max_{S\in \mathbb{R}^{n}} & \  \langle  J,X\rangle \\
\text{ s.t: } tr(X)  = 1, & \ X_{i,j}=0 \  \forall (i,j)\in E,  \\ 
   X \succeq & \  0.
\end{aligned}
\end{equation}
This problem is related to our setup, as it is known that the feasible region of the dual program is polyhedral whenever the considered graph is perfect. This striking results coincides with the fact that it is precisely for these graphs where the theta number coincides with the \textit{independence number} of the graph.

We apply the ideas developed in Section \ref{revised_sec:2} on two families of graphs. The first class is that of \textit{regular} graphs. Notice that the constraints $ X_{i,j}=0 \  \forall (i,j)\in E$ can be expressed as $\left \langle X,A^{ij}\right \rangle = 0$ where $A^{ij}$ is matrix of all zeros except it has a $1$ in its $ij$,$ji$ entries whenever $G$ contains edge $ij$. It is clear that if $A$ is the adjacency matrix of $G$, we have $A = \sum_{ij,\in E}A^{ij}$. Regular graphs are interesting in our setting as it is easy to check that whenever $G$ is regular graphs, $A$ its adjacency matrix, and $J$ the matrix of all ones, we have $JA = AJ$. The second class of graphs we consider are Erd\H{o}s-R\'enyi random graphs, which are typically not regular and it is not obvious how to combine the $A^{ij}$ to obtain a matrix that commutes with $J$. We will use program (\ref{conmute_generator}) to find such matrices.

Given a finite set $\mathcal{S}$,
we obtain the linear relaxation of program (\ref{theta}):
\begin{equation}\label{lineartheta}\tag{$LTn$}
\begin{aligned}
\max_{X\in \mathbb{S}^{n}} & \  \langle  J,X\rangle \\
\text{ s.t: } tr(X)  = 1, & \ X_{i,j}=0 \  \forall (i,j)\in E,  \\ 
  v^\top Xv \geq &  \ 0 \ \forall \  v\in \mathcal{S}
\end{aligned}
\end{equation}
In Section \ref{revised_sec:5}, we compare the objective value of programs (\ref{theta}) and (\ref{lineartheta}), on Erd\H{o}s-R\'enyi random graphs and $d-$regular graphs. Interestingly, this problem is much more resistant to the the cut generation strategy for solving the corresponding SDP proposed in Algorithm \ref{alg1}. As we will see, generating cuts through the separation oracle of the semidefinite cone fails completely on both Erd\H{o}s-R\'enyi graphs and $d-$regular graphs. On the contrary, setting $\mathcal{S}$ to be the columns of a matrix that simultaneously diagonalizes $J$ and $A$ -where $A$ is the adjacency matrix of $G$ in the case of regular graphs or a matrix given by program (\ref{conmute_generator}) in the case of Erd\H{o}s-R\'enyi graphs - performs significantly better.

In our discussion on the max cut problem we showed that there is a eigenvalue bound for the max cut value that every graph satisfies, and one might wonder if there such a bound for the theta number. This is indeed the case, albeit only for regular graphs.

\begin{remark}
Let $G$ be a $d-regular $ graph with $n$ vertices. Let $W$ be the adjacency matrix of $G$ with largest eigenvalue $\lambda_1$ and smallest eigenvalue $\lambda_n$, then the Lov\'asz theta number $\vartheta(G)$ satisfies

\begin{equation}\label{eq:theta_eigen_bound}
    \vartheta(G) \leq \frac{-n\lambda_n}{\lambda_1 - \lambda_n}
\end{equation}

For a proof of this result, see \cite{Lovasz1979shannon}. 
We conjecture that the objective value of the linear program (\ref{lineartheta}) is also upper bound by $\frac{-n\lambda_n}{\lambda_1 - \lambda_n}$ as this was the case in all the experiments we performed for $d-$regular graphs.

\end{remark}


\subsection{QCQPs}\label{sec:qcqp}

We consider more general SDPs obtained as the Shor relaxation \cite{shor1990dual} of certain QCQPs to test the proposed methodology in three different settings, each highlighting an interesting point. General QCQPs were introduced in Section \ref{revised_sec:2}, but in this section and Section \ref{revised_sec:5} we will consider a more specialized version of them, following \cite{bao2011semidefinite}, of the form

\begin{equation}
\begin{aligned}
    \inf_{x \in \mathbb{R}^n}  x^\top & A_0 x + b_0^\top x +c_0 \\
    \text{s.t: }  x^\top A_i x  + & b_i^\top x   \leq b_i \ \forall i \in [r], \\
    D x & = d,\\
    l \leq & x  \leq u,
\end{aligned}
\end{equation}
where $r$ denotes the number of quadratic constraints and is at least $1$.
$A_i, \ i =\{0,\dots,r\}$ are symmetric matrices, not necessarily PSD, $b_i, \  i = \{0,\dots r\}$ are vectors in $\mathbb{R}^n$, $D$ is a  $q\times n$ real matrix and 
$d \in \mathbb{R}^q$. $l$ and $u$ are vectors in $\mathbb{R}^n$ and we assume that
$-\infty < l\leq u < +\infty$ so that the bounding boxes are non-empty and bounded. If the bounding boxes are of the form $[l,u]^n$ we can do a linear change of variables so that $x \in [0,1]^n$.
Such problems admit the following SDP relaxation:

\begin{equation}\label{shor_qcqp}
\begin{aligned}
&\inf_{x \in \mathbb{R}^n, X\in \mathbb{S}^n}  \left \langle A_0,X \right \rangle + b_0^\top x +c_0 \\
   & \text{s.t: }   \left \langle A_i,X \right \rangle +  b_i^\top x\leq  \ c_i \ \forall i \in [r], \\ 
    & \ \ \ \ \ \ \ \ \ \ \ \ \ \ \ \ D  x  = d,\\
    & \ \ \ \ \ \ \ \ \  0 \leq  x_i  \leq 1 \ \forall i \in [n], \\
&\ \ \ \ \ \ \ \ \   0 \leq  X_{i,j} \leq 1 \ \forall i,j \in [n], \\
 & \ \ \ \ \ \ \ \ \ \ \ \ \ \ \ \begin{bmatrix}
X & x \\
x^\top & 1 
\end{bmatrix} \succeq 0.
\end{aligned}     
\end{equation}

By letting
$
\Hat{A}_i:=\begin{bmatrix}
A_i & b_i \\
b^\top_i & c_i  
\end{bmatrix}, i\in \{0,1,\dots,m\}, \  \Hat{X} := \begin{bmatrix}
X & x \\
x^\top & 1 
\end{bmatrix}
$ and by $\Hat{X}_{n+1}$ the $n+1$'th column of $\Hat{X}$ we can write the previous problem in the \ref{SDP} form
\begin{equation}\label{shor_reformulated}\tag{$QSDP$}
\begin{aligned}
\inf_{\Hat{X}\in \mathbb{S}^{n+1}} & \left \langle \Hat{A}_0,\Hat{X} \right \rangle \\
 \text{s.t }:  & \left \langle \Hat{A}_i, \Hat{X} \right  \rangle \leq  \ 0 \ \forall i \in [r], \\
    &D \Hat{X}_{n+1}  = d,\\
 &0 \leq \Hat{X}_{i,j} \leq 1 \ \forall i,j \in [n+1], \\
 & X_{n+1,n+1} = 1, \\
&\Hat{X}  \succeq 0.
\end{aligned}     
\end{equation}

In Section \ref{revised_sec:5} we test our methodology on random QCQPs using instances generated as in \cite{bao2011semidefinite}.

\subsection{Quadratic knapsack problem}

An interesting point arises whenever the quadratic forms determining the objective and the constraints do not have linear and constant terms, i.e. $b_i = c_i =0 \ \forall i \in \{0,\dots,r\}$. In that case, our methodology takes $\mathcal{S} = \{v_1,\dots,v_{n+1}\}$ to be the eigenvectors of a matrix in $\mathbb{S}^{n+1}$ whose $n+1$'th row and column are $0$. Hence, the constraints 
$v ^\top  \Hat{X} v \geq 0$ in program \ref{shor_reformulated} essentially ignore the last row and column of $\Hat{X}$ and amount to the constraints $u_i^\top  X u_i \geq 0$
where $u_1,\dots,u_n$ are a basis of eigenvectors of an aggregation of the $A_i, \ i \in [r].$ This is a weaker constraint than what we actually want, which is $u_i^\top \left (X-xx^\top x \right ) u_i \geq 0 , \ i \in [n]$.

There are a few approaches we can consider to deal with this issue. 
For instance, we could choose to overlook it entirely and proceed by relaxing \ref{shor_reformulated} to an LP, ignoring that the $b_i$ are $0$. 
Alternatively, if we have a linear constraint $b_i^{\top} x = c $, we may set $\Hat{A}_0 =
\begin{bmatrix}
A_0 & b_i \\
b^\top_i & -2c 
\end{bmatrix}$
which shifts the objective by a constant. Finally, and perhaps more interestingly, we may use the constraints 
$u_i^\top \left (X-xx^\top \right ) u_i \geq 0 , \ i \in [n]$ directly, which can be equivalently rewritten as:
\begin{equation}\label{eq:soc}
u_i^\top  Xu_i \geq u_i^\top \left ( x x^\top \right ) u_i^\top = (u_i^\top x)^2 \ \forall i \in [n].
\end{equation}

These are \textit{second order cone constraints} which result in a second order cone relaxation of program \ref{QCQP} depending on a set $\mathcal{S}$ of vectors $u$ in $\mathbb{R}^n$. Such a program is both a relaxation of \ref{shor_reformulated}, and a strenghening of the linear relaxation that changes the constraint $\Hat{X}\succeq 0$ for 
$u^\top X u \geq 0$ with $u \in \mathcal{S}$, for any finite set $\mathcal{S}$.

We test these different possibilities in Section \ref{revised_sec:5} on instances of the \textit{Quadratic Knapsack problem } \cite{pisinger2007quadratic} which is a QCQP of the form

\begin{equation}\label{quadnapsack}\tag{QKP}
\begin{aligned}
\max_{x \in \mathbb{R}^n} & \  x^\top A_0 x \\
s.t:   \sum_{j=1}^k w_j x_j & \leq C, \  x \in \{0,1\}^n
\end{aligned}
\end{equation}

where $w \in \mathbb{R}^n$, $A_0 \in \mathbb{S}^n$, $C \in \mathbb{R}_+$. It has been noted in the literature that the usual Shor semidefinite relaxation of this program is not very strong \cite{helmberg2000semidefinite,pisinger2007quadratic} and one may add certain valid inequalities which result in the following tighter SDP:

\begin{equation}\label{strenghen_shor_quadnapsack}\tag{QKPSDP}
\begin{aligned}
&\max_{X \in \mathbb{S}^n} \left \langle A,X\right \rangle \\
s.t: &  \sum_{j=1}^n w_{j}X_{ij} - C X_{ii} \leq \  0 \ \forall \ i \in [n],  \\
&X-diag(X)diag(X)^\top \succeq 0.
\end{aligned}
\end{equation}

Using the idea before and a finite set $\mathcal{S}$ one may further relax this problem to obtain the second order cone program

\begin{equation}\label{soc_strenghen_shor_quadnapsack}\tag{QKSSOC}
\begin{aligned}
\max_{X \in \mathbb{S}^n} &\left \langle A,X\right \rangle & \\
s.t:  \sum_{j=1}^n w_{j} X_{ij} & - C X_{ii} \leq  \  0 \ \forall \ i \in [n],  \\
 u^\top Xu \geq & \left(u^\top diag(X)  \right)^2  \ \forall u \in \mathcal{S}.
\end{aligned}    
\end{equation}

\subsubsection{Extended trust region}

In the previous problems it is not obvious how to linearly combine the matrices $A_i, i\in [r]$, that determine the quadratic forms to form the identity matrix, and hence we cannot apply Theorem \ref{teo:dual_bound} directly to arbitrary QCQPs. This motivates us to consider a variation where the identity matrix is explicitly one of the constraint matrices. This is the case of the \textit{generalized trust region} problem \cite{locatelli2016exactness}. That type of QCQPs consists in minimizing a quadratic function over the intersection of the unit ball and some half-spaces:

\begin{equation}\label{ex_trust_reg}\tag{TR}
\begin{aligned}
\min_{x \in \mathbb{R}^n} &  \  x^\top Q x  + 2b^\top x,  \\
 \text{s.t: } & x^\top x \leq 1 \\
 & Dx \leq d
\end{aligned}     
\end{equation}

with $Q \in \mathbb{S}^n$, $b \in \mathbb{R}^n$, $D \in k \times n$ for some $k \in \mathbb{N}$ and $d \in \mathbb{R}^k$. Notice that the constraint $x^\top x \leq 1$ can be written as $x^\top I_n \ x \leq 1 $. In Section \ref{revised_sec:5} we test our methodology on a slightly more general version of this problem, where we keep some 
quadratic constraints. Abusing the language, we still refer to this family of problems as extended trust region problems.

\section{Experimental results}
\label{revised_sec:5}
In this section we present experimental results exhibiting the quality of our linear relaxations for the semidefinite relaxation of max cut, Lov\'asz's theta number and on the SDP relaxations of families of QCQPs described in Section \ref{revised_sec:4}. For each of these problems, we will compare the optimal value of the linear relaxations to the optimal value of the SDP which they respectively relax by means of the quotient of the objective values. We contrast these quotients to the alternative of using Algorithm \ref{alg1}, starting with $\mathcal{S}= e_1,\dots,e_n$ and iteratively generating cuts using the SDP separation oracle. Whenever we fix an semidefinite program with some label $SDP$, we denote by $Iter_k(SDP)$ the linear program obtained at the $k-th$ iteration of Algorithm \ref{alg1}. For instance, $Iter_0(SDP)$ is simply is dropping the semidefinite constraint of the SDP instance.
We define $z_n$ as the optimal value of $Iter_n(SDP)$. For each family of experiments, where we consider a certain $SDP$, we will denote by $z_{\mathcal{S}}$ the objective value of the corresponding linear relaxation obtained by following the ideas of Section \ref{revised_sec:1}. $z_{sdp}$ will denote the objective value of the SDP instance. Although we consider different SDPs, there will not be danger of confusion as we caption of the figures and tables indicate which SDP we are addressing.\\

All of the code used is available at
\url{https://github.com/dderoux/Instance_specific_relaxations}. To solve the resulting optimization programs we have used Mosek \cite{mosek}. \footnote{The experiments were performed on a $32$ GB RAM ThinkPad Lenovo T490s machine running windows 10 with a Intel(R) Core(TM) i7-8665U CPU @ 1.90GHz 2.11 GHz. }

\subsection{Max cut}

\begin{figure}[htbp]
\centering
\begin{adjustbox}{max width=1.3\textwidth}
\begin{subfigure}[b]{0.6\linewidth}
\centering
\includegraphics[width=\linewidth]{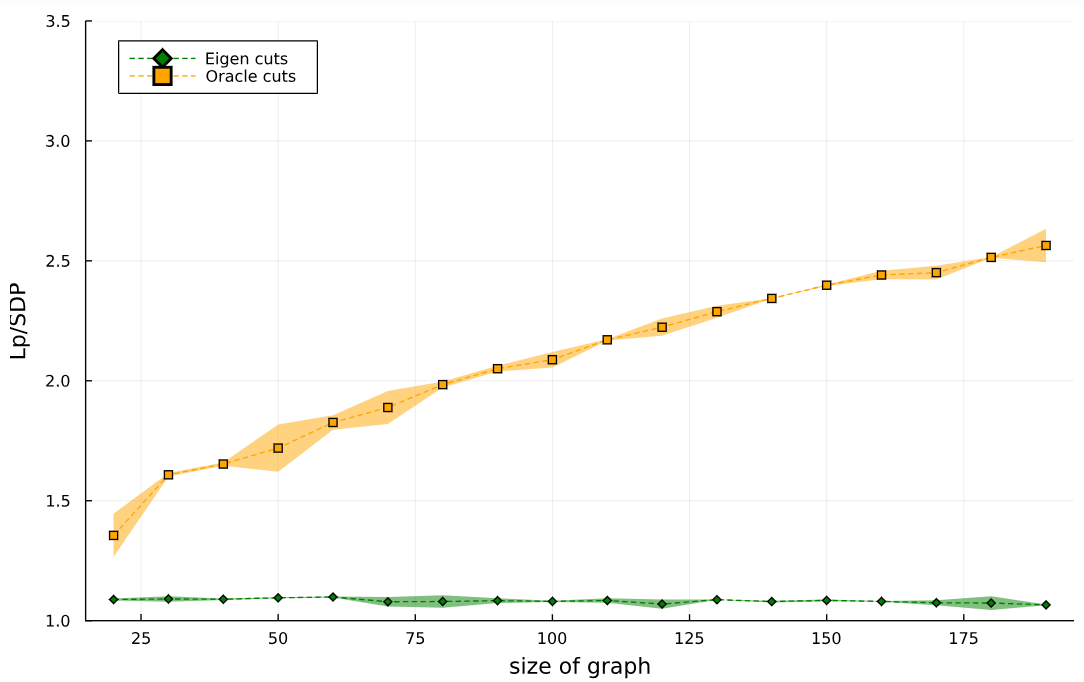}
\caption{$p=0.5$}
\end{subfigure}%
\begin{subfigure}[b]{0.6\linewidth}
\centering
\includegraphics[width=\linewidth]{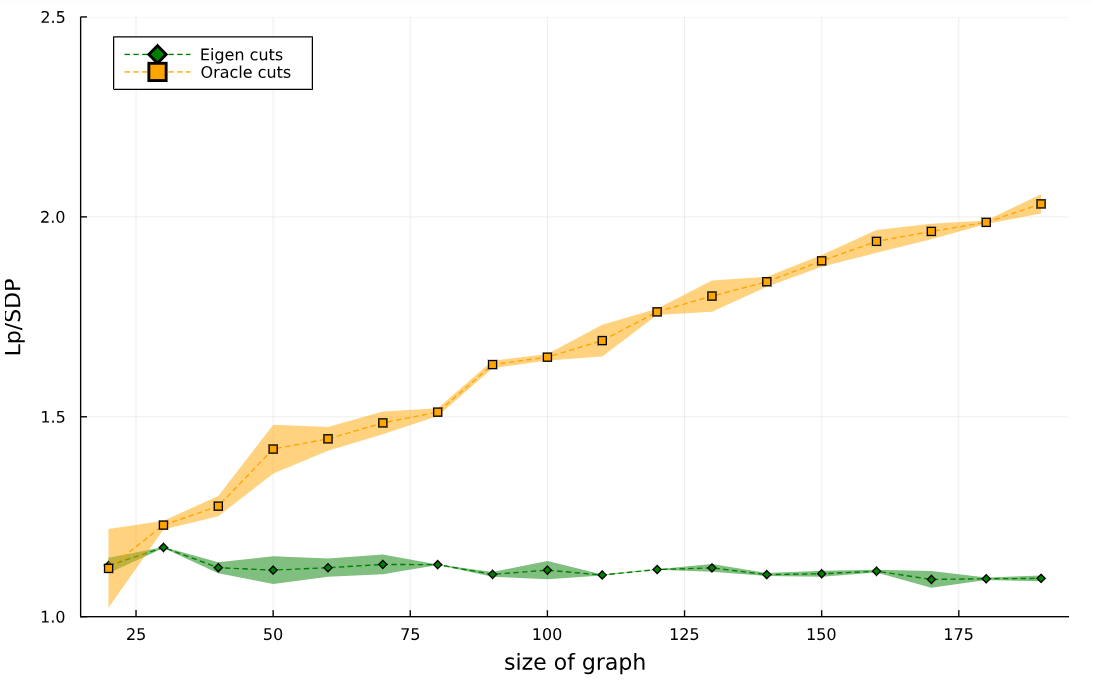}
\caption{$p=0.9$}
\end{subfigure}
\end{adjustbox}

\begin{subfigure}[b]{0.6\linewidth}
\centering
\includegraphics[width=\linewidth]{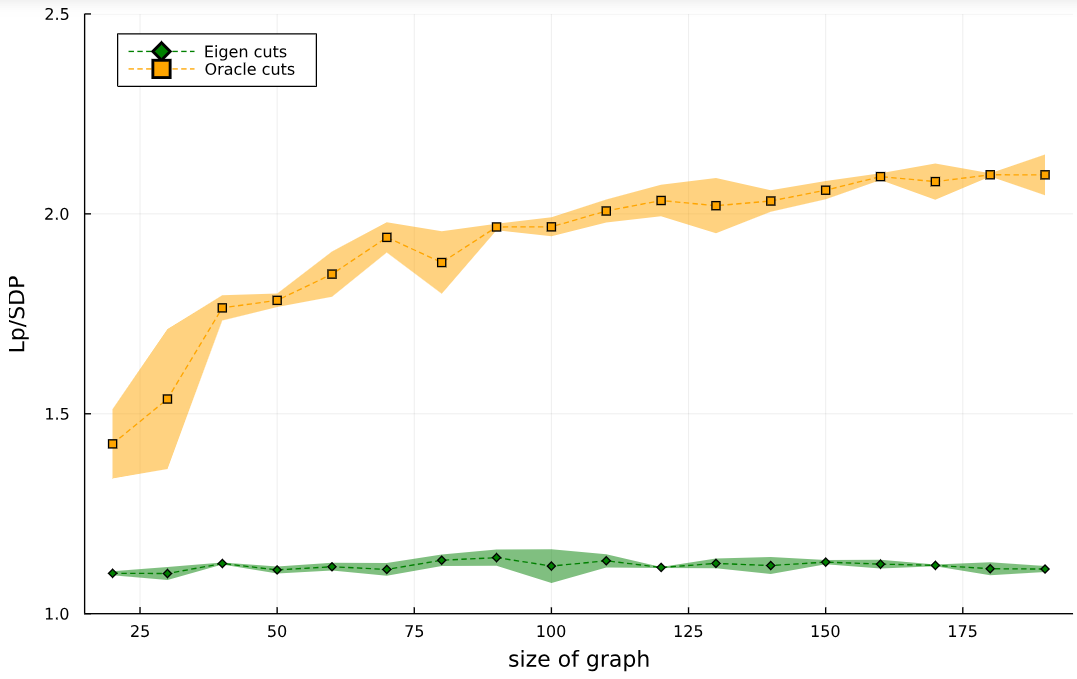}
\caption{$p=\frac{3\log(n)}{n}$}
\end{subfigure}

\caption{ Ratio of $\frac{z_\mathcal{S}}{z_{sdp}}$ (Eigen cuts) and $\frac{z_n}{z_{sdp}}$ (Oracle cuts) for instances of max cut where the graph has been sampled according to the   Erd\H{o}s-R\'enyi random model, for different values of $p$, as $n$ grows.}
\label{fig:mc_Erdos}
\end{figure}

\begin{figure}[htbp]
\centering
\begin{adjustbox}{max width=1.3\textwidth}
\begin{subfigure}[b]{0.6\linewidth}
\centering
\includegraphics[width=\linewidth]{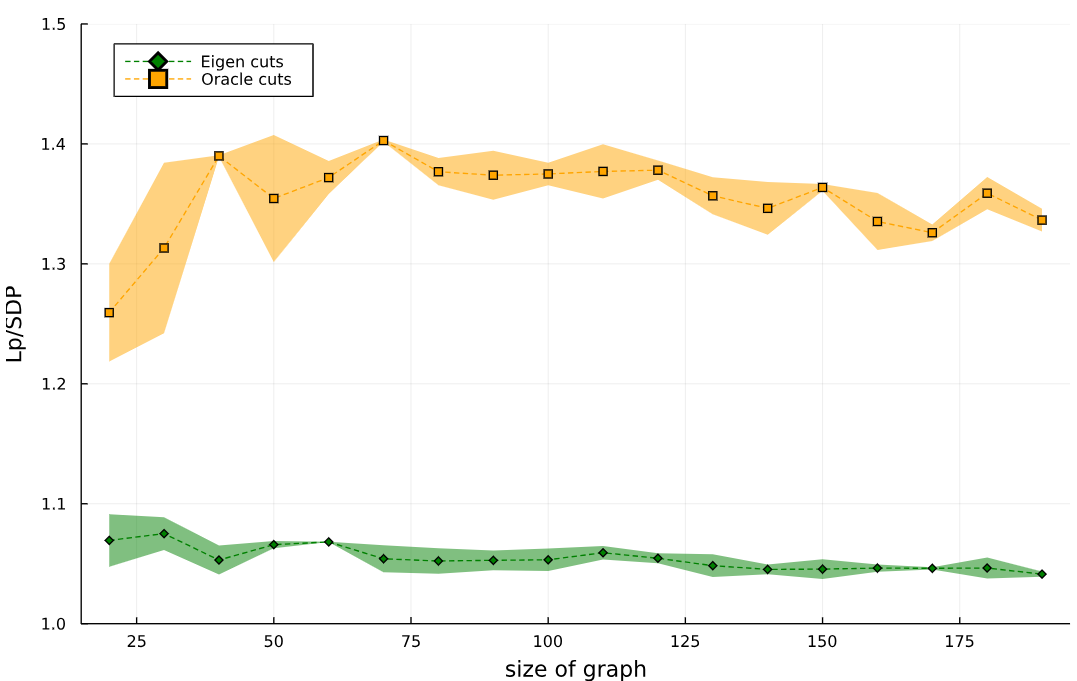}
\caption{$d=5$}
\end{subfigure}%
\begin{subfigure}[b]{0.6\linewidth}
\centering
\includegraphics[width=\linewidth]{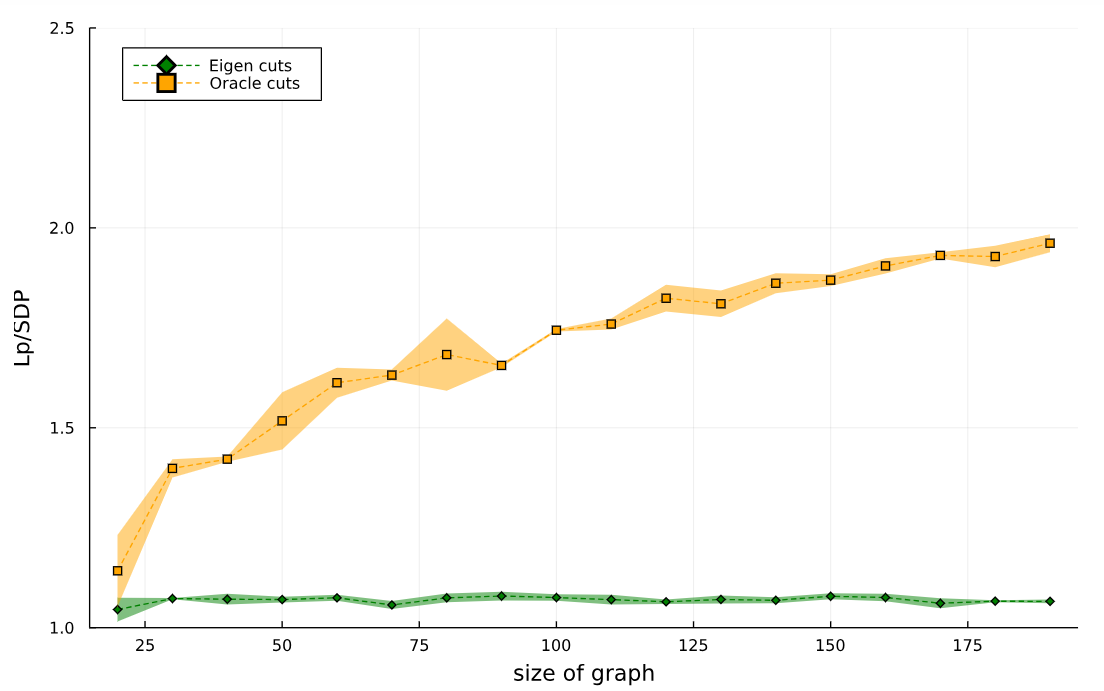}
\caption{$d=\sqrt{n}$}
\end{subfigure}
\end{adjustbox}

\begin{subfigure}[b]{0.6\linewidth}
\centering
\includegraphics[width=\linewidth]{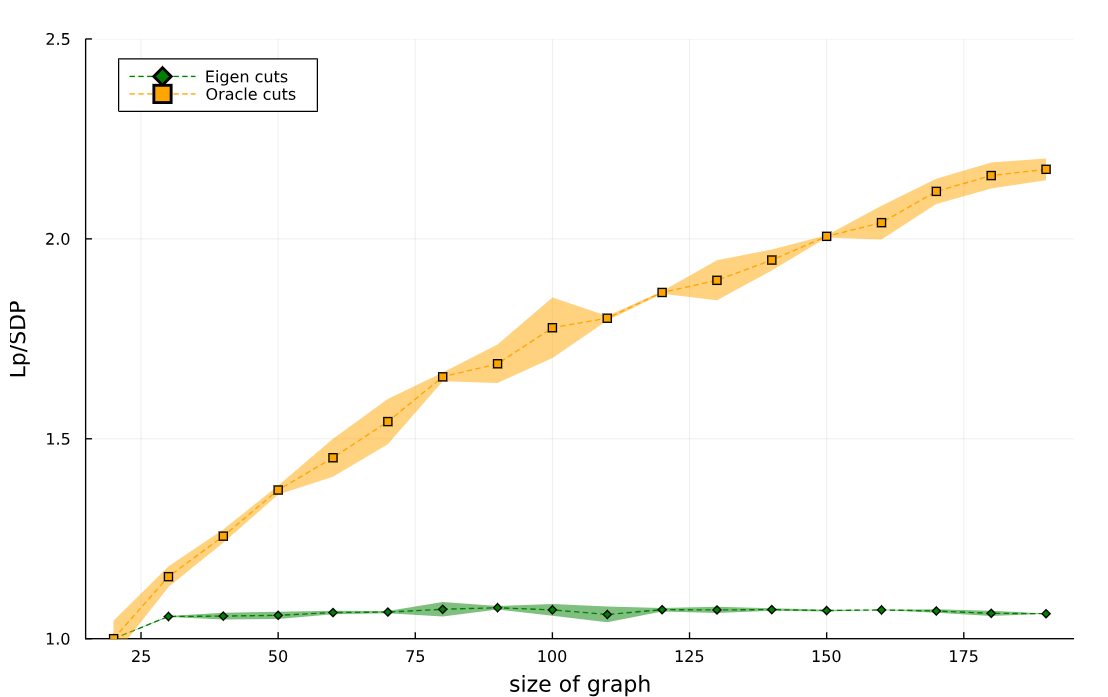}
\caption{$d=\frac{n}{10}$}
\end{subfigure}

\caption{ Ratio of $\frac{z_\mathcal{S}}{z_{sdp}}$ (Eigen cuts) and $\frac{z_n}{z_{sdp}}$ (Oracle cuts) for instances of max cut where the graph is a random $d-$regular graph, for different  values of $d$, as $n$ grows.}
\label{fig:mc_regular}
\end{figure}

Denote by $z_{sdp}$, $z_n$ and $z_\mathcal{S}$ the objective values of programs (\ref{GW}), $Iter_n(GW)$ and (\ref{Srel}) with 
$\mathcal{S}$ chosen as in Subsection \ref{subsec:finding_conmuting}. In this particular case, since the(\ref{GW}) semidefinite program does not have linear constraints beyond the ones of the diagonal, the identity $I$ is the only constraint matrix and it commutes with the objective matrix $W$. This means that $\mathcal{S}$ is simply a 
eigenbasis for the matrix $W$. As we proved in Theorem \ref{teo:eigenbound}, (\ref{Srel}) satisfies the eigenvalue bound for max cut.
For each $n$ ranging from $20$ up to $200$ in steps of $10$, we generate $5$ random graphs 
and plot the maximum, minimum and median of the quotients $\frac{z_\mathcal{S}}{z_{sdp}}$.
We present our results for   Erd\H{o}s-R\'enyi random graphs and $d-$regular random graphs in Figures~\ref{fig:mc_Erdos} and ~\ref{fig:mc_regular} respectively.

\comment{


\begin{figure}
    \centering
    \includegraphics[width=0.9\textwidth]{revised_pictures/mc_  Erd\H{o}s_05.png}
    \caption{ Ratio of $\frac{z_\mathcal{S}}{z_{sdp}}$ (Eigen cuts) and $\frac{z_n}{z_{sdp}}$ (Oracle cuts) for instances of max cut where the graph has been sampled according to the   Erd\H{o}s-R\'enyi random model with $p=0.5$ as $n$ grows.}
    \label{fig:mc_  Erd\H{o}s_05}
\end{figure}

\begin{figure}
    \centering
    \includegraphics[width=0.9\textwidth]{revised_pictures/mc_  Erd\H{o}s_09.png}
    \caption{ Ratio of $\frac{z_\mathcal{S}}{z_{sdp}}$ (Eigen cuts) and $\frac{z_n}{z_{sdp}}$ (Oracle cuts) for instances of max cut where the graph has been sampled according to the   Erd\H{o}s-R\'enyi random model with $p=0.9$ as $n$ grows.}
    \label{fig:mc_  Erd\H{o}s_09}
\end{figure}

\begin{figure}
    \centering
    \includegraphics[width=0.9\textwidth]{revised_pictures/mc_  Erd\H{o}s_log.png}
    \caption{ Ratio of $\frac{z_\mathcal{S}}{z_{sdp}}$ (Eigen cuts) and $\frac{z_n}{z_{sdp}}$ (Oracle cuts) for instances of max cut where the graph has been sampled according to the   Erd\H{o}s-R\'enyi random model with $p=\frac{3\log(n)}{n}$ as $n$ grows.}
    \label{fig:mc_  Erd\H{o}s_log}
\end{figure}

}


\comment{

\begin{figure}
    \centering
    \includegraphics[width=0.9\textwidth]{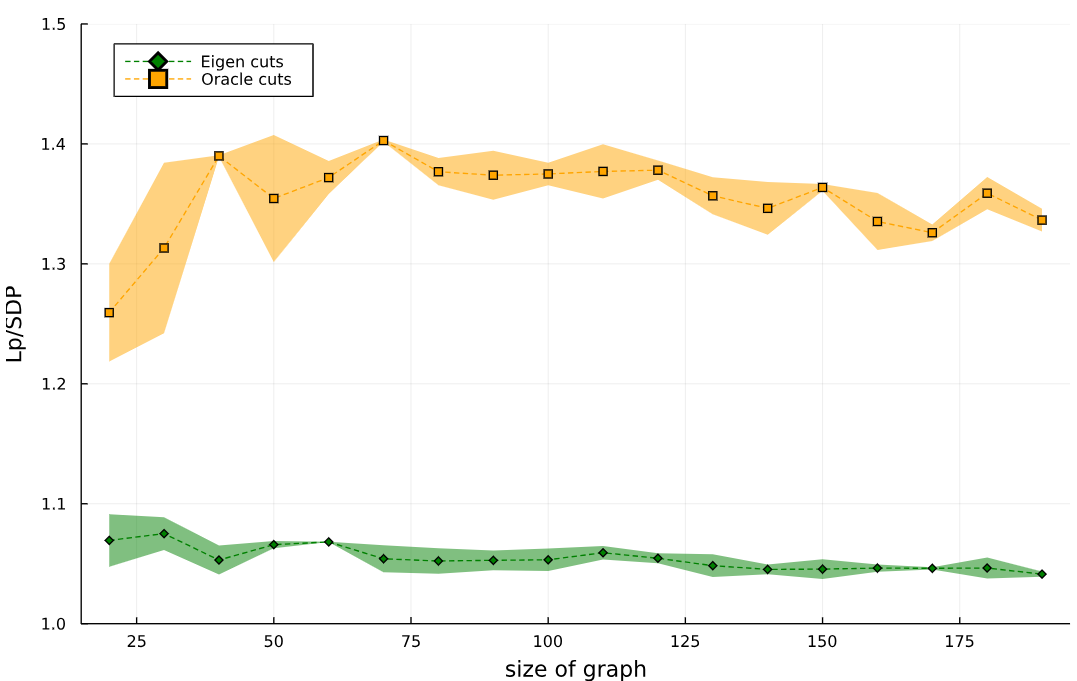}
    \caption{ Ratio of $\frac{z_\mathcal{S}}{z_{sdp}}$ (Eigen cuts) and $\frac{z_n}{z_{sdp}}$ (Oracle cuts) for instances of max cut where the graph is a random $5-$regular graph, and as $n$ grows.}
    \label{fig:mc_regular_5}
\end{figure}

\begin{figure}
    \centering
    \includegraphics[width=0.9\textwidth]{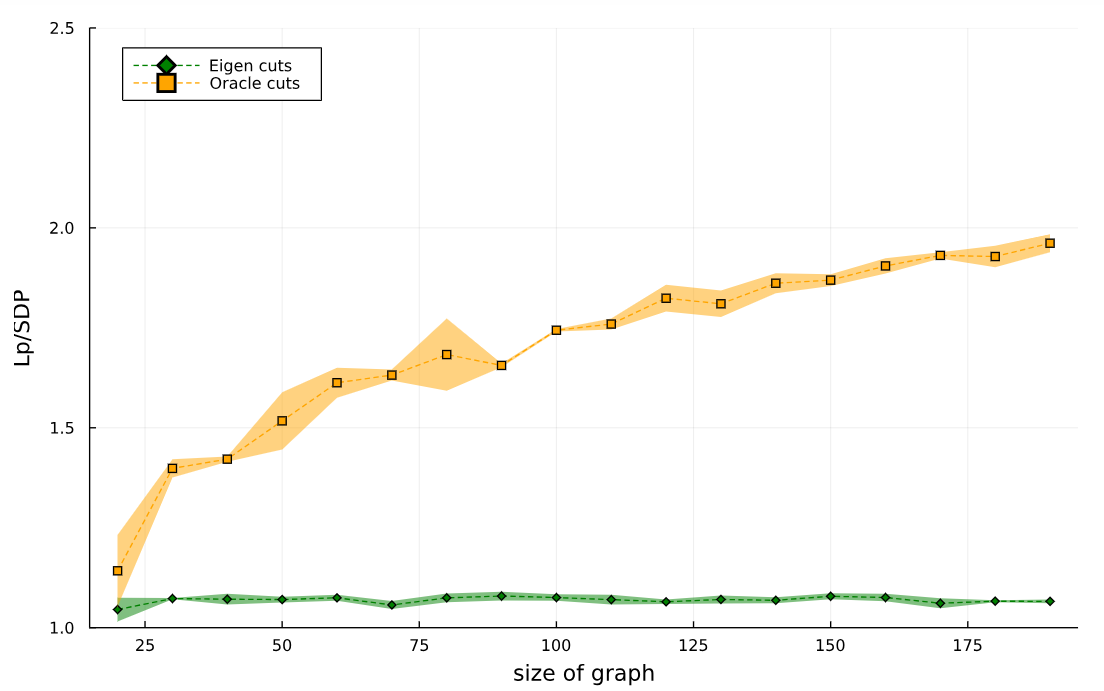}
    \caption{ Ratio of $\frac{z_\mathcal{S}}{z_{sdp}}$ (Eigen cuts) and $\frac{z_n}{z_{sdp}}$ (Oracle cuts) for instances of max cut where the graph is a random $d-$regular graph with $d = \sqrt{n}$, and as $n$ grows.}
    \label{fig:mc_regular_sqrt}
\end{figure}

\begin{figure}
    \centering
    \includegraphics[width=0.9\textwidth]{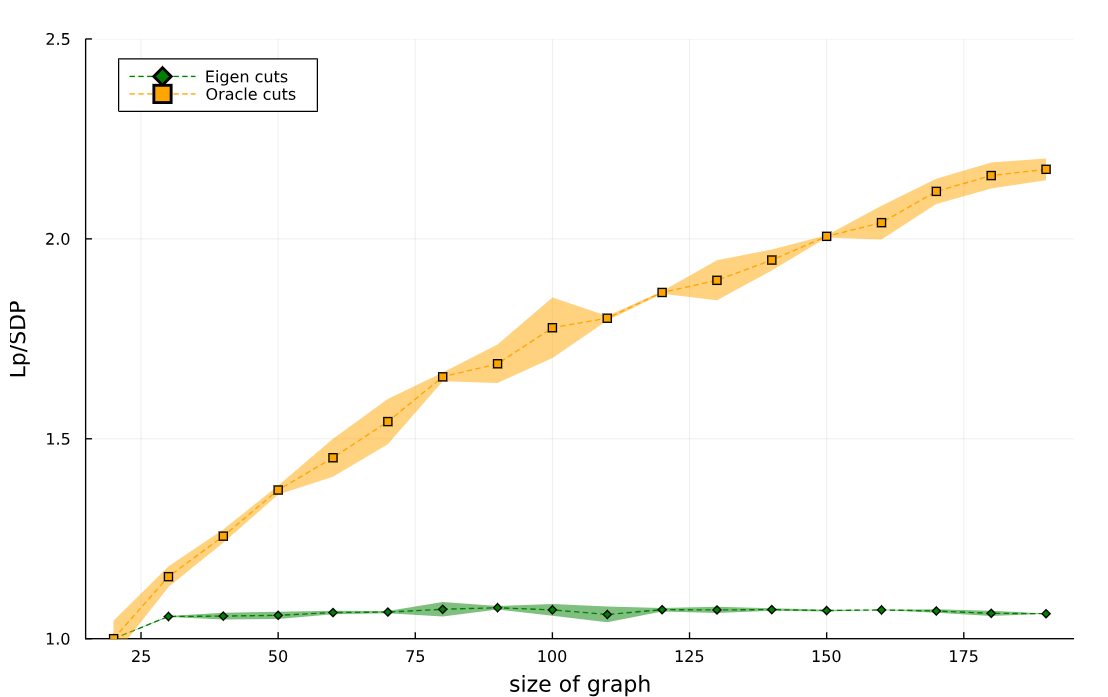}
    \caption{ Ratio of $\frac{z_\mathcal{S}}{z_{sdp}}$ (Eigen cuts) and $\frac{z_n}{z_{sdp}}$ (Oracle cuts) for instances of max cut where the graph is a random $d-$regular graph with $d = \frac{n}{10}$, and as $n$ grows.}
    \label{fig:mc_regular_linear}
\end{figure}

}

\subsubsection{Comparison with the Eigenvalue Bound}

In Tables \ref{table3} and \ref{table4}, we compare the performance of $z_{\mathcal{S}}$ and the eigenvalue bound $\chi(G):= -n\lambda_n(G)$ on the graphs $\mathcal{G}(n,k,l)$ which we introduced Section \ref{revised_sec:3}, for different values of $n$, $k$ and $l$. Since all of our experiments are random, we present averaged values over $5$ instances, as well as the standard deviations of our results. 
Notice that the eigenvalue bound fails to give a small upper bound on the max cut value for this family of graphs. For the case $n=400$, the bound fails completely, by giving a worse bound that the trivial upper bound of $m$ for max cut. However, the linear program succeeds, in all of our experiments, to have a quotient of at most $1.04$ within the optimal value of the Goemans and Williamson relaxation.

\begin{table}[H]
\caption{Ratio of $\chi(G)$ to $z_{sdp}$ and ratio  of $z_{SP_{\mathcal{S}}}$ to $z_{sdp}$ for $k=4$ and $l=5$.}\label{tab3}
\centering
\begin{tabular}{|p{1cm}|p{3cm}|p{3cm}|}
\hline
$ \mathbf{n}$  & $ \mathbf{\chi(G)/z_{sdp}}$\textbf{: average(sd)} & $ \mathbf{z_{\mathcal{S}}/z_{sdp}}$\textbf{: average(sd)}   \\
\hline
$64$  & $1.241\ (0.008)$  &  $1.020 \ (0.002)$   \\
\hline
$100$  & $1.417 \ (0.007)$ &  $1.017 \ (0.002)$  \\
\hline
$196$  & $1.760 \ (0.003)$   & $1.012 \ (0.001)$ \\
\hline 
$400$ & $2.289 \ (0.003)$  & $1.010 \ (0.001)$  \\
\hline
\end{tabular}
\label{table3}
\end{table}

\begin{table}[H]
\caption{Ratio of $\chi(G)$ to $Z_{sdp}$ and ratio  of $z_{\mathcal{S}}$ to $z_{sdp}$ for $k=6$ and $l=10$. }\label{tab4}
\centering
\begin{tabular}{|p{1cm}|p{3cm}|p{3cm}|}
\hline
$ \mathbf{n}$  & $ \mathbf{\chi(G)/z_{sdp}}$\textbf{: average(sd)} & $ \mathbf{z_{\mathcal{S}}/z_{sdp}}$\textbf{: average(sd)}   \\
\hline
$64$   & $1.137 \ (0.008)$ &$1.029\ (0.002)$ \\
\hline
$100$    & $1.278\ (0.007)$ &  $1.024 \ (0.001)$  \\
\hline
$196$    & $1.546 \ (0.005)$ &  $1.020\ (0.001)$  \\
\hline 
$400$ &   $1.962 \ (0.002)$  & $1.013 \ (0.001)$  \\
\hline
\end{tabular}
\label{table4}
\end{table}

\comment{

\begin{figure}\label{logpgraph}
\includegraphics[width=0.5\textwidth]{logn.jpg}[width=0.5]
\caption{A figure caption is always placed below the illustration.
Please note that short captions are centered, while long ones are
justified by the macro package automatically.} \label{fig1}
\end{figure}

\begin{figure}\label{p02}
\includegraphics[width=0.5\textwidth]{p02.jpg}
\caption{A figure caption is always placed below the illustration.
Please note that short captions are centered, while long ones are
justified by the macro package automatically.} \label{fig1}
\end{figure}

\begin{figure}\label{np05}
\includegraphics[width=0.5\textwidth]{np05.jpg}
\caption{A figure caption is always placed below the illustration.
Please note that short captions are centered, while long ones are
justified by the macro package automatically.} \label{fig1}
\end{figure}

\begin{table}\label{table1}
\caption{Quotient of $SPGW(\mathcal{S})$ and $DSDGW(\mathcal{S})$ for different values of $n$ and $p$ . }\label{tab1}
\centering
\begin{tabular}{|l|l|l|l|l|l|}
\hline
$n$ & $p=0.1$ & $p=0.3$ & $p=0.5$  & $p=0.7$ & $p=0.9$\\
\hline
50 &  &  & & & \\
\hline
100 &  &  & & & \\
\hline
200 &  &  & & & \\

\hline
\end{tabular}
\end{table}

\begin{table}\label{table2}
\caption{Quotient of $SPGW(\mathcal{S})$ and $DSDGW(\mathcal{S})$ for different values of $n$ and $p$ . }\label{tab2}
\centering
\begin{tabular}{|l|l|l|l|l|l|}
\hline
$n$ & $p=0.1$ & $p=0.3$ & $p=0.5$  & $p=0.7$ & $p=0.9$\\
\hline
50 &  &  & & & \\
\hline
100 &  &  & & & \\
\hline
200 &  &  & & & \\

\hline
\end{tabular}
\end{table}
}

\subsection{Lov\'asz theta number}

Denote by $z_{sdp}$, $z_n$ and $z_\mathcal{S}$ the objective values of programs \ref{theta}, $Iter_n(Tn)$ and \ref{lineartheta} with 
$\mathcal{S}$ chosen as in Subsection \ref{subsec:finding_conmuting}, respectively. For each $n$ ranging from $20$ to $200$ in steps of $10$, we generate $5$ random graphs 
and plot the maximum, minimum and median of the quotients $\frac{z_\mathcal{S}}{z_{sdp}}$ and $\frac{z_n}{z_{sdp}}$ for these five instances.
In the following subsections, we present these plots for   Erd\H{o}s-R\'enyi and random d-regular graphs.

\subsubsection{  Erd\H{o}s-R\'enyi random graphs}

\begin{figure}[htbp]
\centering
\begin{adjustbox}{max width=1.3\textwidth}
\begin{subfigure}[b]{0.6\linewidth}
\centering
\includegraphics[width=\linewidth]{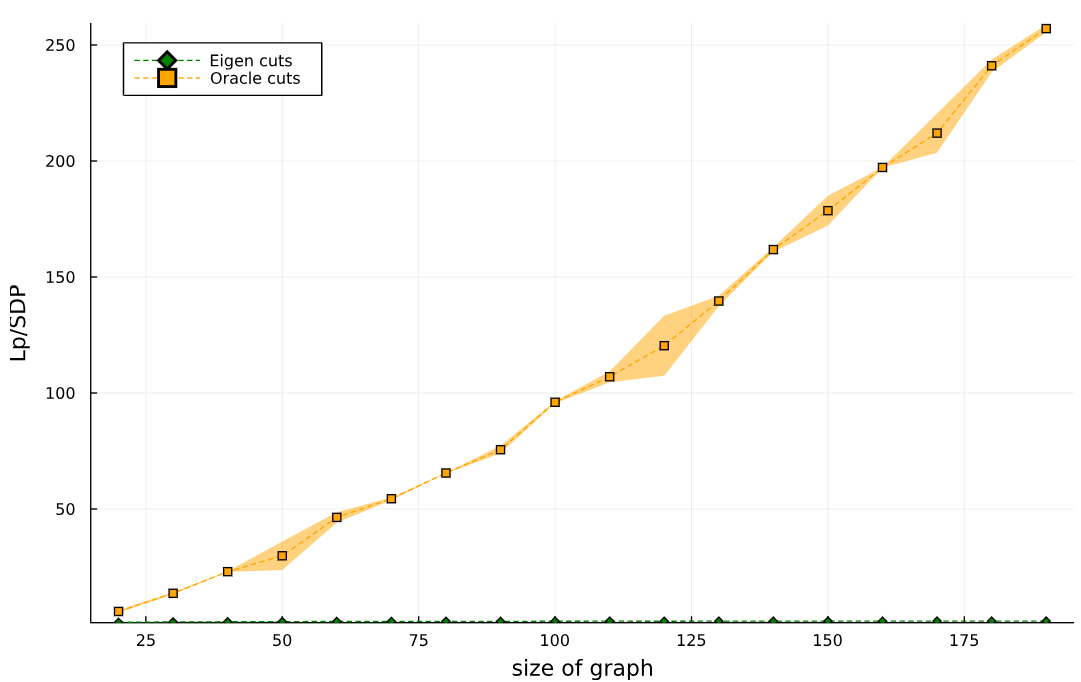}
\caption{$p=0.5$}
\end{subfigure}%
\begin{subfigure}[b]{0.6\linewidth}
\centering
\includegraphics[width=\linewidth]{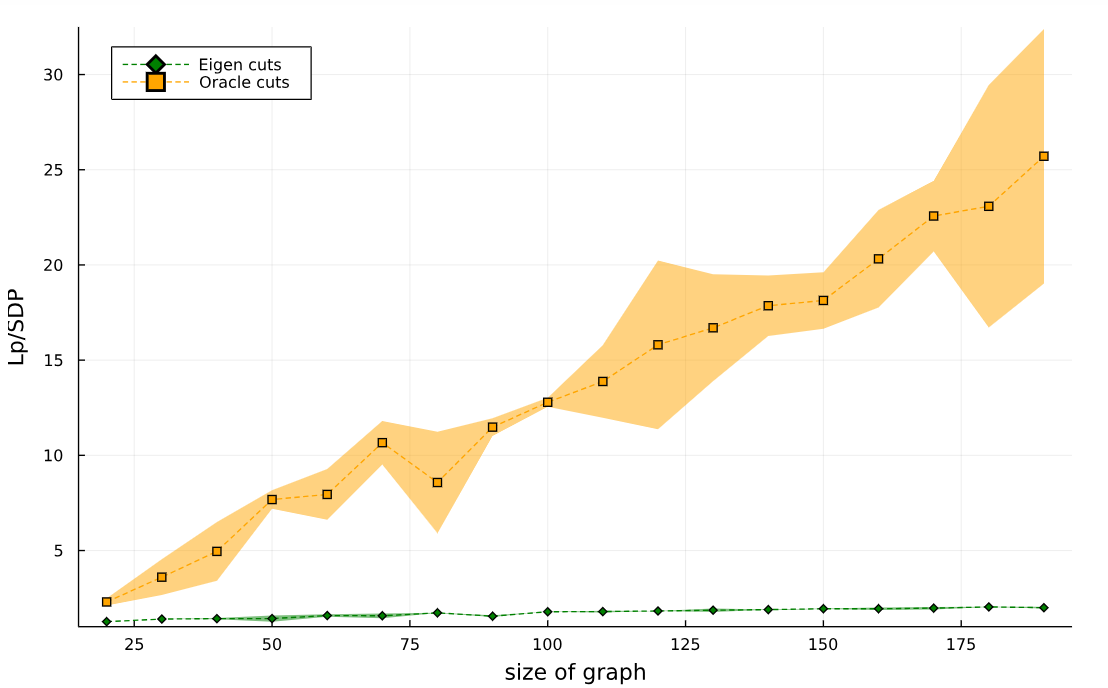}
\caption{$p=0.9$}
\end{subfigure}
\end{adjustbox}

\begin{subfigure}[b]{0.6\linewidth}
\centering
\includegraphics[width=\linewidth]{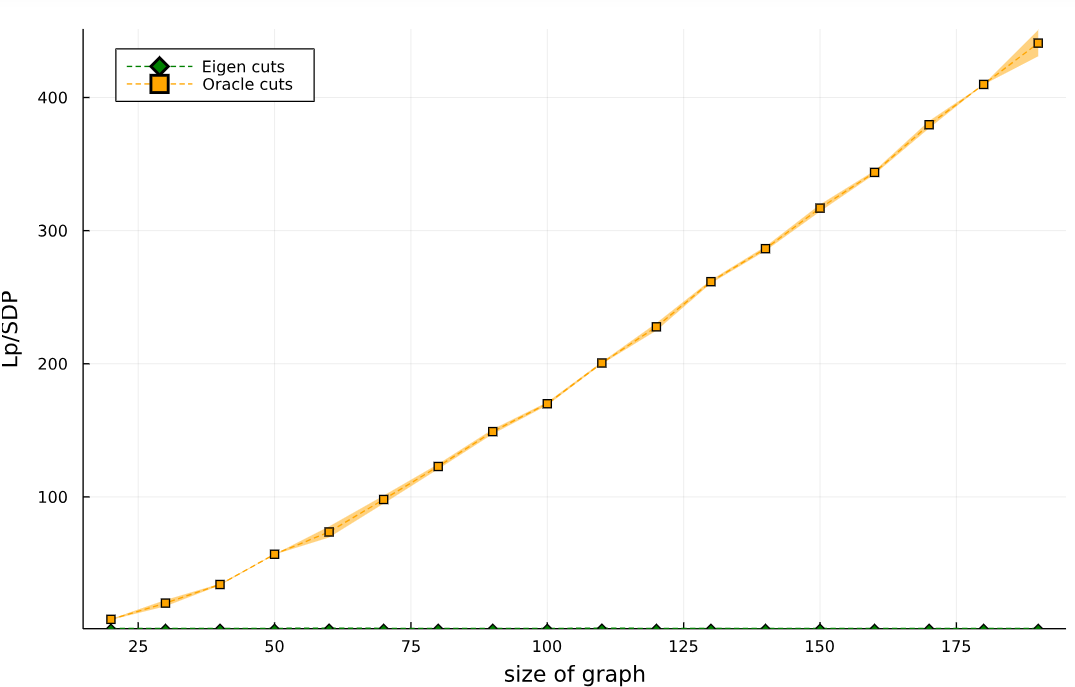}
\caption{$p=\frac{3\log(n)}{n}$}
\end{subfigure}

\caption{Quotients for the Lov\'asz theta number $\frac{z_\mathcal{S}}{z_{sdp}}$ (Eigen cuts) and $\frac{z_n}{z_{sdp}}$ (Oracle cuts) as $n$ grows for   Erd\H{o}s-R\'enyi random graphs with different values of $p$.}
\label{fig:theta_Erdos}
\end{figure}

In Figure \ref{fig:theta_Erdos} We plot the mentioned quotients for   Erd\H{o}s-R\'enyi random graph while we vary $p$, the probability of connecting two edges.

\comment{

\begin{figure}[h]
    \centering
    \includegraphics[width=0.9\textwidth]{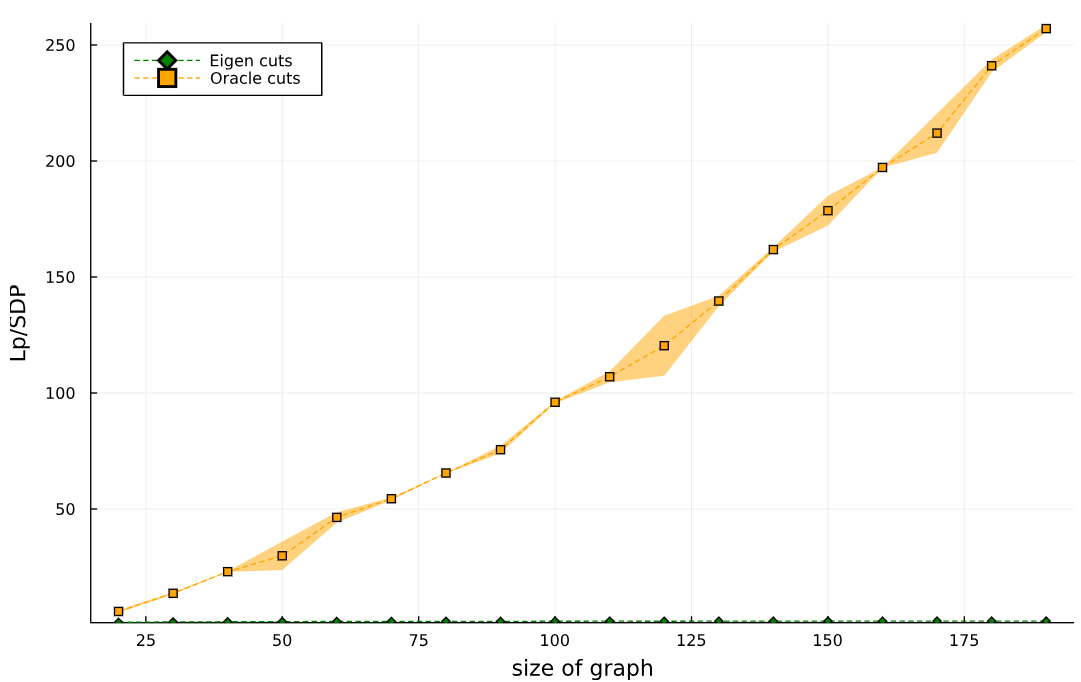}

    \caption{ Quotients $\frac{z_\mathcal{S}}{z_{sdp}}$ (Eigen cuts) and $\frac{z_n}{z_{sdp}}$ (Oracle cuts) as $n$ grows for   Erd\H{o}s-R\'enyi random graphs with $p=0.5$ for the Lov\'asz theta number.}
    \label{fig:theta_Erdos_05}
\end{figure}

\begin{figure}[h]
    \centering
    \includegraphics[width=0.9\textwidth]{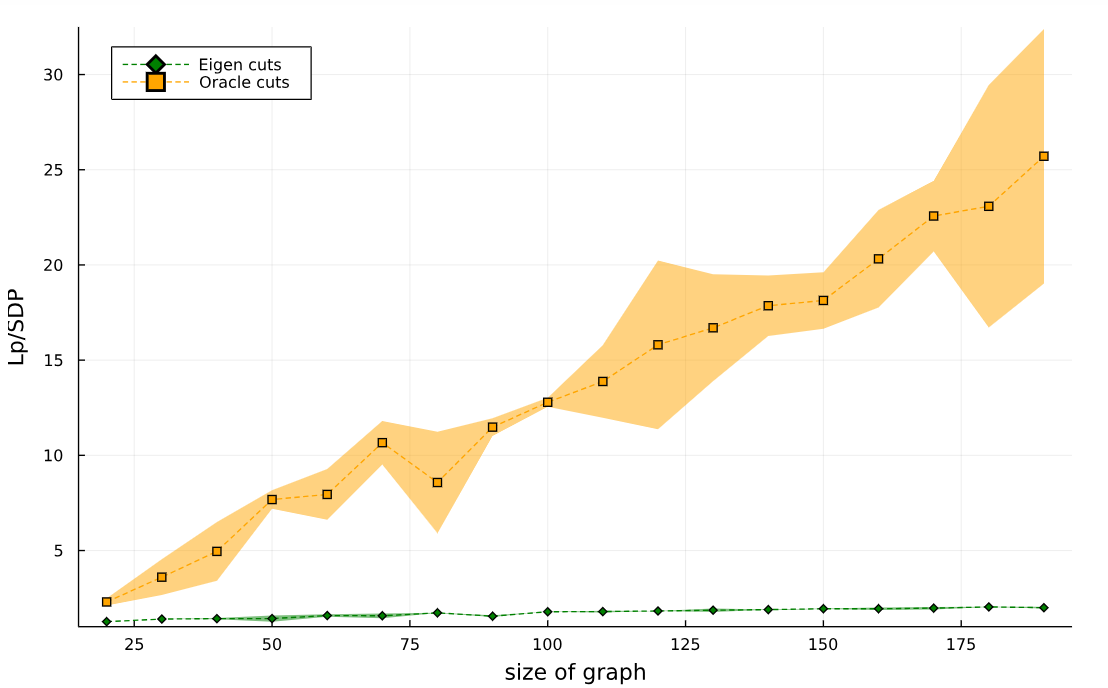}
    \caption{ Quotients $\frac{z_\mathcal{S}}{z_{sdp}}$ (Eigen cuts) and $\frac{z_n}{z_{sdp}}$ (Oracle cuts) as $n$ grows for   Erd\H{o}s-R\'enyi random graphs with $p=0.9$ for the Lov\'asz theta number.}
    \label{fig:theta_Erdos_09}
\end{figure}

\begin{figure}[h]
    \centering
    \includegraphics[width=0.9\textwidth]{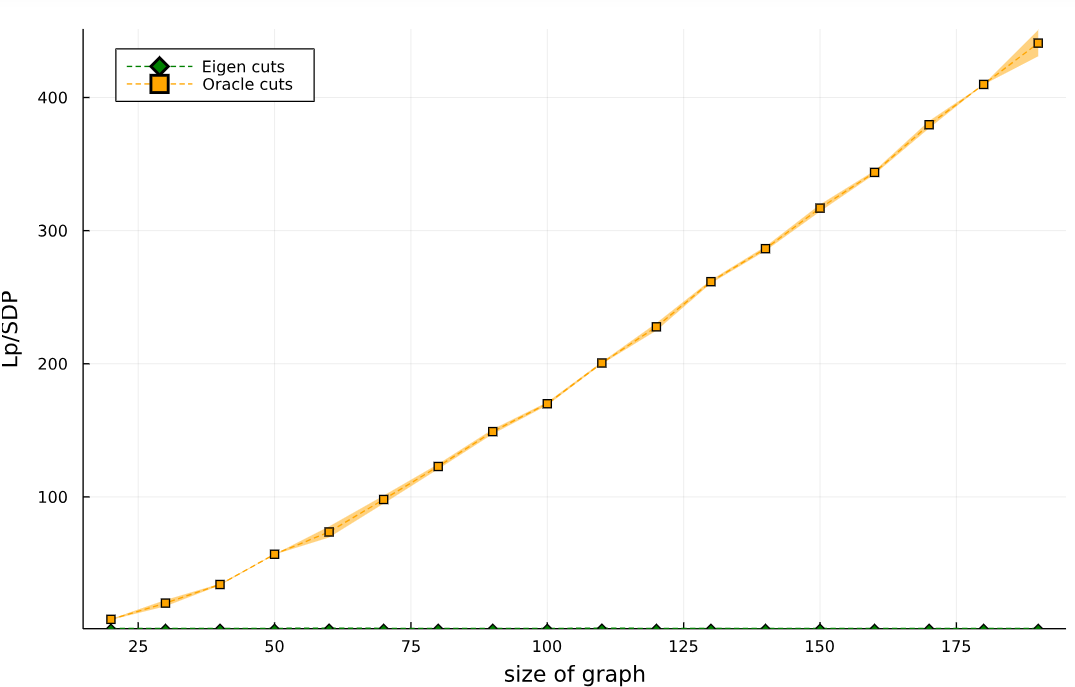}
    \caption{Quotients $\frac{z_\mathcal{S}}{z_{sdp}}$ (Eigen cuts) and $\frac{z_n}{z_{sdp}}$ (Oracle cuts) as $n$ grows for   Erd\H{o}s-R\'enyi random graphs with $p=\frac{3\log(n)}{n}$, for the Lov\'asz theta number.
}
    \label{fig:theta_Erdos_log}
\end{figure}

 }

\subsubsection{d-regular random graphs}

\begin{figure}[htbp]
\centering
\begin{adjustbox}{max width=1.3\textwidth}
\begin{subfigure}[b]{0.6\linewidth}
\centering
\includegraphics[width=\linewidth]{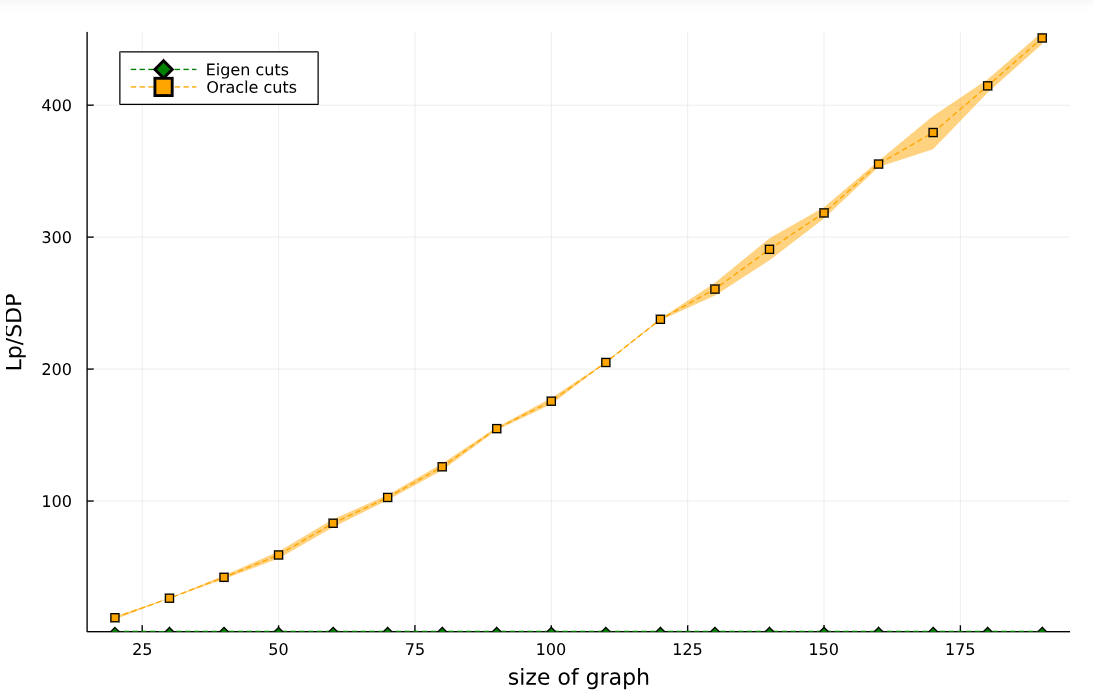}
\caption{$d=5$}
\end{subfigure}%
\begin{subfigure}[b]{0.6\linewidth}
\centering
\includegraphics[width=\linewidth]{revised_pictures/lovas_reg_sqrt.png}
\caption{$d=\sqrt{n}$}
\end{subfigure}
\end{adjustbox}

\begin{subfigure}[b]{0.6\linewidth}
\centering
\includegraphics[width=\linewidth]{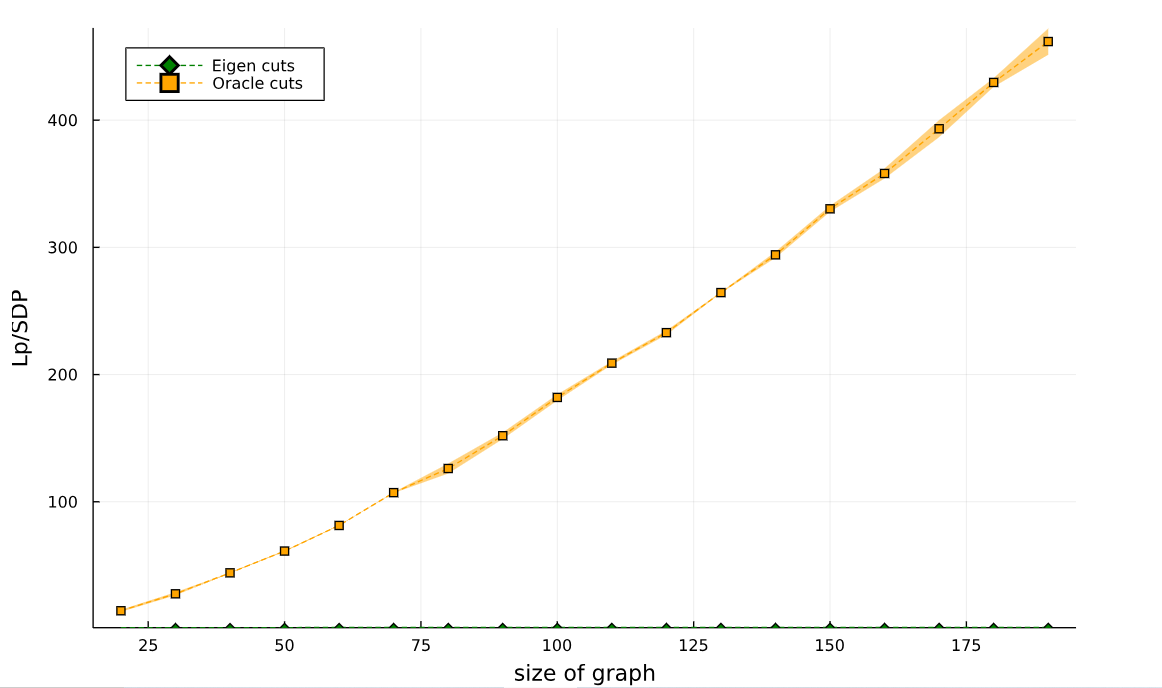}
\caption{$d=\frac{n}{10}$}
\end{subfigure}

\caption{Quotients for the Lov\'asz theta number $\frac{z_\mathcal{S}}{z_{sdp}}$ and $\frac{z_n}{z_{sdp}}$ as $n$ grows for random $d-$regular graphs with different values of $d$, as $n$ grows.}
\label{fig:theta_regular}
\end{figure}

In Figure \ref{fig:theta_regular} we plot the mentioned quotients for $d-$regular random graph while we vary $d$.

\comment{

\begin{figure}[h]
    \centering
    \includegraphics[width=0.9\textwidth]{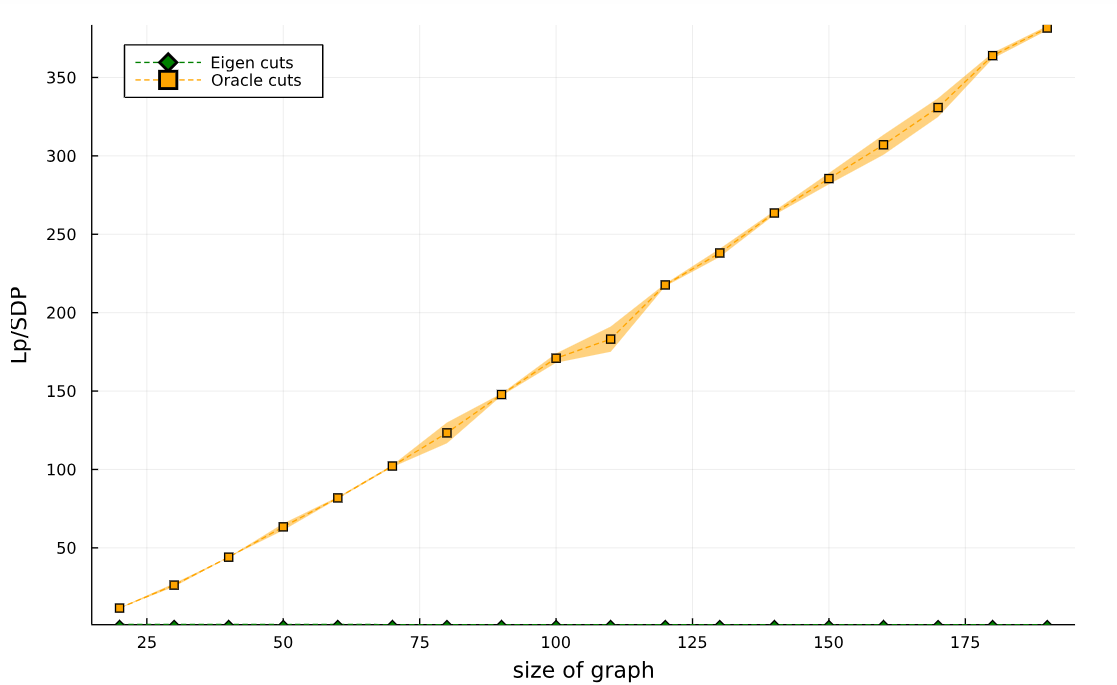}
    \caption{Quotients $\frac{z_\mathcal{S}}{z_{sdp}}$ (Eigen cuts) and $\frac{z_n}{z_{sdp}}$ (Oracle cuts) as $n$ grows for random $d-$regular graphs with $d=5$, as $n$ grows, for the Lov\'asz theta number.}
    \label{fig:theta_reg_5}
\end{figure}

\begin{figure}[h]
    \centering
    \includegraphics[width=0.9\textwidth]{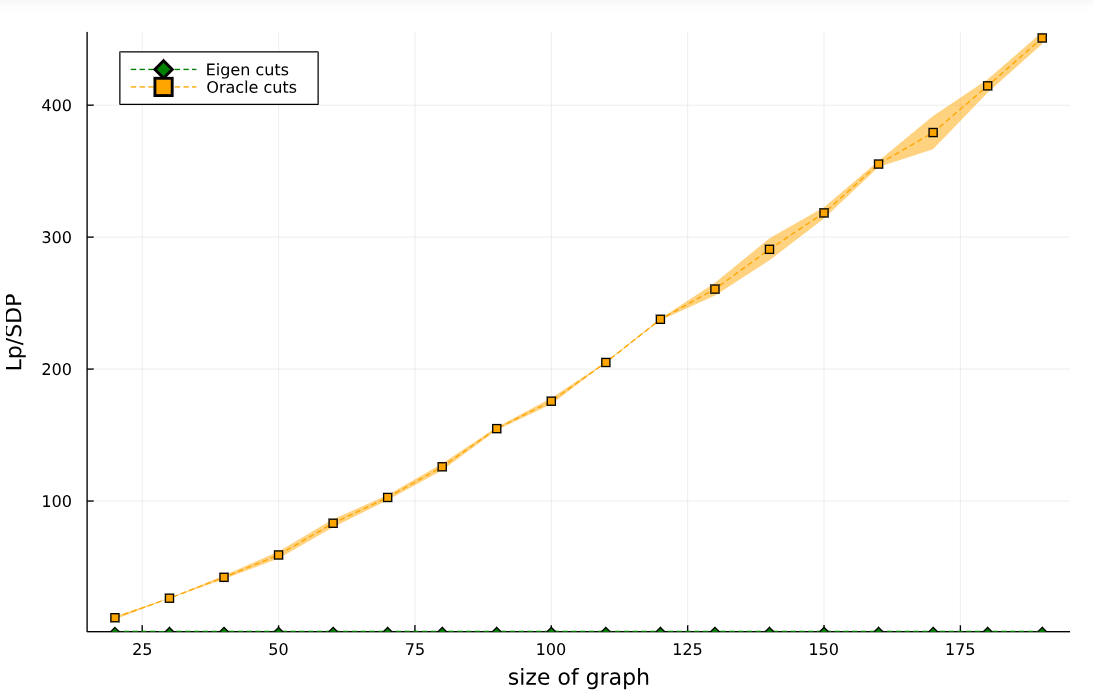}
     \caption{Quotients $\frac{z_\mathcal{S}}{z_{sdp}}$ (Eigen cuts) and $\frac{z_n}{z_{sdp}}$ (Oracle cuts) as $n$ grows for random $d-$regular graphs with $d=\sqrt{n}$, as $n$ grows, for the Lov\'asz theta number.}
    \label{fig:theta_reg_sqrt}
\end{figure}

\begin{figure}[h]
    \centering
    \includegraphics[width=0.9\textwidth]{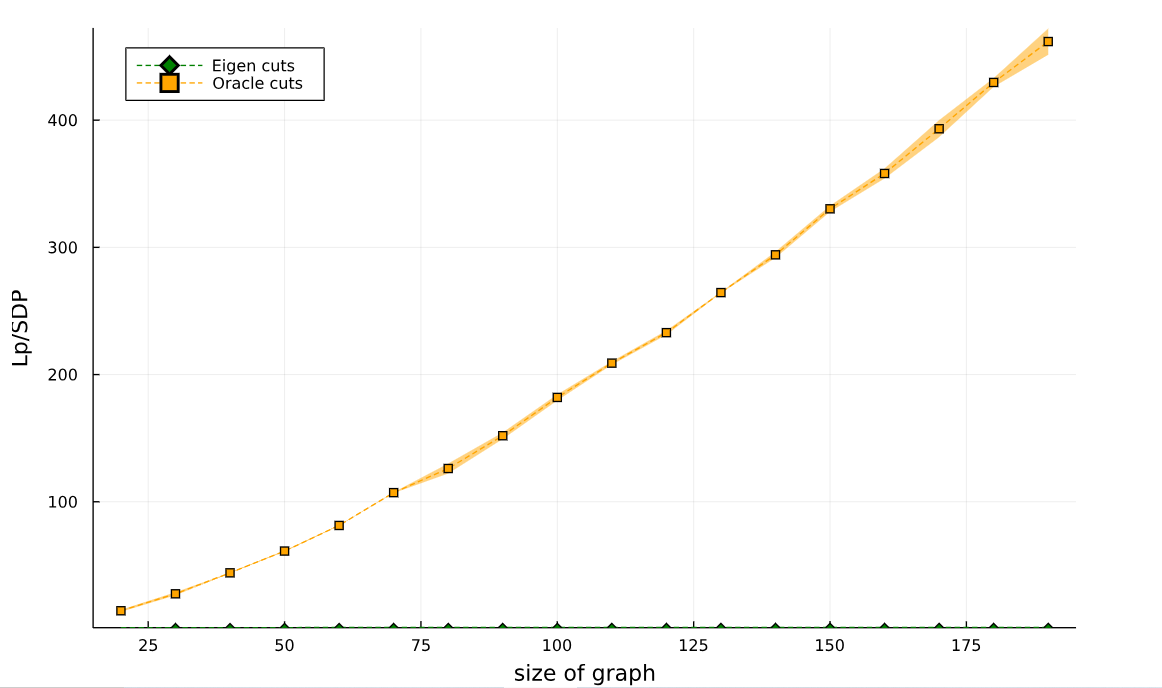}
     \caption{Quotients $\frac{z_\mathcal{S}}{z_{sdp}}$ and $\frac{z_n}{z_{sdp}}$ as $n$ grows for random $d-$regular graphs with $d=\frac{n}{10}$, as $n$ grows, for the Lov\'asz theta number.}
    \label{fig:theta_reg_linear}
\end{figure}

}

\subsection{Quadratically constrained quadratic problems}

In this subsection we test the proposed methodology on the different QCQPs introduced in Section \ref{revised_sec:4}.

\subsubsection{Random QCQPs}

We generate random QCQPs following the review \cite{bao2011semidefinite}, where the authors compare various SDP relaxations of QCQPs in terms of percentage distance to the objective and solution time. For these instances, the $x$ variables are bounded in an unit box $[0,1]^n$ and the number of variables is varied from $20$ up to $100$ in steps of $10$. The vectors $c,d$ in $\mathbb{R}^{r+1}$ and $\mathbb{R}^q$ respectively and the matrices $D \in \mathbb{R}^{q\times n}$ and $A_i \in \mathbb{S}^n$, $i\in \{0,\dots,r\}$ have entries drawn uniformly and independently at random from an uniform distribution supported in $[-1,1]$. The vector $b\in \mathbb{R}^{r+1}$ has entries sampled uniformly at random from an uniform distribution supported in $[0,100]$.
Since QCPQs are highly sensitive to the number of quadratic constraints, we test different combinations of number of quadratic and linear constraints, according to the following combinations:

\begin{itemize}
    \item QCQPs with $r=1, q = \frac{n}{10}.$
    \item  QCQPs with $r=1, q = \frac{n}{5}.$
    \item  QCQPs with $r=\frac{n}{2}, q = \frac{n}{10}$
    \item  QCQPs with $r=n, q = \frac{n}{10}.$
\end{itemize}

Furthermore, we consider different densities $\Delta$ for the matrix $A_0$,  which corresponds to the percentage of nonzero elements of the matrix, on average. For a given combination of these parameters and a value of $n$ we generate $5$ random instances and solve the following optimization programs for each:

\begin{itemize}
    \item Problem \ref{shor_reformulated}. We denote the objective value of this semidefinite program by $z_{sdp}$.
    \item The linear relaxation $L_{\mathcal{S}}$ of \ref{shor_reformulated} where we let $\mathcal{S}$ the elements of a eigenvector basis of the matrix $A_0$. We denote the objective value of this problem by $z_\mathcal{S}$.
    \item The LP $Iter_n$(\ref{shor_reformulated}). We denote by $z_n$ the objective value of this program.
    \item The LP $Iter_0$(\ref{shor_reformulated}). We denote by $z_0$ the objective value of this program.
    \item The second order cone program obtained by dropping the constraint $\Hat{X}\succeq 0$ from  \ref{shor_reformulated} adding the constraints (\ref{eq:soc}) with $\mathcal{S}$ the elements of a eigenvector basis of the matrix $A_0$. We denote the objective value of this problem by $z_{soc}$.
\end{itemize}

\begin{figure}[htbp]
\centering
\begin{adjustbox}{max width=1.3\textwidth}
\begin{subfigure}[b]{0.6\linewidth}
\centering
\includegraphics[width=\linewidth]{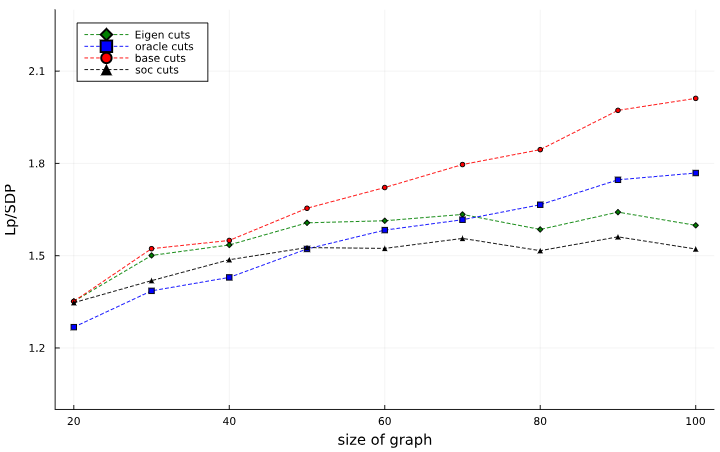}
\caption{$r=1, \ q= \frac{n}{5}$}
\end{subfigure}%
\begin{subfigure}[b]{0.6\linewidth}
\centering
\includegraphics[width=\linewidth]{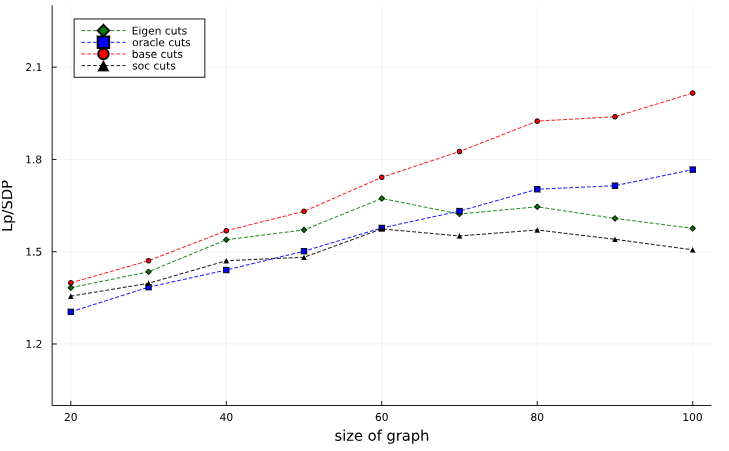}
\caption{$r=1, \ q= \frac{n}{10}$}
\end{subfigure}
\end{adjustbox}

\vspace{1cm} 

\begin{adjustbox}{max width=1.3\textwidth}
\begin{subfigure}[b]{0.6\linewidth}
\centering
\includegraphics[width=\linewidth]{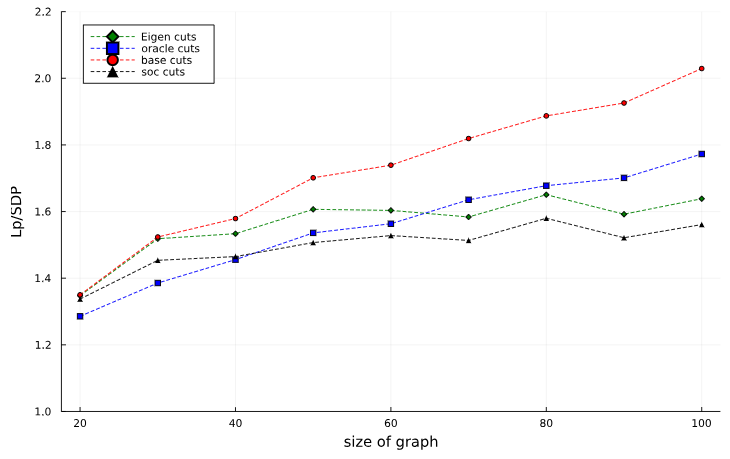}
\caption{$r=\frac{n}{2}, \ q= \frac{n}{10}$}
\end{subfigure}%
\begin{subfigure}[b]{0.6\linewidth}
\centering
\includegraphics[width=\linewidth]{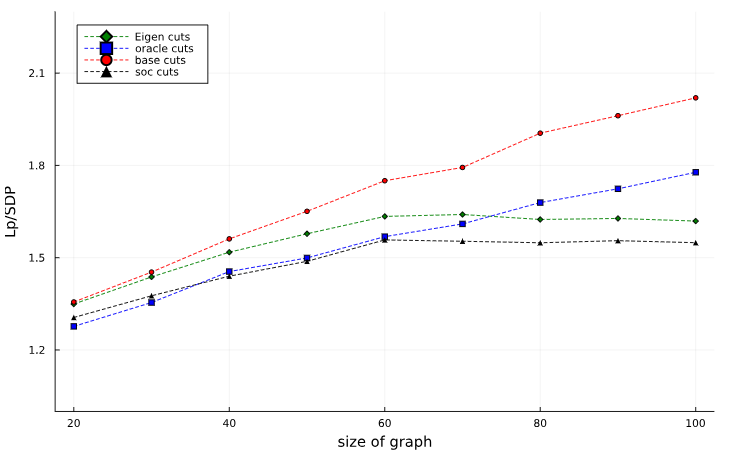}
\caption{$r=n, \ q= \frac{n}{10}$}
\end{subfigure}
\end{adjustbox}
\caption{Quality of the ratios $\frac{z_\mathcal{S}}{z_{sdp}}$ (eigen cuts), $\frac{z_n}{z_{sdp}}$, (oracle cuts), $\frac{Z_0}{z_{sdp}}$ (base cuts) and $\frac{z_{soc}}{z_{sdp}}$ for random QCQP instances with density $0.25$. }
\label{fig:qcqp_025}
\end{figure}

\begin{figure}[htbp]
\centering
\begin{adjustbox}{max width=1.3\textwidth}
\begin{subfigure}[b]{0.6\linewidth}
\centering
\includegraphics[width=\linewidth]{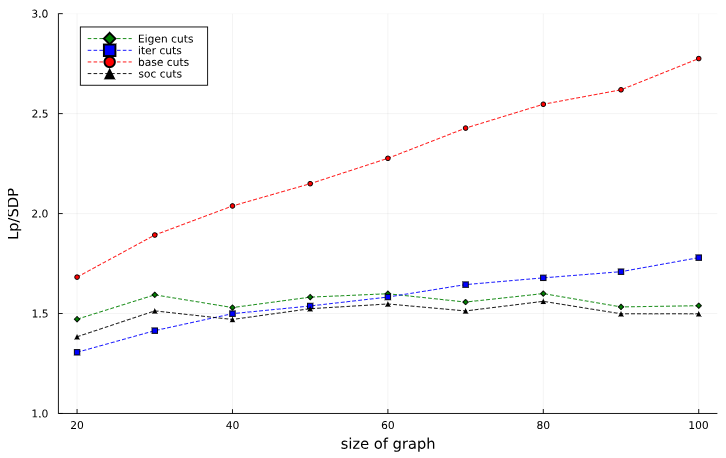}
\caption{$r=1, \ q= \frac{n}{5}$}
\end{subfigure}%
\begin{subfigure}[b]{0.6\linewidth}
\centering
\includegraphics[width=\linewidth]{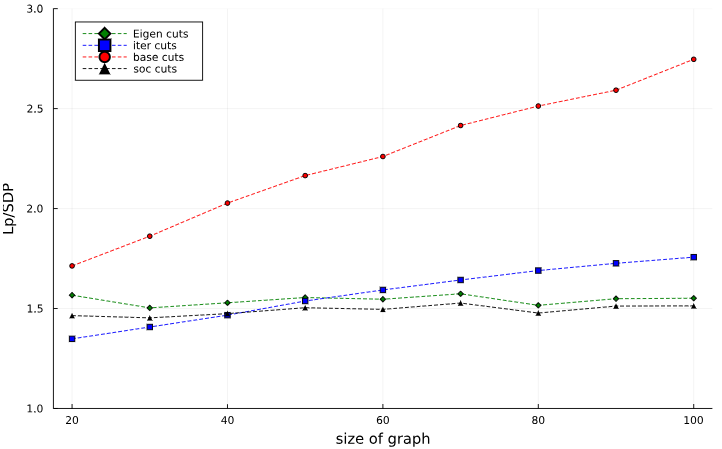}
\caption{$r=1, \ q= \frac{n}{10}$}
\end{subfigure}
\end{adjustbox}

\vspace{1cm} 

\begin{adjustbox}{max width=1.3\textwidth}
\begin{subfigure}[b]{0.6\linewidth}
\centering
\includegraphics[width=\linewidth]{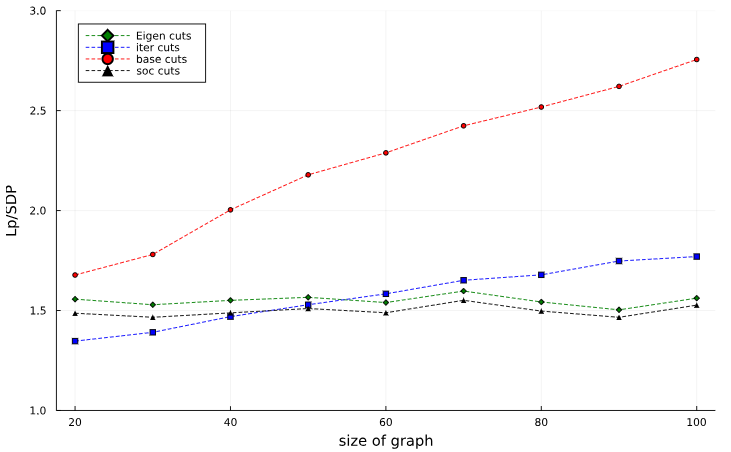}
\caption{$r=\frac{n}{2}, \ q= \frac{n}{10}$}
\end{subfigure}%
\begin{subfigure}[b]{0.6\linewidth}
\centering
\includegraphics[width=\linewidth]{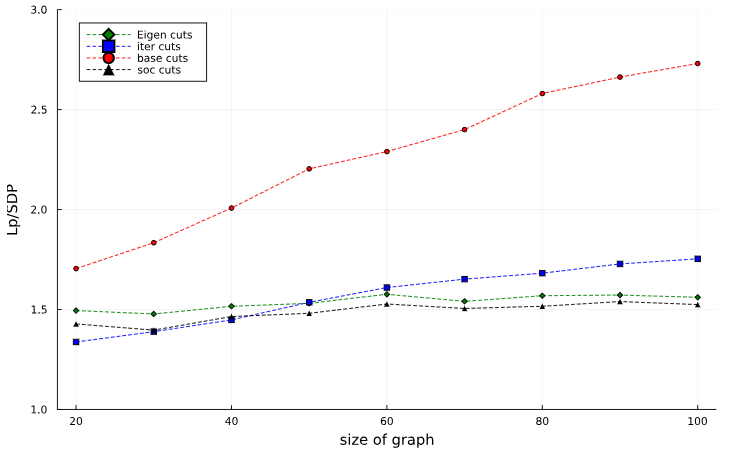}
\caption{$r=n, \ q= \frac{n}{10}$}
\end{subfigure}
\end{adjustbox}
\caption{Quality of the ratios $\frac{z_\mathcal{S}}{z_{sdp}}$ (eigen cuts), $\frac{z_n}{z_{sdp}}$, (oracle cuts), $\frac{z_0}{z_{sdp}}$ (base cuts) and $\frac{z_{soc}}{z_{sdp}}$ for random QCQP instances with density $0.5$. }
\label{fig:qcqp_05}
\end{figure}

\begin{figure}[htbp]
\centering
\begin{adjustbox}{max width=1.3\textwidth}
\begin{subfigure}[b]{0.6\linewidth}
\centering
\includegraphics[width=\linewidth]{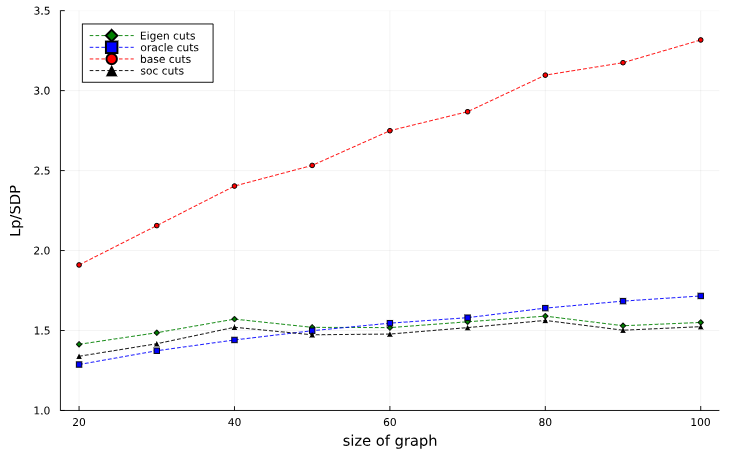}
\caption{$r=1, \ q= \frac{n}{5}$}
\end{subfigure}%
\begin{subfigure}[b]{0.6\linewidth}
\centering
\includegraphics[width=\linewidth]{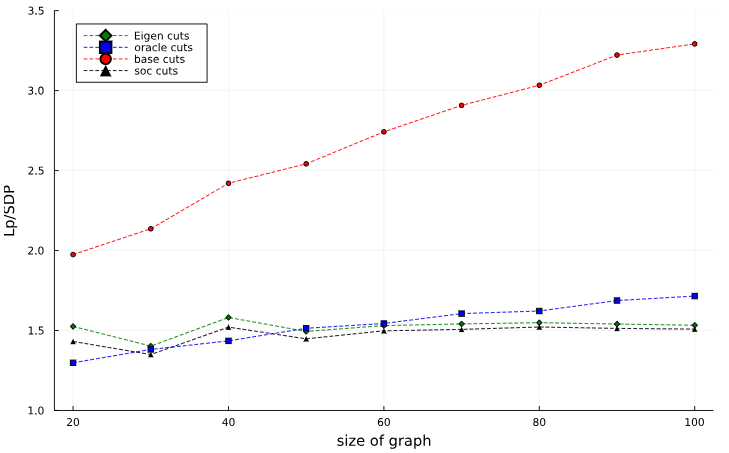}
\caption{$r=1, \ q= \frac{n}{10}$}
\end{subfigure}
\end{adjustbox}

\vspace{1cm} 

\begin{adjustbox}{max width=1.3\textwidth}
\begin{subfigure}[b]{0.6\linewidth}
\centering
\includegraphics[width=\linewidth]{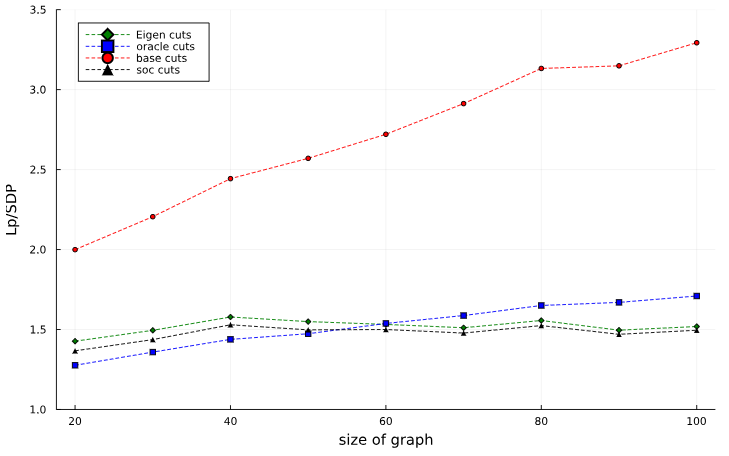}
\caption{$r=\frac{n}{2}, \ q= \frac{n}{10}$}
\end{subfigure}%
\begin{subfigure}[b]{0.6\linewidth}
\centering
\includegraphics[width=\linewidth]{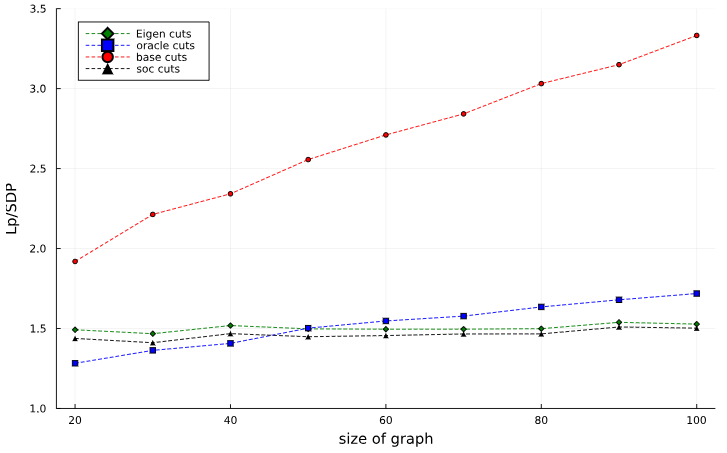}
\caption{$r=n, \ q= \frac{n}{10}$}
\end{subfigure}
\end{adjustbox}
\caption{Quality of the ratios $\frac{z_\mathcal{S}}{z_{sdp}}$ (eigen cuts), $\frac{z_n}{z_{sdp}}$, (oracle cuts), $\frac{Z_0}{z_{sdp}}$ (base cuts) and $\frac{z_{soc}}{z_{sdp}}$ for random QCQP instances with density $0.75$. }  
\label{fig:qcqp_075}
\end{figure}

\begin{figure}[htbp]
\centering
\begin{adjustbox}{max width=1.3\textwidth}
\begin{subfigure}[b]{0.6\linewidth}
\centering
\includegraphics[width=\linewidth]{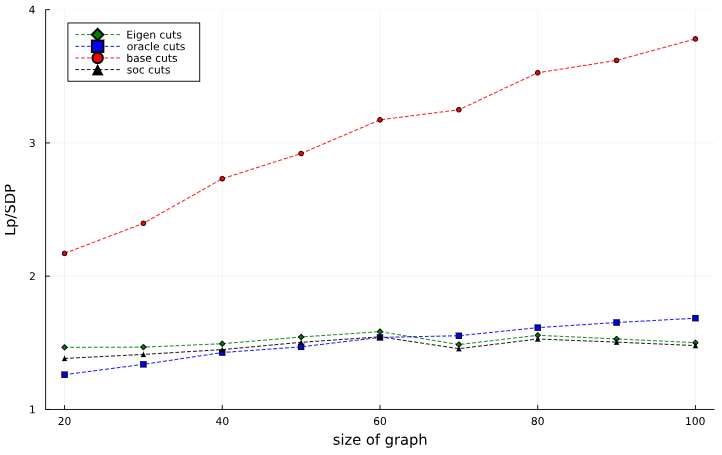}
\caption{$r=1, \ q= \frac{n}{5}$}
\end{subfigure}%
\begin{subfigure}[b]{0.6\linewidth}
\centering
\includegraphics[width=\linewidth]{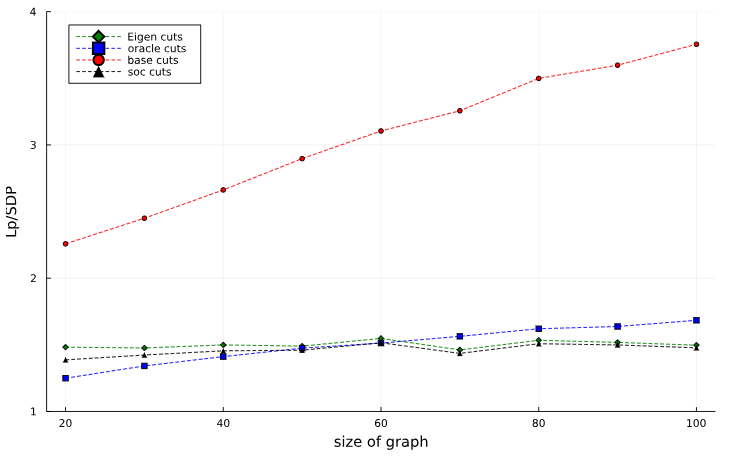}
\caption{$r=1, \ q= \frac{n}{10}$}
\end{subfigure}
\end{adjustbox}

\vspace{1cm} 

\begin{adjustbox}{max width=1.3\textwidth}
\begin{subfigure}[b]{0.6\linewidth}
\centering
\includegraphics[width=\linewidth]{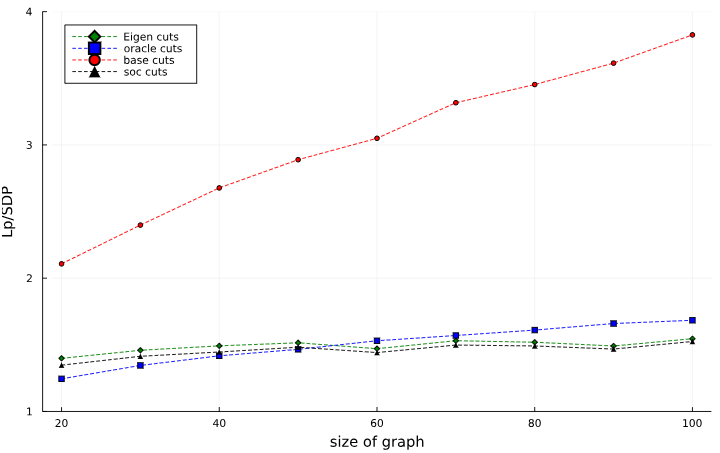}
\caption{$r=\frac{n}{2}, \ q= \frac{n}{10}$}
\end{subfigure}%
\begin{subfigure}[b]{0.6\linewidth}
\centering
\includegraphics[width=\linewidth]{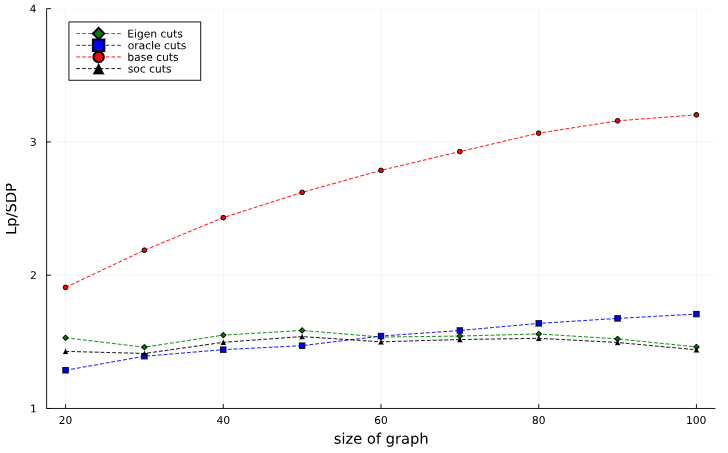}
\caption{$r=n, \ q= \frac{n}{10}$}
\end{subfigure}
\end{adjustbox}
\caption{Quality of the ratios $\frac{z_\mathcal{S}}{z_{sdp}}$ (eigen cuts), $\frac{z_n}{z_{sdp}}$, (oracle cuts), $\frac{Z_0}{z_{sdp}}$ (base cuts) and $\frac{z_{soc}}{z_{sdp}}$ for random QCQP instances with density $1$. } 
\label{fig:qcqp_1}
\end{figure}

We average the values of ratios $\frac{z_{\mathcal{S}}}{z_{sdp}}$,  $\frac{z_n}{z_{sdp}}$,  $\frac{z_0}{z_{sdp}}$ and $\frac{z_{soc}}{z_{sdp}}$  over the five instances, and plot the results in Figures 
\ref{fig:qcqp_025}, \ref{fig:qcqp_05}, \ref{fig:qcqp_075} and \ref{fig:qcqp_1}.
We observe that due to randomness, it will not be possible to linearly combine the matrices $A_1,\dots,A_m,A_{m+1}$ so that they commute with the objective matrix $A_0$. Therefore, program \ref{conmute_generator} will typically return the $0$ matrix, and $\mathcal{S}$ will simply be a basis of eigenvectors of $A_0$.

For these instances, the quality of all the relaxations is encouraging, with the trivial relaxation obtained by dropping the semidefinite constraint getting a ratio of at most $4$ in all of our experiments. The second order cone relaxation is typically the better as soon as $n$ exceeds $50$. Whenever the density increases, we notice that the ratios $\frac{z_\mathcal{S}}{z_{sdp}}$ and $\frac{z_{soc}}{z_{sdp}}$ get closer and closer, hinting at that the second order cone relaxation is not much stronger than the linear relaxation. Although the LP $Iter_n$ achieves a better ratio for small $n$, this is no longer true for larger values of $n$. In addition, notice that for a value of $n$ this LP requires solving $n$ LPs and $n$ eigenvector decompositions.

\subsubsection{Extended trust region problems}

\begin{figure}[htbp]
\centering
\begin{adjustbox}{max width=1.3\textwidth}
\begin{subfigure}[b]{0.6\linewidth}
\centering
\includegraphics[width=\linewidth]{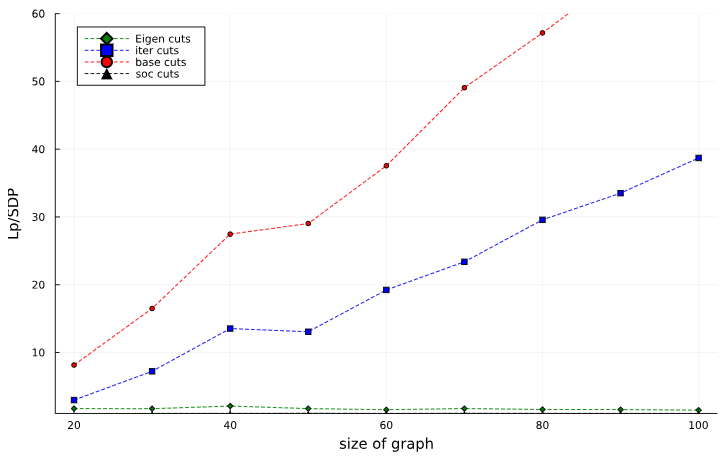}
\caption{$r=1, \ q= \frac{n}{5}$}
\end{subfigure}%
\begin{subfigure}[b]{0.6\linewidth}
\centering
\includegraphics[width=\linewidth]{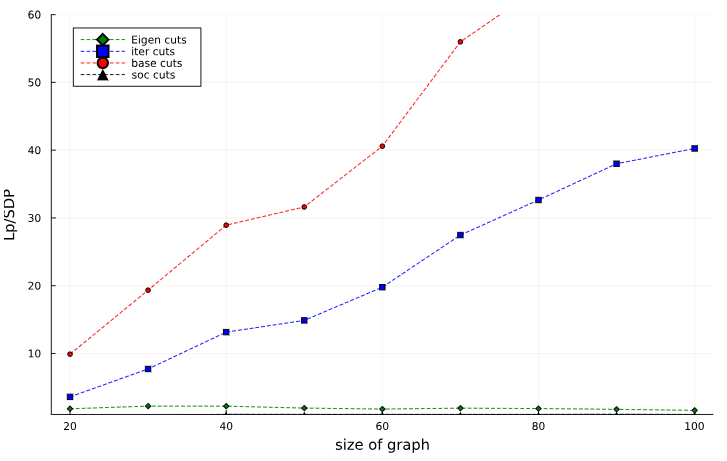}
\caption{$r=1, \ q= \frac{n}{10}$}
\end{subfigure}
\end{adjustbox}

\vspace{1cm} 

\begin{adjustbox}{max width=1.3\textwidth}
\begin{subfigure}[b]{0.6\linewidth}
\centering
\includegraphics[width=\linewidth]{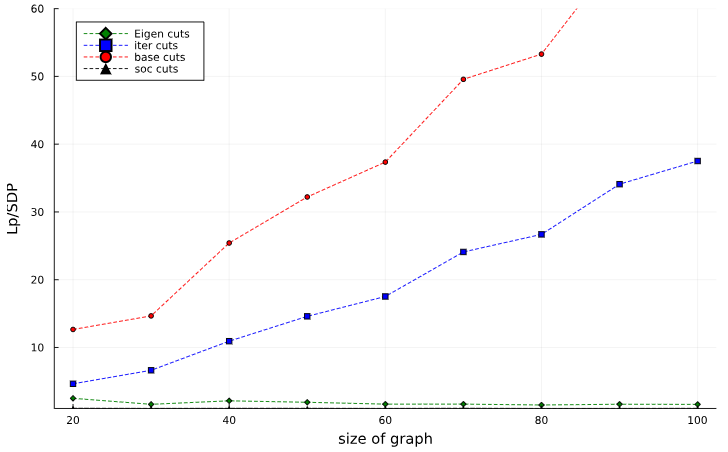}
\caption{$r=\frac{n}{2}, \ q= \frac{n}{10}$}
\end{subfigure}%
\begin{subfigure}[b]{0.6\linewidth}
\centering
\includegraphics[width=\linewidth]{revised_pictures/trust_pictures/trust_025den_rn_qn10.PNG}
\caption{$r=n, \ q= \frac{n}{10}$}
\end{subfigure}
\end{adjustbox}
\caption{Quality of the ratios $\frac{z_\mathcal{S}}{z_{sdp}}$ (eigen cuts), $\frac{z_n}{z_{sdp}}$, (oracle cuts), $\frac{Z_0}{z_{sdp}}$ (base cuts) and $\frac{z_{soc}}{z_{sdp}}$ for instances of the extended trust region problem with density $0.25$. }
\label{fig:trust_025}
\end{figure}

\begin{figure}[htbp]
\centering
\begin{adjustbox}{max width=1.3\textwidth}
\begin{subfigure}[b]{0.6\linewidth}
\centering
\includegraphics[width=\linewidth]{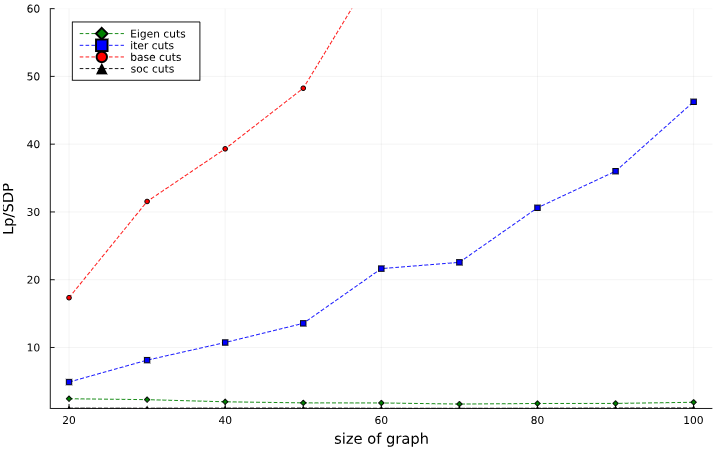}
\caption{$r=1, \ q= \frac{n}{5}$}
\end{subfigure}%
\begin{subfigure}[b]{0.6\linewidth}
\centering
\includegraphics[width=\linewidth]{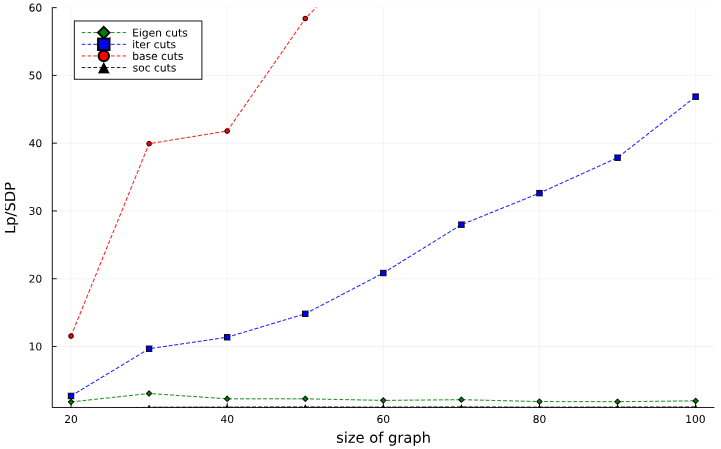}
\caption{$r=1, \ q= \frac{n}{10}$}
\end{subfigure}
\end{adjustbox}

\vspace{1cm} 

\begin{adjustbox}{max width=1.3\textwidth}
\begin{subfigure}[b]{0.6\linewidth}
\centering
\includegraphics[width=\linewidth]{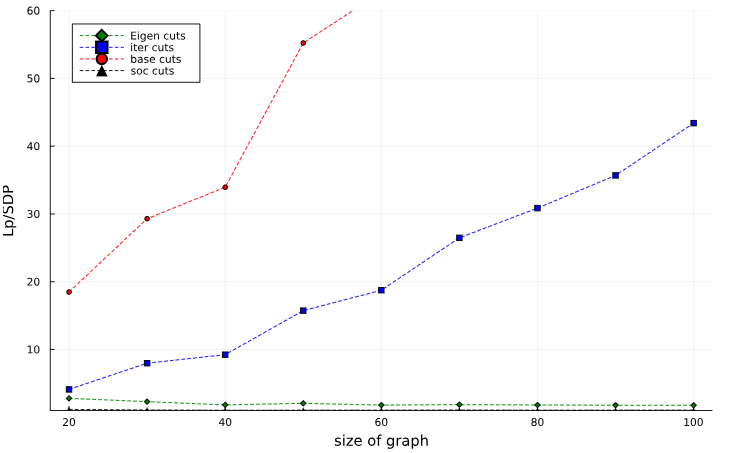}
\caption{$r=\frac{n}{2}, \ q= \frac{n}{10}$}
\end{subfigure}%
\begin{subfigure}[b]{0.6\linewidth}
\centering
\includegraphics[width=\linewidth]{revised_pictures/trust_pictures/trust_05den_rn_qn10.PNG}
\caption{$r=n, \ q= \frac{n}{10}$}
\end{subfigure}
\end{adjustbox}
\caption{Quality of the ratios $\frac{z_\mathcal{S}}{z_{sdp}}$ (eigen cuts), $\frac{z_n}{z_{sdp}}$, (oracle cuts), $\frac{Z_0}{z_{sdp}}$ (base cuts) and $\frac{z_{soc}}{z_{sdp}}$ for instances of the extended trust region problem with density $0.5$. }
\label{fig:trust_05}
\end{figure}

\begin{figure}[htbp]
\centering
\begin{adjustbox}{max width=1.3\textwidth}
\begin{subfigure}[b]{0.6\linewidth}
\centering
\includegraphics[width=\linewidth]{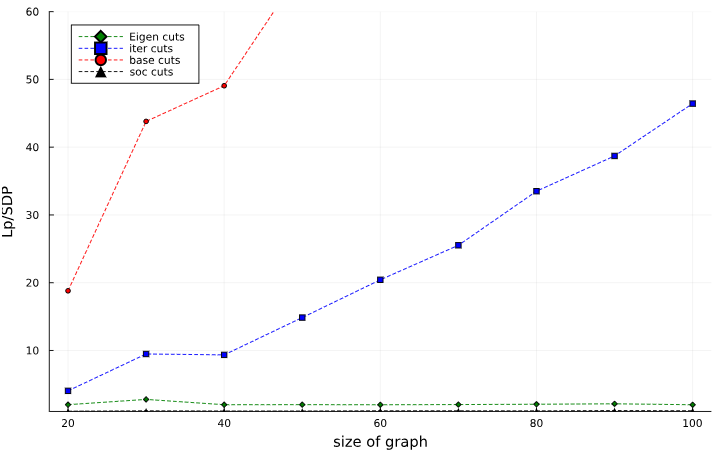}
\caption{$r=1, \ q= \frac{n}{5}$}
\end{subfigure}%
\begin{subfigure}[b]{0.6\linewidth}
\centering
\includegraphics[width=\linewidth]{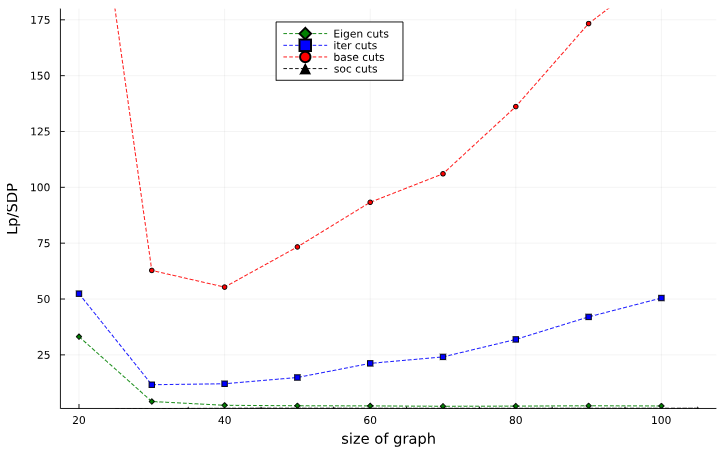}
\caption{$r=1, \ q= \frac{n}{10}$}
\end{subfigure}
\end{adjustbox}

\vspace{1cm} 

\begin{adjustbox}{max width=1.3\textwidth}
\begin{subfigure}[b]{0.6\linewidth}
\centering
\includegraphics[width=\linewidth]{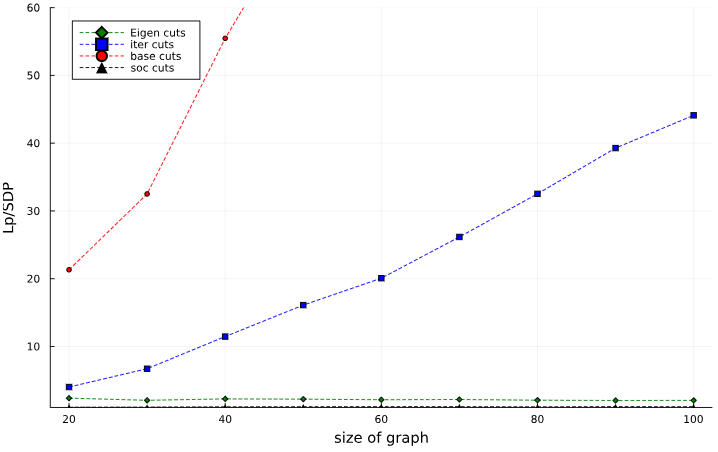}
\caption{$r=\frac{n}{2}, \ q= \frac{n}{10}$}
\end{subfigure}%
\begin{subfigure}[b]{0.6\linewidth}
\centering
\includegraphics[width=\linewidth]{revised_pictures/trust_pictures/trust_075den_rn_qn10.PNG}
\caption{$r=n, \ q= \frac{n}{10}$}
\end{subfigure}
\end{adjustbox}
\caption{Quality of the ratios $\frac{z_\mathcal{S}}{z_{sdp}}$ (eigen cuts), $\frac{z_n}{z_{sdp}}$, (oracle cuts), $\frac{z_0}{z_{sdp}}$ (base cuts) and $\frac{z_{soc}}{z_{sdp}}$ for instances of the extended trust region problem with density $0.75$. }
\label{fig:trust_075}
\end{figure}

\begin{figure}[htbp]
\centering
\begin{adjustbox}{max width=1.3\textwidth}
\begin{subfigure}[b]{0.6\linewidth}
\centering
\includegraphics[width=\linewidth]{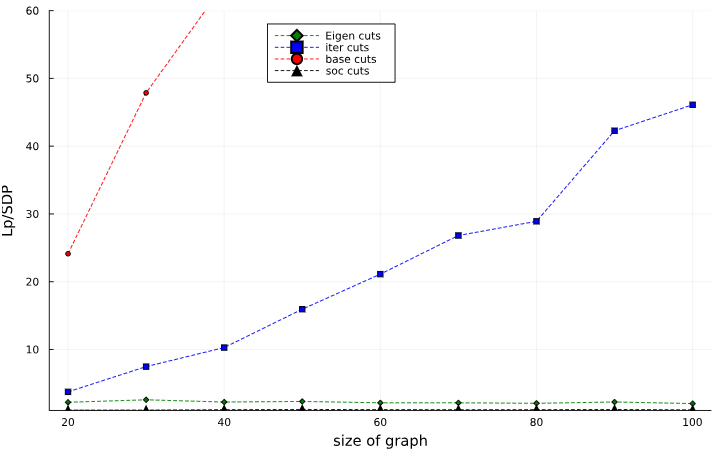}
\caption{$r=1, \ q= \frac{n}{5}$}
\end{subfigure}%
\begin{subfigure}[b]{0.6\linewidth}
\centering
\includegraphics[width=\linewidth]{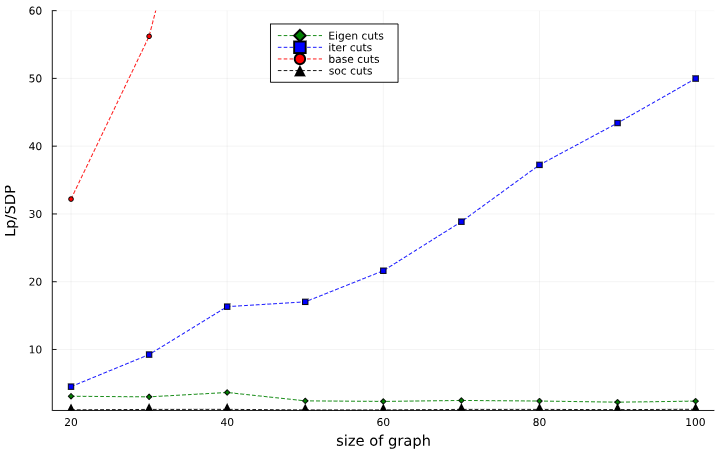}
\caption{$r=1, \ q= \frac{n}{10}$}
\end{subfigure}
\end{adjustbox}

\vspace{1cm} 

\begin{adjustbox}{max width=1.3\textwidth}
\begin{subfigure}[b]{0.6\linewidth}
\centering
\includegraphics[width=\linewidth]{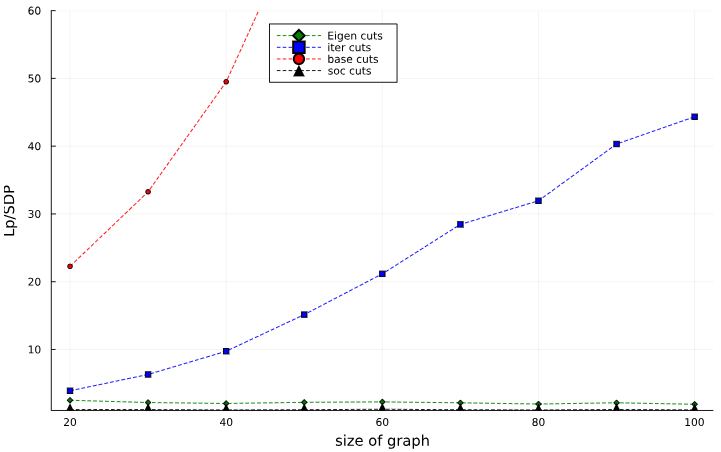}
\caption{$r=\frac{n}{2}, \ q= \frac{n}{10}$}
\end{subfigure}%
\begin{subfigure}[b]{0.6\linewidth}
\centering
\includegraphics[width=\linewidth]{revised_pictures/trust_pictures/trust_1den_rn_qn10.PNG}
\caption{$r=n, \ q= \frac{n}{10}$}
\end{subfigure}
\end{adjustbox}
\caption{Quality of the ratios $\frac{z_\mathcal{S}}{z_{sdp}}$ (eigen cuts), $\frac{z_n}{z_{sdp}}$, (oracle cuts), $\frac{z_0}{z_{sdp}}$ (base cuts) and $\frac{z_{soc}}{z_{sdp}}$ for instances of the extended trust region problem with density $1$. } 
\label{fig:trust_1}
\end{figure}

We now consider instances of the extended trust region problem with extra quadratic constraints, as presented in Subsection \ref{sec:qcqp}. These instances are the same as in the previous subsection, but with the added quadratic constraint $x^\top I_n x \leq 1$. We present our results in figures \ref{fig:trust_025}, \ref{fig:trust_05}, \ref{fig:trust_075} and \ref{fig:trust_1}.

The results for these experiments are similar across the different densities. In all of our experiments, the second order cone relaxation and the linear relaxation \ref{LSDP} of the extended trust region problem are very strong with the ratio to the SDP relaxation being very close to $1$. Moreover, this ratio does not get worse as $n$ increases, quite in sharp contrast to the base relaxation $Iter_0$ of objective value $z_0$ and the LP $Iter_n$, which gets a ratio worse than $50$ whenever $n$ exceeds $100$. for these instances, program \ref{LSDP} specialized to the extended trust region problem and program \ref{soc_strenghen_shor_quadnapsack}
certify the dual bounds provided by Theorem \ref{teo:dual_bound}, which we believe is the reason of the effectiveness of these relaxations. 

\subsubsection{Quadratic knapsack problem}

We now consider instances of the quadratic knapsack problem as presented in Subsection \ref{sec:qcqp}. In this family of problems, the linear term $b_0$ in the objective is $0$, and therefore we can consider the different strategies mentioned in Section \ref{sec:qcqp}. Hence, for each instance we solve $5$ programs, as follows:

\begin{itemize}
    \item Problem \ref{strenghen_shor_quadnapsack}. We denote the objective value of this semidefinite program by $z_{sdp}$.
    \item The linear relaxation $L_{\mathcal{S}}$ of \ref{strenghen_shor_quadnapsack} where we let $\mathcal{S}$ the elements of a eigenvector basis of the matrix $A_0$. We denote the objective value of this problem by $z_\mathcal{S}$.
    \item The LP $Iter_n$(\ref{shor_reformulated}). We denote by $z_n$ the objective value of this program.
    \item The LP $Iter_0$(\ref{shor_reformulated}). We denote by $z_0$ the objective value of this program.
    \item The second order cone relaxation of \ref{strenghen_shor_quadnapsack} given by program \ref{soc_strenghen_shor_quadnapsack}.
    
\end{itemize}

\begin{figure}[htbp]
\centering
\begin{adjustbox}{max width=1.3\textwidth}
\begin{subfigure}[b]{0.6\linewidth}
\centering
\includegraphics[width=\linewidth]{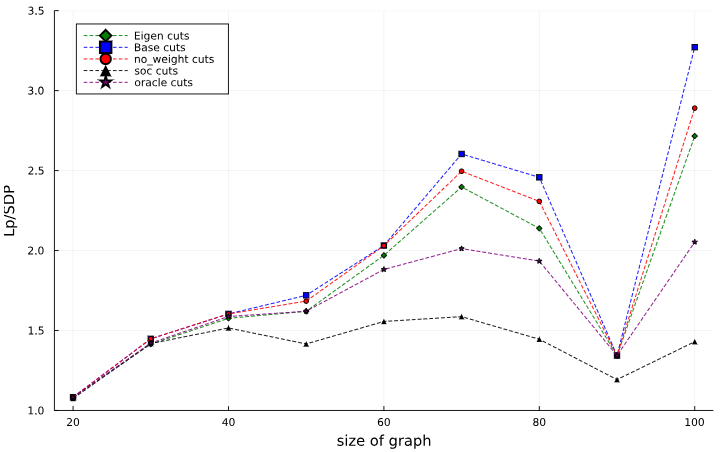}
\caption{ $\Delta = 0.05$}
\end{subfigure}%
\begin{subfigure}[b]{0.6\linewidth}
\centering
\includegraphics[width=\linewidth]{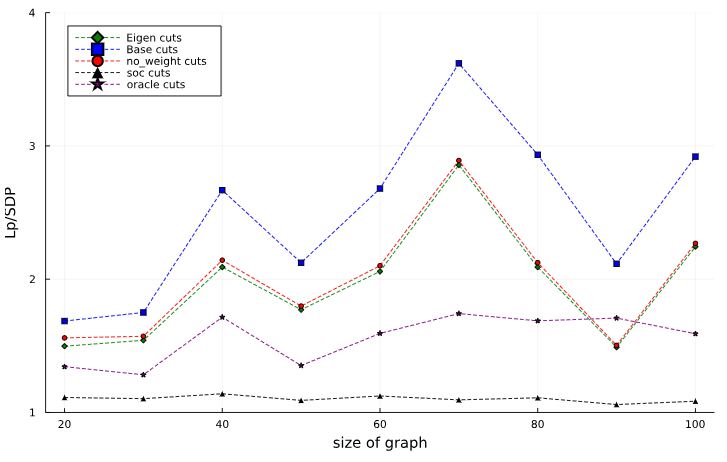}
\caption{$\Delta = 0.25$}
\end{subfigure}
\end{adjustbox}

\vspace{1cm} 

\begin{adjustbox}{max width=1.3\textwidth}
\begin{subfigure}[b]{0.6\linewidth}
\centering
\includegraphics[width=\linewidth]{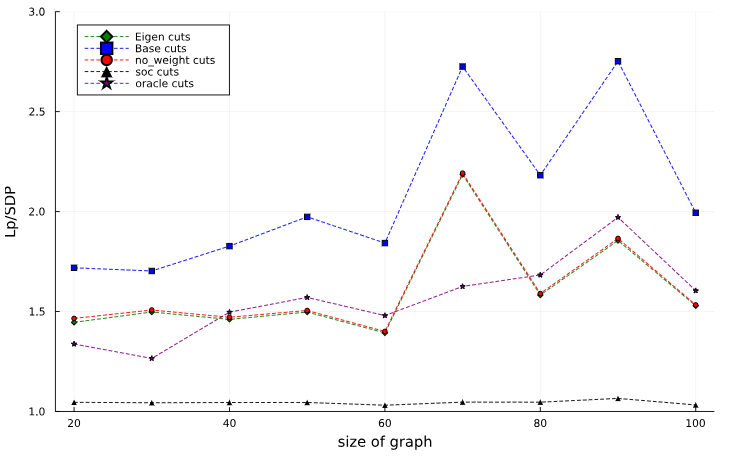}
\caption{$\Delta = 0.5$}
\end{subfigure}%
\begin{subfigure}[b]{0.6\linewidth}
\centering
\includegraphics[width=\linewidth]{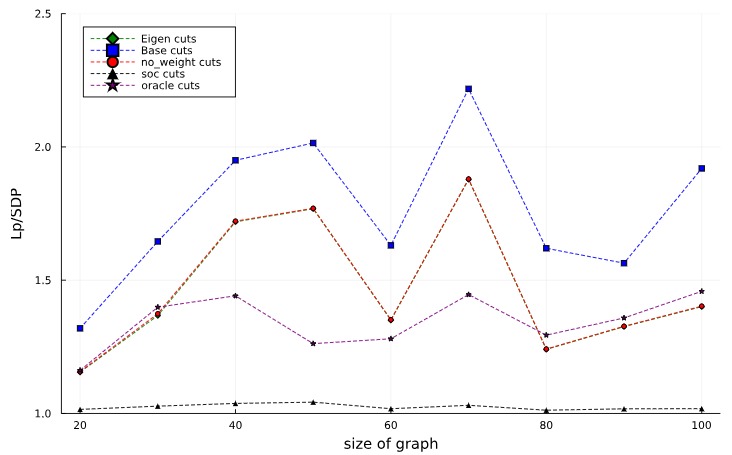}
\caption{$\Delta = 0.75$}
\end{subfigure}
\end{adjustbox}

\vspace{1cm} 

\begin{adjustbox}{max width=1.3\textwidth}
\begin{subfigure}[b]{0.6\linewidth}
\centering
\includegraphics[width=\linewidth]{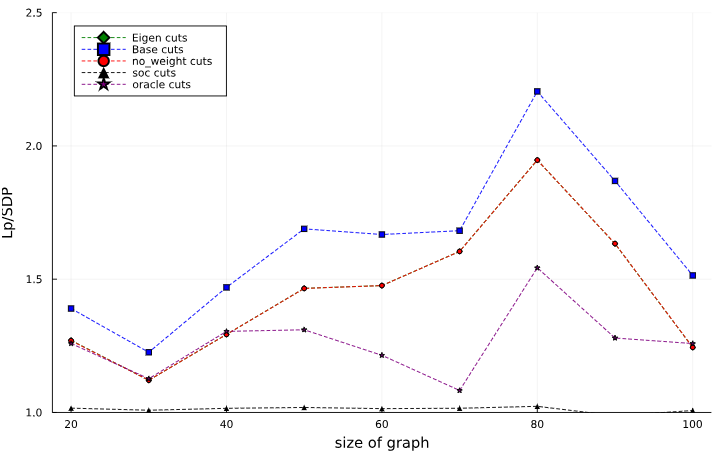}
\caption{$\Delta = 0.95$}
\end{subfigure}%
\begin{subfigure}[b]{0.6\linewidth}
\centering
\includegraphics[width=\linewidth]{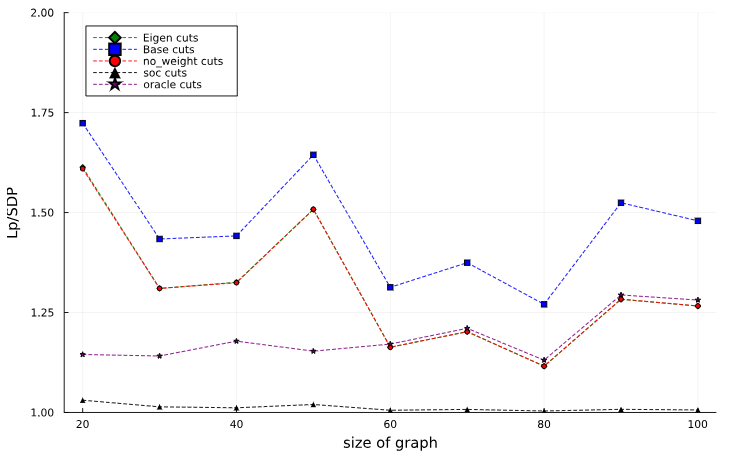}
\caption{$\Delta = 1$}
\end{subfigure}
\end{adjustbox}
\caption{Quality of the ratios $\frac{z_\mathcal{S}}{z_{sdp}}$ (eigen cuts), $\frac{z_n}{z_{sdp}}$, (oracle cuts), $\frac{z_0}{z_{sdp}}$ (base cuts) and $\frac{z_{soc}}{z_{sdp}}$ for instances of the quadratic knapsack problem with density different densities. }

\label{fig:knapsack}
\end{figure}

The instances were generated following \cite{pisinger2007quadratic}, who specify instances that have become the standard to computationally test this optimization problem. Namely, we first set a \textit{density} value $\Delta \in [0,1]$, which corresponds to the percentage of nonzero elements of the matrix $A_0$. Each weight
$w_j , \ j \in [n]$ is uniformly randomly distributed in $[1,50]$. The $ij$ entry of $A_0$ equals the $ji$ entry and is nonzero with probability $\Delta$, in which case it is uniformly distributed in $[1,100], \ i,j \in [n]$. The capacity $C$ of the knapsack is taken uniformly at random from the interval $[50,\sum_{j=1}^n w_j ]$. We present our results in Figure \ref{fig:knapsack}.

For this family of problems, all relaxations are within reasonable bounds of the SDP objective value. 
It is nonetheless appealing that the second order cone relaxation performs very well, with the ratio to the objective of the SDP nearly $1$, regardless of the value of $n$. The relaxation \ref{LSDP} seems to perform similarly to $Iter_n$.

\subsection{Computational time considerations}

Algorithm \ref{alg1} offers a meta-algorithm to solve semidefinite programs. Ideally, choosing appropriate starting sets $\mathcal{S}$ to initialize the algorithm will result in better solving times. It is critical then that solving program \ref{LSDP} or a second order cone strenghening takes significantly less time than solving the SDP. In what follows, we report solving times of the different programs proposed.

For the max cut and the Lov\'asz theta number we consider Erd\H{o}s-R\'enyi random graphs on $270$ and $200$ vertices respectively. The probability of adding an edge between two vertices is set to $p=0.75$. We repeat the experiments for $3$ instances and report the average solving time and worst ratio of the LP to the SDP objective value among the three instances.

\begin{itemize}
    \item[Max cut]: The worst ratio found was $1.08$. The average solving time of the SDP was $0.47$ seconds. The average solving time of the LP was $9.77$ seconds.
    \item[Theta number]: The worst ratio found was $6.7$. The average solving time of the SDP was $2994$ seconds. The average solving time of the LP was $39$ seconds.
    
\end{itemize}

We proceed by reporting the solving times 
for the quadratic knapsack, random QCQPs and the Extended Trust Region problem. We consider problems with $270$ variables. For the Trust Region and random QCQPs we set the number of quadratic constraints to $10$, and the number of linear constraints to $20$. For each problem, we generate $3$ instances as described previously, setting the density $\Delta$ to $0.75$. We report the average solving time, worst ratio of the LP to the SDP objective and worst ratio of the SOC to the SDP value among the three instances.

\begin{itemize}
    \item[Trust region]: The average solving time of the SDP was $7476$ seconds. The average solving time of the LP was $216$ seconds. The average solving time of the SOC was $8.9$ seconds. The worst ratio found for the LP was $2.49$, and the worst ratio found for the SOC was $1.17$.
    \item[Random QCQPs]: The average solving time of the SDP was $6510$ seconds. The average solving time of the LP was $13$ seconds. The average solving time of the SOC was $21$ seconds. The worst ratio found for the LP was $1.49$, and the worst ratio found for the SOC was $1.48$.
    \item[Knapsack]: The average solving time of the SDP was $7422$ seconds. The average solving time of the LP was $15$ seconds. The average solving time of the SOC was $20$ seconds. The worst ratio found for the LP was $1.844$, and the worst ratio found for the SOC was $1.01$.
    
\end{itemize}

It is noteworthy that solving the max cut SDP is faster by $4$ orders of magnitude than the all of the other semidefinite programs considered in this paper. In addition, it is quite surprising that the SOC relaxations of the QCQPs have solving times comparable to that of the LPs. In particular, the solving time of the SOC is two orders of magnitude faster than the LP for the trust region problems. 
We point out that very strong, fast and scalable, specialized algorithms for semidefinite programs such as the max cut problem and the Lov\'az theta number exist, such as \cite{garstka2021cosmo,wang2023decomposition,wen2013feasible}, and therefore alternatives such as an outer approximation algorithm as \ref{alg1} might not be appealing for these problems.

\section{Summary and future work}
In this work, we introduced a generic technique to obtain linear and second order cone relaxations of semidefinite programs with provable guarantees based on the commutativity of the constraints and objective matrices. 
We believe that other algebraic properties of these matrices can be exploited to obtain further stronger relaxations. 
Although we believe solving semidefinite programs with linear programs is an interesting topic is its own right, we posit that our ideas can be exploited in settings where linear approximations of convex regions is an essential component of state-of-the-art algorithms, such as in copositive programming \cite{bundfuss2009adaptive} and outer approximation algorithms for semidefinite integer programs \cite{lubin2016extended}.

On the theoretical side, the main remaining question regarding the max cut problem is if the proposed linear program \ref{Srel} provides a better-than-$2$ approximation algorithm. From our computational tests, we are not aware of any instance where the approximation factor is worse than $1.8$.

For the Lov\'asz theta number, the main theoretical question is if our proposed linear program satisfies the 
same inequalities that $\vartheta(\Bar{G})$ does. Namely,
\[
\alpha(G) \leq \vartheta(\Bar{G}) \leq \chi(G)
\]
where $\alpha(G) $ and $ \chi(G)$ are the clique and chromatic numbers of $G$, respectively. It would be interesting as well to find out if program \ref{theta} satisfies the bound (\ref{eq:theta_eigen_bound}) for $d-$regular graphs.
Finally, the second order cone relaxations for the knapsack and extended trust region problems performed well in terms of both solving time and objective value. It would be then worthwhile to explore the specialization of Algorithm \ref{alg1} to these problems, and to compare its behaviour to state of the art algorithms for those problems.

 \appendix

 \section{Missing Proofs }
\label{AppendixA}

Here we present a proof of Observation \ref{OptimalS}.

\begin{proof}
We first describe the dual of program $L_{\mathcal{S}}$ for a generic set $\mathcal{S} =\{s_1,\dots,s_k\} $. This program is given by:

\begin{equation}\label{DLP}\tag{$DL_{\mathcal{S}}$}
\begin{aligned}
&\max_{y\in \mathbb{R}^{n}, \ \alpha \in \mathbb{R}^n_+} \  b^\top y \\
\text{ s.t: }& C-  \sum_{i=1}^ry_iA_i  = \sum_{i=1}^k\alpha_is_is_i^\top.
\end{aligned}
\end{equation}

Notice that for any set $\mathcal{S} = \{s_1,\dots,s_k\}$, program \eqref{DLP}
is a restriction of $(\ref{DSDP})$ as the matrices $C-\sum_{i}y_iA_i$ are restricted to belong to the convex cone generated by the PSD matrices $s_is_i^T, i\in [k]$, rather than the whole set of positive semidefinite matrices. It follows that the optimal value of (\ref{DSDP}) upper bounds the optimal value of $(DL_{\mathcal{S}})$ for any set $\mathcal{S}$.
By hypothesis, both (\ref{SDP}) and its dual are strictly feasible and therefore solvable by strong conic duality. Hence, we let $\mathcal{S}^*$ be the elements of a basis of $\mathbb{R}^n$ of orthonormal eigenvectors of an optimal solution $S^*$ of program \ref{DSDP}. The dual of $L_{\mathcal{S}^*}$ is then
$\max_{y\in \mathbb{R}^{n} \ \alpha \in \mathbb{R}^n_+} b^\top y$ subject to
$C-  \sum_{i=1}^my_iA_i  = \sum_{i=1}^k\alpha_iv_iv_i^\top$. Hence, letting $y_i = y_i^*$ and $\alpha_i = \beta_i$ gives a feasible solution to $DL_{\mathcal{S}^*}$ which matches the optimal value of $\ref{DSDP}$ and hence is optimal.
To conclude, observe that strong linear duality holds and therefore $L_{\mathcal{S}^*}$ is solvable and its optimal value equals that of 
$\ref{DSDP}$.

\end{proof}

We now present a direct proof of Theorem \ref{teo:eigenbound}.

\begin{proof}[Proof of Theorem  \ref{teo:eigenbound}]

The inequality holds if we are able to show a feasible solution of the dual program of \eqref{Srel} whose objective value equals $-n\lambda_n$. 
Let $\mathcal{S} = \{v_1,v_2,\dots,v_n\}= \mathcal{E}(W)$.
Consider an eigenvector $v$ of $W$ with corresponding eigenvalue $\lambda$, so that $Wv = \lambda v$. Observe that $\lambda-\lambda_n \geq 0$ since $\lambda_n$ is the most negative eigenvalue of $W$. Clearly, the vector $v$ is a eigenvector of the matrix $W-\lambda_nI_n$ with corresponding eigenvalue $\lambda-\lambda_n$. By the spectral theorem, we have $W - \lambda_n I_{n} = \sum_{i=1}^n (\lambda_i-\lambda_n)v_iv_i^T$, where $v_1,\dots,v_n$ are an orthonormal eigenbasis of $W$. In other words, we have:

\[
 \lambda_n I_{n} + \sum_{i=1}^n (\lambda_i-\lambda_n)v_iv_i^T = W.
\]

This yields the desired feasible solution  $\Lambda = -\lambda_nI_n$ which has an objective value  $\frac{1}{2}m-\frac{n}{4}\lambda_n$ for \ref{dSrel}.
\end{proof}

We next derive the dual of program \ref{Srel}.

\begin{lemma}
The dual of program  \ref{Srel} is given by program \ref{MCdSrel}.
\comment{
\begin{equation}\tag{$DSP_{\mathcal{S}}$}
\begin{aligned}
\min_{\lambda, \alpha,\delta, \beta,\Lambda} \ &\frac{1}{2}m - \frac14 \left  [ Tr(\Lambda) - \sum_{i\neq j}\delta_{ij} - \sum_{i\neq j}\alpha_{ij}   \right]    \\
\text{ s.t: } & W - \Lambda = \sum_{i=1}^k \beta_iv_iv_i^\top, \\
 & \delta_{ij} \geq 0 \ \forall i\neq j \in [n], \\
  & \alpha_{ij} \geq 0 \ \forall i\neq j \in [n], \\
  & \lambda_{i} \in \mathbb{R} \ \forall i \in [n], \\
   & \beta_{i} \geq 0 \ \forall i \in [n], \\
    \Lambda \in \mathbb{S}^n, \Lambda_{ij} = & \delta_{ij}-\alpha_{ij}\ \forall i\neq j \in [n],
   \Lambda_{ii} = \lambda_{i} \  \forall i \in [n],
   \end{aligned}
\end{equation}

where by $k$ the cardinality of $\mathcal{S}$.
}
\end{lemma}

\begin{proof}
For this proof we ignore the constant $\frac{1}{2}m$ in the objective together with the multiplicative term $\frac{1}{4}$. Introduce dual variables $\lambda_i \in \mathbb{R}$ for $i \in [n]$ corresponding to the constraints $X_{ii} =1$, $\beta_i \in \mathbb{R}^n_+$ for $i \in [k]$ corresponding to $v^\top X v \geq 0$ and
$\alpha_{i,j},\delta_{i,j}\geq 0$ for $i\neq j \in [n]$, corresponding to $X_{ij} \leq 1$ and $X_{ij}\geq -1$ respectively, for 
$i\neq j \in [n]$ (in fact we need only to consider the indices $i<j$ since $X$ is symmetric but we will ignore this as it only complicates the proof). Multiplying the dual variables with the constraints accordingly gives the inequality

\[
\sum_{i=1}^n \lambda_i X_{ii} + \sum_{i=1}^k \beta_i \langle X, v_iv_i^\top \rangle -\sum_{i\neq j}\alpha_{ij}X_{ij} + \sum_{i\neq j}\delta_{ij}X_{ij} \geq \sum_{i=1}^n \lambda_i -\sum_{i\neq j}\delta_{ij} - \sum_{i\neq j}\alpha_{ij} 
\]
    
Let $\Lambda_{ij} = \delta_{ij}-\alpha_{ij}$ for $i\neq j$ and $\Lambda_{ii} = \lambda_i$ for all $i\in [n]$. This gives the inequality

\[
\langle \Lambda ,X \rangle +  \sum_{i=1}^k \beta_i \langle X, v_iv_i^\top \rangle \geq \sum_{i=1}^n \Lambda_{ii} - \sum_{i\neq j} \Lambda_{ij}.
\]

If we let $\Lambda + \sum_{i=1}^k \beta_iv_iv_i^\top = W$ we get 

\[
\left \langle -W, X \right \rangle \leq \sum_{i\neq j} \Lambda_{ij} -  \sum_{i=1}^n \Lambda_{ii}
\]
and this completes the proof.
 
\end{proof}

\section{Solutions to underlying problems }
\label{AppendixB}

In Section \ref{revised_sec:5}, we focused on comparing the objective of the semidefinite relaxation to that of a LP or second order cone relaxation. However, the semidefinite relaxation \ref{GW} is the relaxation of combinatorial problem, and hence a natural question is whether the linear programs proposed give good solutions for the actual underlying problem. As far as we are aware, there is no algorithm to round the  Lov\'asz theta number to obtain 
an independent sub-graph in a graph, but we can certainly round solutions of the LPs for the max cut problem and for a problem which is not combinatorial: the sparse PCA problem. 

In what follows, we present experimental results showing the quality of actual graph cuts obtained using programs \ref{Srel} and \ref{DSDGW} and the subsequent rounding using vectors for $\mathcal{S}$ which we describe in the next subsection. We present as well the ratios
$\frac{\frac{1}{2}m+\frac{1}{4}Z_{SP_\mathcal{S}}}{\frac{1}{2}m+\frac{1}{4}Z_{SD_\mathcal{S}}}$ and $\frac{\frac{1}{2}m+\frac{1}{4}Z_{SP_\mathcal{S}}}{\frac{1}{2}m+\frac{1}{4}Z_{GW}}$ which we call \textit{LP-gap} and \textit{optimality gap}, respectively.

These experiments are done for for Erd\H{o}s-R\'enyi random graphs, $16$ graphs taken from TSPLIB, $14$ graphs from the network repository. We compare them with graph cut values obtained by Mirka and Williamson on the same graph instances \footnote{The results included here were obtained by direct communication with the authors, and will be included in a future version of \cite{mirka2022experimental}. }. For the Trevisan's algorithm see \cite{trevisan2012max}. The simple and the sweep algorithms are modifications of Trevisan's algorithm presented in \cite{mirka2022experimental}. The greedy algorithm for max cut is folklore, and the specifics are detailed in the previous reference as well. In the second subsection, we include results for the sparse PCA problem.

\subsection{Finding cuts from $SD_\mathcal{S}$}\label{cutGeneration}

 A particular advantage of program \ref{DSDGW} is that its solutions are also feasible for $GW$ and hence can be employed using the rounding algorithm in ~\cite{goemans1995improved} to obtain feasible cuts, as the next observation shows.
 
\begin{observation}\label{obsCuts}
Let $X$ be feasible for program \ref{DSDGW} where $\mathcal{S}$ is finite and contains the standard basis $e_1,\dots,e_n$ of $\mathbb{R}^n$ . Then, we can find a cut of value at least 
$0.878\left(\frac{1}{2}m+\frac{1}{4}z_{SD_{\mathcal{S}}}\right)$.
\end{observation}
\begin{proof}
Notice that since that for $i$ in $1,\dots,n$, the matrix $e_i(e_i)^{T}$ is the matrix of all zeros with a $1$ on its $ii$ entry. Further, since $W$ has $0$ on its diagonal we may assume that any optimal solution to \ref{DSDGW} has all diagonal entries equal to $1$. It follows that such a solution is feasible for the Goemans and Williamson semidefinite program, and hence we can use the rounding procedure described in \cite{goemans1995improved} to obtain a cut with the claimed value.
\end{proof}

In our experiments, we will be using different vectors for programs \ref{DSDGW} and \ref{Srel}, so to avoid any confusion we denote by $\mathcal{S}'$ the set of vectors used for the relaxation $SD_{\mathcal{S}'}$.
An interesting source of vectors for $\mathcal{S}'$ for program $SD_{\mathcal{S}'}$ are the eigenvectors of an optimal solution $\Hat{X}$ to \ref{Srel} where $\mathcal{S}= \mathcal{E}(W)$. Although $\Hat{X}$ is not PSD in general, we can take the eigenvectors $x_i, i\in [k]$ that correspond to positive eigenvalues of $\Hat{X}$, and 
let $\mathcal{S}' = \{x_i,\dots,x_k\}$.
We observe that the computational cost of this procedure comes from solving the LP (or the SDP), whereas producing a random vector in the unit ball to find a cut is computationally cheap. Hence we produce $100$ random vectors and report the value of the best cut found using those vectors in all of our experiments. We note that in \cite{mirka2022experimental} the same method is used to find cuts from a solution to \ref{GW}.

\subsection{Erd\H{o}s-R\'enyi random graphs}
Letting $\mathcal{S}=\mathcal{E}(W)$, we know, thanks to Corollary \ref{cor:1} ,that the LP gap and hence the optimality gap converge to $1$ as $n$ grows with high probability for Erd\H{o}s-R\'enyi graphs when $np$ is not very small. We empirically evaluate the size of cuts produced by \ref{DSDGW} and the subsequent rounding and present these results in table \ref{table  ErdosRenyiMaxCut}, together with the results obtained by Mirka and Williamson on the same graphs. Surprisingly, our procedure generates the best cut value on $15$ out of the $20$ instances reported in \cite{mirka2022experimental}. Furthermore, we obtain cuts better than the ones produced by the Goemans and Williamson rounding procedure -which first solves a semidefinite program- on $18$ out of the $20$ instances. Our results are reported in Table \ref{table  ErdosRenyiMaxCut}.

\begin{table}
\caption{Optimality gap, LP-gap, and other algorithms for Erd\H{o}s-R\'enyi random graphs for the max cut problem. } 
\centering
\begin{tabular}{|p{1.3cm}|p{0.8cm}|p{1.2cm}|p{0.8cm}|p{0.8cm}|p{0.8cm}|p{0.9cm}|p{0.8cm}|p{0.8cm}|p{0.8cm}|p{0.8cm}|}
\hline
\textbf{Graph}& \textbf{m}  & \textbf{Optimality gap} & \textbf{LP Gap} & \textbf{LP cut value} & \textbf{Greedy} & \textbf{Trevisan} & \textbf{Simple} & \textbf{sweep} & \textbf{GW} & \textbf{OPT} \\
\hline
$G(50,0.1)$ & $146$ & $1.076$  &   $ 1.275$ & $114$ & $104$ &$\mathbf{116}$ &$112$ & $113$ & $113$ & $116$ \\
\hline
$G(50,0.25) $ & $313$& $1.021$ &  $1.119$  & $208$ & $199$ & $\mathbf{210}$ & $209$ & $\mathbf{210}$ & $199$ & $211$ \\
\hline
$G(50,0.5) $ & $613$  & $1.019$   & $1.071$ &$ \mathbf{373}$   & $357$ & $363$ & $372$ & $\mathbf{373}$  & $371$ & $377$ \\
\hline 
$G(50,0.75) $ & $915$  & $1.014$   & $1.053$ &$ \mathbf{522}$   & $510$ & $510$ & $510$ & $520$ & $513$ & $524$ \\
\hline 
$G(100,0.1) $ & $426$ & $1.054$  & $  1.213$ & $\mathbf{304}$ & $272$ &$284$ &$294$ &$296$ & $301$& $311$ \\
\hline
$G(100,0.25) $ & $1243$ & $1.023$  & $1.102$  & $\mathbf{772}$ & $732$ & $747$ & $766$ & $766$ &$769$ & $782$\\
\hline
$G(100,0.5) $ & $2449$ &$1.013$  & $1.048$ & $\mathbf{1394}$ & $1350$ & $1356$ & $1391$ & $1391$ & $1375$ & $1416$\\
\hline
$G(100,0.75) $ & $3725$ & $1.006$  & $1.032$  & $\mathbf{2025}$ & $1978$ & $2008$  & $2014$ & $2019$ &$1982$ & $2035$ \\
\hline
$G(200,0.1)$ & $1947$ & $1.031$  & $1.151$ & $  \mathbf{1275}$&  $1204$ &$1242$ &$1257$ &$1257$ & $1267$ & $1296$\\
\hline
$G(200,0.25)$ & $4998$  & $1.015$  & $1.074$  & $\mathbf{2926}$ & $2796$  & $2892$ & $2920$ & $2922$ & $2861 $& $2975$\\
\hline
$G(200,0.5)$ & $10026$  & $1.007$  & $ 1.031$  & $5473$ & $5388$ & $5397$ & $5441$ & $5451$ & $\mathbf{5489}$ & $5542$\\
\hline
$G(200,0.75)$ & $14818$ & $1.005$  &$1.025$  & $\mathbf{7860}$ & $7731$ & $7743$  & $7848$ & $7852$ & $7835$ & $7904$ \\
\hline
$G(350,0.1)$  & $6011$ & $1.018$ & $1.105$  &  $3645$ & $3542$  & $3548$ &$3661$ &$\mathbf{3667}$& $3513$ & $3735$\\
\hline
$G(350,0.25)$  & $15260$ & $1.009$  & $1.050$  & $\mathbf{8613}$ & $8344$  & $8426$ & $8535$ & $8553$ & $8349$ & $8709$\\
\hline
$G(350,0.5)$  & $30423$ & $1.006$  & $ 1.027$  & $\mathbf{16327}$ & $16110$& $16225$  & $16253$ & $16298$ & $15904$ & $16482$\\
\hline
$G(350,0.75)$  & $45696$ & $1.003$  & $ 1.017$  & $\mathbf{23818}$  & $23613$ & $23678$ & $23791$ & $23811$ &$23674$ &  $23967$\\
\hline
$G(500,0.1)$ & $12506$  & $1.0143$  & $1.087$  & $7326$  & $7174$ & $7174$ & $7314$ &$7314$ & $\mathbf{7372}$ & $7532$\\
\hline
$G(500,0.25)$ & $31127$ & $1.007$  & $1.040$  & $\mathbf{17177}$  &  $16833$& $17014$ &$17045$ &$17075$ & $16766$ & $17399$\\
\hline
$G(500,0.5)$ & $62167$ & $1.004$  & $1.022$  & $\mathbf{32978}$ &  $32557$ & $32862$ & $32952$ & $32960$ & $32713$ & $33234$\\
\hline
$G(500,0.75)$ & $93343$  & $1.003$  & $1.014$  & $\mathbf{48326}$  & $47995$ & $47995$ & $48244$ & $48255$ & $47597$ & $48576$\\
\hline
\end{tabular}
\label{table  ErdosRenyiMaxCut}
\end{table}

\subsection{Relevant instances}

We test our algorithm on $16$ complete graphs from TSPLIB \cite{reinelt1991tsplib}, an online library of sample instances for the Travelling Salesman Problem and related graph problems. These graph are complete weighted graphs and hence we do not report the number of edges of each graph. In Table \ref{tableSpecialCompleteGraphs} we present the optimality gap and the LP gap found for these graphs. We report as well the size of the cuts obtained following the cut generation technique presented in \ref{cutGeneration}. Our algorithm finds the best cut in $7$ of $16$ instances, and a better (or equal) cut than the GW relaxation on $14$ out of the $16$ instances. 
We then present the same quotients on $14$ graph instances taken from the Network Repository \cite{nr-aaai15} in Table \ref{tableSpecialSparseGraphs}. Since of these graphs are weighted and some are not, we do not report the number of edges of each graph. Our algorithm finds the best cut in $8$ of the $14$ instances, and a better (or equal) cut than the GW relaxation on $10$ out of the $14$ instances.

\begin{table}
\caption{Optimality gap, LP-gap, and other algorithms for some graphs on the TSPLIB graph database \cite{reinelt1991tsplib} for the max cut problem. } 
\centering
\begin{tabular}{|p{1.3cm}|p{1.2cm}|p{0.8cm}|p{0.8cm}|p{0.8cm}|p{0.9cm}|p{0.8cm}|p{0.8cm}|p{0.8cm}|}\hline
\textbf{Graph}  & \textbf{Optimality gap} & \textbf{LP Gap} & \textbf{LP cut value} & \textbf{Greedy} & \textbf{Trevisan} & \textbf{Simple} & \textbf{Sweep} & \textbf{GW}  \\
\hline
bayg29  & $1.027$  &   $1.139$ & $\mathbf{4.269}\times \mathbf{10^4}$& $3.837 \times 10^4$ &$4.225 \times 10^4$ &$\mathbf{4.269} \times \mathbf{10^4}$ & $\mathbf{4.269} \times \mathbf{10^4}$ & $\mathbf{4.269} \times \mathbf{10^4}$ \\
\hline
bays29  & $1.025$ &  $1.139$  & $\mathbf{5.399}\times \mathbf{10^4}$ & $4.831\times 10^4$ &$5.393\times 10^4$ &$5.369\times 10^4$ & $\mathbf{5.399}\times \mathbf{10^4}$ & $5.386\times 10^4$  \\
\hline
berlin52  & $1.049$   & $1.186$ &$\mathbf{4.706}\times \mathbf{10^5}$   & $4.532\times 10^5$ &$4.616\times 10^5$ &$4.465\times 10^5$ & $4.681\times 10^5$ & $4.522\times 10^5$  \\
\hline 
bier127 & $1.036$  & $1.170$ & $2.323\times 10^7$ & $2.162\times 10^7$ &$2.300\times 10^7$ &$2.322\times 10^7$ & $\mathbf{2.330}\times \mathbf{10^7}$ & $2.320\times 10^7$\\
\hline
brazil58 & $1.031$  & $1.181$  & $2.313\times 10^6$& $\mathbf{2.319}\times \mathbf{10^6}$ &$\mathbf{2.319}\times \mathbf{10^6}$ &$2.315\times 10^6$ & $2.315\times 10^6$ & $2.180\times 10^6$ \\
\hline
brg180 &$1.009$  & $1.029$ & $4.563\times 10^7$ & $4.118\times 10^7$ &$\mathbf{4.616}\times \mathbf{10^7}$ &$4.531\times 10^7$ & $4.551\times 10^7$ & $4.330\times 10^7$\\
\hline
ch130 & $1.021$  & $1.127$  & $\mathbf{1.888}\times \mathbf{10^6}$& $1.777\times 10^6$ &$1.885\times 10^6$ &$\mathbf{1.888}\times \mathbf{10^6}$ & $\mathbf{1.888}\times \mathbf{10^6}$ & $1.887\times 10^6$ \\
\hline
ch150 & $1.024$  & $1.109$ & $2.525\times 10^6$& $2.500\times 10^6$ &$2.521\times 10^6$ &$\mathbf{2.526}\times \mathbf{10^6}$ & $\mathbf{2.526}\times \mathbf{10^6}$ & $2.434\times 10^6$ \\
\hline
d198 & $1.055$  & $1.279$  & $1.289\times 10^7$ & $9.635\times 10^6$ &$1.286\times 10^7$ &$1.292\times 10^7$ & $\mathbf{1.293}\times \mathbf{10^7}$ & $\mathbf{1.293}\times \mathbf{10^7}$\\
\hline
eil101 & $1.0218$  & $1.133$  & $\mathbf{1.071}\times \mathbf{10^5}$ & $1.052\times 10^5$ &$1.070\times 10^5$ &$1.063\times 10^5$ & $1.064\times 10^5$ & $1.058\times 10^5$ \\
\hline
gr120 & $1.011$  & $1.158$  & $2.156\times 10^6$ & $2.123\times 10^6$ &$2.147\times 10^6$ &$2.156\times 10^6$ & $\mathbf{2.157}\times \mathbf{10^6}$ & $2.154\times 10^6$\\
\hline
gr137 & $1.013$  &$1.192$  & $3.068\times 10^7$& $2.241\times 10^7$ &$3.044\times 10^7$ &$3.066\times 10^7$ & $\mathbf{3.070}\times \mathbf{10^7}$ & $\mathbf{3.070}\times \mathbf{10^7}$\\
\hline
gr202 & $1.030$ & $1.180$  &  $\mathbf{1.599}\times \mathbf{10^7}$& $1.372\times 10^7$ &$1.533\times 10^7$ &$1.559\times 10^7$ & $1.593\times 10^7$ & $1.581\times 10^7$ \\
\hline
gr96 & $1.022$  & $1.130$  & $1.165\times 10^7$ & $8.967\times 10^6$ &$1.156\times 10^7$ &$\mathbf{1.166}\times \mathbf{10^7}$ & $\mathbf{1.166}\times \mathbf{10^7}$ & $1.157\times 10^7$\\
\hline
kroA100 & $1.007$  & $1.156$  & $\mathbf{5.897}\times \mathbf{10^6}$ & $5.848\times 10^6$ &$5.850\times 10^6$ &$\mathbf{5.897}\times \mathbf{10^6}$ & $\mathbf{5.897}\times \mathbf{10^6}$ & $\mathbf{5.897}\times \mathbf{10^6}$\\
\hline
a280 & $1.018$  & $1.138$  & $3.209 \times 10^6$& $2.447\times 10^6$ &$3.151\times 10^6$ &$\mathbf{3.21}\times \mathbf{10^6}$ & $\mathbf{3.21}\times \mathbf{10^6}$ & $2.970\times 10^6$\\
\hline
\end{tabular}
\label{tableSpecialCompleteGraphs}
\end{table}

\comment{
\begin{table}[H]
\caption{ Size of cuts obtained via the generation method of Observation \ref{obsCuts} on a subset of TSPLIB graph database. } 
\centering
\begin{tabular}{|p{1cm}|p{1.2cm}|p{1.2cm}|p{1.5cm}|}
\hline
Graph  & Greedy  & GW   & cut value \\
\hline
bayg29  & $1.027$  &   $1.139$ & $4.269\times 10^4$   \\
\hline
bays29  & $1.025$ &  $1.139$  & $5.399\times 10^4$  \\
\hline
berlin52  & $1.049$   & $1.186$ &$4.706\times 10^5$    \\
\hline 
bier127 & $1.036$  & $1.170$ & $2.323\times 10^7$ \\
\hline
brazil58 & $1.031$  & $1.181$  & $2.313\times 10^6$ \\
\hline
brg180 &$1.009$  & $1.029$ & $4.563\times 10^7$ \\
\hline
ch130 & $1.021$  & $1.127$  & $1.888\times 10^6$ \\
\hline
ch150 & $1.024$  & $1.109$ & $2.525\times 10^6$ \\
\hline
d198 & $1.055$  & $1.279$  & $1.289\times 10^7$ \\
\hline
eil101 & $1.0218$  & $1.133$  & $1.071\times 10^5$ \\
\hline
gr120 & $1.011$  & $1.158$  & $2.156\times 10^6$ \\
\hline
gr137 & $1.013$  &$1.192$  & $3.068\times 10^7$ \\
\hline
gr202 & $1.030$ & $1.180$  &  $1.599\times 10^7$ \\
\hline
gr96 & $1.022$  & $1.130$  & $1.165\times 10^7$ \\
\hline
kroA100 & $1.007$  & $1.156$  & $5.897\times 10^6$ \\
\hline
a280 & $1.018$  & $1.138$  & $3.209 \times 10^6$ \\
\hline
\end{tabular}
\label{tableSpecialCompleteGraphsCuts}
\end{table}
}

\begin{table}
\caption{Optimality gap, LP-gap, and other algorithms for some graphs of the Network repository graph database \cite{nr-aaai15} for the max cut problem.} 
\centering
\begin{tabular}{|p{1.3cm}|p{1.2cm}|p{0.8cm}|p{0.8cm}|p{0.8cm}|p{1cm}|p{0.8cm}|p{0.8cm}|p{0.8cm}|p{0.8cm}|}
\hline
\textbf{Graph}  & \textbf{Optimality gap} & \textbf{LP Gap} & \textbf{LP cut value} & \textbf{Greedy} & \textbf{Trevisan} & \textbf{Simple} & \textbf{Sweep} & \textbf{GW} \\
\hline
ENZYMES8  & $1.034$  &   $1.269$ & $1.230\times 10^2$ & $1.170\times 10^2$ &$ \mathbf{1.260}\times  \mathbf{10^2}$ &$ \mathbf{1.260}\times  \mathbf{10^2}$ & $ \mathbf{1.260}\times  \mathbf{10^2}$ & $ \mathbf{1.260}\times  \mathbf{10^2}$ \\
\hline
eco-stmarks  & $1.095$ &  $1.393$  & $ \mathbf{1.756}\times  \mathbf{10^3}$& $8.891\times 10^2$ &$1.190\times 10^3$ &$9.354\times 10^2$ & $9.354\times 10^2$ & $9.601\times 10^2$ \\
\hline
johnson16-2-4  & $1.000$   & $1.000$ &$3.012\times 10^3$ & $ \mathbf{3.036}\times  \mathbf{10^3}$ &$ \mathbf{3.036}\times  \mathbf{10^3}$ &$2.958\times 10^3$ & $2.986\times 10^3$ & $2.918\times 10^3$ $116$   \\
\hline 
hamming6-2 & $1.000$  & $1.000$ & $ \mathbf{9.920}\times  \mathbf{10^2}$& $ \mathbf{9.920}\times  \mathbf{10^2}$ & $\mathbf{9.920}\times  \mathbf{10^2}$ &$9.680\times 10^2$ & $9.690$ & $9.760\times 10^2$ \\
\hline
ia-infect-hyper & $1.020$  & $1.081$  & $\mathbf{1.254}\times \mathbf{10^3}$& $1.213\times 10^3$ &$1.233\times 10^3$ &$1.227\times 10^3$ & $1.227$ & $1.211\times 10^3$\\
\hline
soc-dolphins &$1.090$  & $1.279$ & $1.160\times 10^2$& $1.120\times 10^2$ &$1.120\times 10^2$ &$1.190\times 10^2$ & $\mathbf{1.210}\times \mathbf{10^2}$ & $1.150\times 10^2$ \\
\hline
email-enron-only & $1.113$  & $1.279$  & $4.060\times 10^2$& $3.920\times 10^2$ &$\mathbf{4.130}\times \mathbf{10^2}$ &$3.710\times 10^2$ & $3.800\times 10^2$ & $3.960\times 10^2$\\
\hline
dwt\_209 & $1.054$  & $1.176$ & $\mathbf{5.410}\times \mathbf{10^2}$ & $5.250\times 10^2$ &$5.270\times 10^2$ &$5.250\times 10^2$ & $5.270$ & $5400\times 10^2$\\
\hline
inf-USAir97 & $1.332$  & $1.683$  & $1.011\times 10^2$ & $9.961\times 10^1$ &$9.820\times 10^1$ &$8.184\times 10^1$ & $9.337\times 10^1$ & $\mathbf{1.074}\times \mathbf{10^2}$\\
\hline
ca-netscience & $1.180$  & $1.334$  & $5.750\times 10^2$ & $5.830\times 10^2$ &$5.880\times 10^2$ &$5.270\times 10^2$ & $5.270\times 10^2$ & $\mathbf{6.110}\times \mathbf{10^2}$\\
\hline
ia-infect-dublin& $1.110$  & $1.247$  & $\mathbf{1.673}\times \mathbf{10^3}$ & $1.648\times 10^3$ &$1.659\times 10^3$ &$1.550\times 10^3$ & $1.558\times 10^3$ & $1.664\times 10^3$\\
\hline
road-chesapeake & $1.106$  &$ 1.313$  & $\mathbf{1.250}\times \mathbf{10^2}$ & $1.230\times 10^2$ &$1.230\times 10^2$ &$1.210\times 10^2$ & $1.230\times 10^2$ & $\mathbf{1.250}\times \mathbf{10^2}$\\
\hline
Erd\H{o}s991 & $1.294$ & $1.560$  &  $\mathbf{9.610}\times \mathbf{10^2}$ & $9.330\times 10^2$ &$9.340\times 10^2$ &$7.350\times 10^2$ & $7.580\times 10^2$ & $9.240\times 10^2$ \\
\hline
dwt\_503& $1.049$  & $1.174$  & $1.805\times 10^3$& $1.822\times 10^3$ &$1.822\times 10^3$ &$\mathbf{1.921}\times \mathbf{10^3}$ & $\mathbf{1.921}\times \mathbf{10^3}$ & $1.909\times 10^3$\\
\hline
\end{tabular}
\label{tableSpecialSparseGraphs}
\end{table}

\comment{
\begin{table}[H]
\caption{ Size of cuts obtained via the generation method of Observation \ref{obsCuts} on a subset of the Network repository graph database. \cite{nr-aaai15}} 
\centering
\begin{tabular}{|p{2cm}|p{1.2cm}|p{1.2cm}|p{1.5cm}|}
\hline
Graph  & Optimality gap & LP gap & cut value \\
\hline
ENZYMES8  & $1.034$  &   $1.269$ & $1.230\times 10^2$   \\
\hline
eco-stmarks  & $1.095$ &  $1.393$  & $1.756\times 10^3$  \\
\hline
johnson16-2-4  & $1.000$   & $1.000$ &$3.012\times 10^3$    \\
\hline 
hamming6-2 & $1.000$  & $1.000$ & $9.920\times 10^2$ \\
\hline
ia-infect-hyper & $1.020$  & $1.081$  & $1.254\times 10^3$ \\
\hline
soc-dolphins &$1.090$  & $1.279$ & $1.16\times 10^2$ \\
\hline
email-enron-only & $1.113$  & $1.279$  & $4.060\times 10^2$ \\
\hline
dwt\_209 & $1.054$  & $1.176$ & $5.410\times 10^2$ \\
\hline
inf-USAir97 & $1.332$  & $1.683$  & $1.011\times 10^2$ \\
\hline
ca-netscience & $1.180$  & $1.334$  & $5.750\times 10^2$ \\
\hline
ia-infect-dublin& $1.110$  & $1.247$  & $1.673\times 10^3$ \\
\hline
road-chesapeake & $1.106$  &$ 1.313$  & $1.250\times 10^2$ \\
\hline
Erd\H{o}s991 & $1.294$ & $1.560$  &  $9.610\times 10^2$ \\
\hline
dwt\_503& $1.049$  & $1.174$  & $1.805\times 10^3$ \\
\hline
\hline
\end{tabular}
\label{tableSpecialSparseGraphsCuts}
\end{table}
}

\subsection{Sparse PCA}

Principal component analysis (PCA) is a popular tool in the statistical and machine learning literature used for dimensionality reduction, data visualisation and analysis. The core idea is to find linear combinations of the variables that correspond to directions of maximal variance, called the principal components. Finding these can be accomplished by means of a singular value decomposition. For more details about applications we refer the reader to  \cite{abdi2010principal}.
One of the main disadvantages of PCA is that the weights in the linear combination of the variables are typically non-zero, thus hindering interpretation and applicability to certain problems, such as biology or finance. In these cases it is desirable to have components that are linear combination of just a few variables.  Such components are called \textit{sparse} components, and many different techniques have been proposed to obtain them.  Cadima and Jolliffe \cite{cadima1995loading} propose an ad-hoc technique consisting in setting to $0$ loadings that are small enough. Zou, Hastie, and Tibshirani \cite{zou2006sparse} write the PCA problem as a regression optimization problem, and then impose an $\ell_1$ penalization term to encourage sparse solutions.
Following the ideas of the previous section, we relax \label{space} by dropping the constraint $X\succeq 0$ and imposing $v^\top Xv \geq 0 $ for all $v\in \mathcal{S}$ where we set $\mathcal{S} = \mathcal{E}(C)$. 
This yields the linear program 

\begin{equation}\label{Lspca}\tag{$LSPCA$}
\begin{aligned}
\max_{X\in \mathbb{S}^{n}} & \  \langle  C,X\rangle \\
\text{ s.t: } tr(X)  = 1, & \ \Vec{1}^\top |X| \Vec{1} \leq k,  \\ 
  v^\top Xv \geq & \ 0 \ \forall v\in \mathcal{S} \\
  X_{ii} \geq 0, \  X_{ii}  +X_{jj} & -2\alpha X_{ij}\geq 0 \  \forall \ i,j\in [n] \\
  -1 \geq X_{ij} & \geq 1, \ \forall i,j \in [n].
\end{aligned}
\end{equation}

The linear constraints on $X$ added on the last two lines are valid for $X$ positive semidefinite since the cone of positive semidefinite matrices is self dual, as long as $\alpha \in [0,\sqrt{2}]$. We mention that these constraints are suggested in \cite{wang2021polyhedral}.

We test the quality of our relaxation in terms of sparsity of the recovered components in the examples presented in \cite{d2004direct} and in terms of explained variance. Explained variance is the typical way to evaluate the performance of a PCA decomposition. However, we point out that there does not seem to be a consensus in the literature for what the "explained variance" for a sparse PCA decomposition is. The reason, in a nutshell, is that components recovered in the sparse case are not mutually orthogonal \cite{camacho2020all}. In this paper, the author propose a set of corrected formulas for the the sparse pca which reduce to the usual explained variance formula when the PCs are orthogonal.

\subsection{Synthetic experiments and Pit props data set}

To evaluate the recovery of sparse principal components with their semidefinite relaxation, \cite{d2004direct} use their program on a synthetic data set and on the Pit pros data set. In this subsection, we compare our linear relaxation
\ref{Lspca} to their semidefinite program by checking the sparse components that both methods produce. D'Aspermont et al.  \cite{d2004direct} generate a synthetic matrix $C$ with sparse components and empirically check that their proposed SDP can indeed recover the components. We repeat this experiment and show that the linear relaxation \ref{Lspca} obtained by setting $\mathcal{S}=\mathcal{E}(C)$ recovers as well the components. We show the results in Table \ref{tableSyntheticPCA}.
In the artificial example, three hidden factors are created:

\[
V_1 \sim \mathcal{N}(0,290), \ V_2 \sim \mathcal{N}(0,300), \ V_3 = -0.3V_1 + 0.925V_2 + \varepsilon, \ \varepsilon \sim \mathcal{N}(0,300)
\]
with $V_1$, $V_2$ and $\varepsilon$ independent. Then, $10$ observed variables are generated as follows:

\[
X_i = V_j + \varepsilon_i^j, \ \varepsilon_i^j \sim \mathcal{N}(0,1),
\]
with $j=1$ for $i=1,2,3,4$, $j=2$ for $i=5,6,7,8$ and $j=3$ for $i=9,10$ and $\{\varepsilon_i^j\}$ independent for all $i \in [10]$ and $j \in [3]$. To recover the sparse components, a solution
$X_1$ for program \ref{Lspca} is found and truncated to keep only the dominant -sparse-
eigenvector $x_1$. Then, the covariance matrix $C$ is deflated to obtain

\[
C_2 = C-(x_1^\top C x_1)x_1x_1^\top
\]
and iterated to obtain further components.
As mentioned in \cite{d2004direct}, the ideal solution is to use
$(X_1,X_2,X_3,X_4)$ for the first principal component to recover factor $V_1$ and only
$(X_5,X_6,X_7,X_8)$ for the second component to recover $V_2$.  
We replicate the results of Table $1$ in  \cite{d2004direct} using the true covariance matrix $C$ and the oracle knowledge that the sparcity $k=4$. We then run our linear relaxation by setting $\mathcal{S}= \mathcal{E}(C)$. We report the results in Table \ref{tableSyntheticPCA}. Observe that the components that our linear relaxation are sparse, and have the same support that the ones found by the SDP and are very close in norm (ignoring the signs, which are irrelevant to this application).

\begin{table}[H]
\caption{Loadings for the first two principal components on the synthetic data set with $k=4$ for both PCs.}
\centering
\begin{tabular}{|rrrrrrrrrrr|}
\hline
  & $X_1$ & $X_2$  & $X_3$ & $X_4$ & $X_5$ & $X_6$  & $X_7$ & $X_8$ & $X_9$ & $X_{10}$ \\
\hline
SPCA  PC1  &  $0$  &$0$ & $0$ &$0$ & $0.5$ &$0.5$ & $0.5$ & $0.5$ & $0$ & $0$  \\
SPCA  PC2  & $0.5$ & $0.5$ & $0.5$ & $0.5$ & $0$ &$0$ & $0$ & $0$ & $0$&$0$  \\
\hline
LSPCA PC1 & $0$ &$0$ & $0$ & $0$ & $-0.598$ &$-0.596$ & $-0.457$ & $-0.28$ & $0$ & $0$ \\
LSPCA PC2 & $ 0.482$ & $0.366$ & $0.762$ & $0.226$ & $0$ &$0$ & $0$ & $0$ & $0$&$0$ \\
 \hline
\end{tabular}
\label{tableSyntheticPCA}
\end{table}

\subsubsection{Pit props dataset}

We next consider the Pit pros dataset, introduced in \cite{jeffers1967two}. This dataset consists of $180$ observations of $13$ measured variables. It is a regularly used dataset in the PCA literature, and is notorious for having hard- to-interpret principal components. 
We replicate the results of Table $2$ in \cite{d2004direct}, where they present two sets of experiments. First, they set $k=5$ for the first component and then $k=2$ for components $2$ and $3$. In the second set of experiments, they set $k=6$ for the first component and then $k=2$ for components $2$ and $3$. For our linear relaxation, We let $\mathcal{S}= \mathcal{E}(C)$ and use the same values of $k$. We report the results in Table \ref{tablePitPros} for $k=5,2,2$. These two tables show that our methods does recover sparse components.

\begin{table}[H]
\caption{Loadings for first three principal components, for the Pit props dataset, $k=5,2,2$.}
\centering
    \addtolength{\rightskip}{-2cm}
\begin{tabular}{|lrrrrrrrrrrrrr|}
\hline
 & topdiam & length & moist & testsg &  ovensg & ringtop  &  ringbud & bowmax &  bowdist & whorls & clear &knots &  diaknot\\
\hline
SPCA PC1  &  $0.56$  &$ 0.583$ & $0$ &$0$ & $0$ &$0$ & $0.263$ & $ 0.098$ & $0.371$ & $0.362$  & 0 &0 &0   \\
SPCA PC2  & $0$ & $0$ & $ 0.707$ & $ 0.707$ & $0$ &$0$ & $0$ & $0$ & $0$&$0$ & $0$& $0$& $0$  \\
SPCA PC3  & $0$ & $0$ & $0$ & $0$ & $0$ &$0.793$ & $ 0.610$ & $0$ & $0$&$0$ & $0$ & $0$ &$-0.012$ \\
\hline
LSPCA PC1 & $0.645$ &$0.437$ & $0$ & $0$ & $0$ &$0.019$ & $0.268$ & $0$ & $0.352$ & $ 0.443$ & $0$& $0$ & $0$\\
LSPCA PC2 & $0$ & $0$ & $0.640$ & $0.767$ & $0$ &$0$ & $0$ & $0$ & $0$&$0$ & $0$ &$0$ & $0$ \\
LSPCA PC2 & $0$ & $0$ & $0$ & $0.0464$ & $0$ &$0.931$ & $0.359$ & $0$ & $0$&$0$ & $0$  & $0$ & $0$ \\
\hline
\end{tabular}
\label{tablePitPros}
\end{table}

\subsubsection{Variance explained}

\begin{figure}
    \centering
    \includegraphics[scale=0.9]{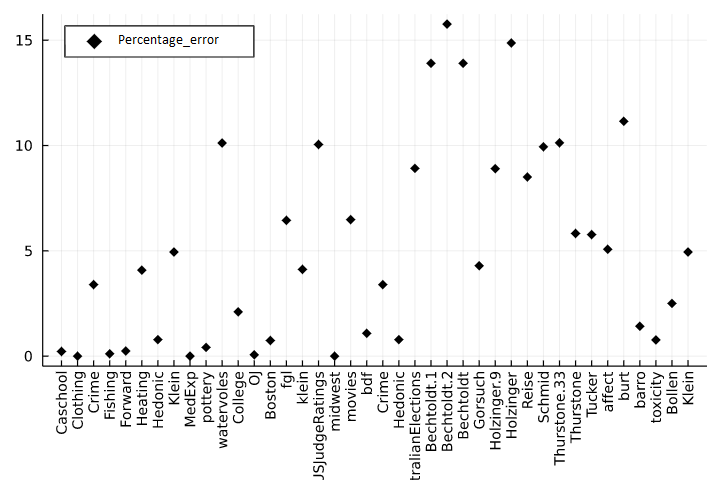}
    \caption{Relative error in percentages between the explained variances for the two SDP method and the LP method to recover sparse components. }
    \label{explainedVariance}
 \end{figure}

In this subsection we evaluate the quality of our linear relaxation in terms of the variance explained by the recovered principal components, which is the typical metric to evaluate the quality of principal components. To compute the explained variance, we use the "The fraction of total variance computed" as defined in \cite{camacho2020all}, which is a corrected formula for explained variance when the components are not orthogonal. We compute this value for both the SDP relaxation and our LP relaxation for $40$ data sets contained in the \textit{Rdataset} repository, which contains the data sets preinstalled in the core of $R$ as well of some of the data sets contained in some of the most popular $R$ libraries. We use only data sets which contain between $8$ and $20$ continuous variables. For each data set, we compute $4$ sparse components. We define 
$ev(SDP)$ to be the largest explained variance by sparse components found by the SDP method by varying the target sparsity $k$ ranging from $1$ to $30\%$ of the number of variables in the data set. We define $ev(LP)$ similarly but this time using the linear relaxation. We point out that the number $k$ for which these maximal values are obtained need not be the same.
We present our results in figure \ref{explainedVariance}. Each point corresponds to the relative error (in percentages) between the explained variances for the two methods, computed as 

\[
100\cdot \frac{ev(SDP)-ev(LP)}{ev(SDP)}.
\]

We note that among the $40$ data sets used, only $1$ has an error larger than $15\%$, all but $6$ have an error larger than $10\%$ and more than half ($23$ out of $40$) have an error of less than $5\%$.

\bibliographystyle{authordate1}
\bibliography{lit.bib}   

\begin{thebibliography}{}

\bibitem[\protect\citename{Abdi \& Williams, }2010]{abdi2010principal}
Abdi, Herv{\'e}, \& Williams, Lynne~J. 2010.
\newblock Principal component analysis.
\newblock {\em Wiley interdisciplinary reviews: computational statistics}, {\bf 2}(4), 433--459.

\bibitem[\protect\citename{Ahmadi, }2012]{ahmadi2012algebraic}
Ahmadi, Amir~Ali. 2012.
\newblock Algebraic relaxations and hardness results in polynomial optimization and Lyapunov analysis.
\newblock {\em arXiv preprint arXiv:1201.2892}.

\bibitem[\protect\citename{Ahmadi \& Hall, }2017]{ahmadi2017sum}
Ahmadi, Amir~Ali, \& Hall, Georgina. 2017.
\newblock Sum of squares basis pursuit with linear and second order cone programming.
\newblock {\em Algebraic and geometric methods in discrete mathematics}, {\bf 685}, 27--53.

\bibitem[\protect\citename{Ahmadi \& Majumdar, }2019]{ahmadi2019dsos}
Ahmadi, Amir~Ali, \& Majumdar, Anirudha. 2019.
\newblock DSOS and SDSOS optimization: more tractable alternatives to sum of squares and semidefinite optimization.
\newblock {\em SIAM Journal on Applied Algebra and Geometry}, {\bf 3}(2), 193--230.

\bibitem[\protect\citename{Alizadeh, }1995]{alizadeh1995interior}
Alizadeh, Farid. 1995.
\newblock Interior point methods in semidefinite programming with applications to combinatorial optimization.
\newblock {\em SIAM journal on Optimization}, {\bf 5}(1), 13--51.

\bibitem[\protect\citename{Alon \& Sudakov, }2000]{alon2000bipartite}
Alon, Noga, \& Sudakov, Benny. 2000.
\newblock Bipartite subgraphs and the smallest eigenvalue.
\newblock {\em Combinatorics Probability and Computing}, {\bf 9}(1), 1--12.

\bibitem[\protect\citename{ApS, }2022]{mosek}
ApS, MOSEK. 2022.
\newblock {\em The MOSEK optimization toolbox for MATLAB manual. Version 10.0.}

\bibitem[\protect\citename{Baltean-Lugojan {\em et~al.}, }2019]{baltean2019scoring}
Baltean-Lugojan, Radu, Bonami, Pierre, Misener, Ruth, \& Tramontani, Andrea. 2019.
\newblock Scoring positive semidefinite cutting planes for quadratic optimization via trained neural networks.
\newblock {\em preprint: http://www. optimization-online. org/DB\_ HTML/2018/11/6943. html}.

\bibitem[\protect\citename{Bao {\em et~al.}, }2011]{bao2011semidefinite}
Bao, Xiaowei, Sahinidis, Nikolaos~V, \& Tawarmalani, Mohit. 2011.
\newblock Semidefinite relaxations for quadratically constrained quadratic programming: A review and comparisons.
\newblock {\em Mathematical programming}, {\bf 129}, 129--157.

\bibitem[\protect\citename{Barahona \& Mahjoub, }1986]{barahona1986cut}
Barahona, Francisco, \& Mahjoub, Ali~Ridha. 1986.
\newblock On the cut polytope.
\newblock {\em Mathematical programming}, {\bf 36}(2), 157--173.

\bibitem[\protect\citename{Ben-Tal \& Den~Hertog, }2014]{ben2014hidden}
Ben-Tal, Aharon, \& Den~Hertog, Dick. 2014.
\newblock Hidden conic quadratic representation of some nonconvex quadratic optimization problems.
\newblock {\em Mathematical Programming}, {\bf 143}, 1--29.

\bibitem[\protect\citename{Ben-Tal \& Nemirovski, }2001]{ben2001lectures}
Ben-Tal, Aharon, \& Nemirovski, Arkadi. 2001.
\newblock {\em Lectures on modern convex optimization: analysis, algorithms, and engineering applications}.
\newblock SIAM.

\bibitem[\protect\citename{Bertsimas \& Cory-Wright, }2020]{bertsimas2020polyhedral}
Bertsimas, Dimitris, \& Cory-Wright, Ryan. 2020.
\newblock On polyhedral and second-order cone decompositions of semidefinite optimization problems.
\newblock {\em Operations Research Letters}, {\bf 48}(1), 78--85.

\bibitem[\protect\citename{Bhardwaj {\em et~al.}, }2015]{bhardwaj2015deciding}
Bhardwaj, Avinash, Rostalski, Philipp, \& Sanyal, Raman. 2015.
\newblock Deciding polyhedrality of spectrahedra.
\newblock {\em SIAM Journal on Optimization}, {\bf 25}(3), 1873--1884.

\bibitem[\protect\citename{Braun {\em et~al.}, }2015]{braun2015approximation}
Braun, G{\'a}bor, Fiorini, Samuel, Pokutta, Sebastian, \& Steurer, David. 2015.
\newblock Approximation limits of linear programs (beyond hierarchies).
\newblock {\em Mathematics of Operations Research}, {\bf 40}(3), 756--772.

\bibitem[\protect\citename{Brouwer {\em et~al.}, }1989]{distRegGraph}
Brouwer, Andries~E., Cohen, \& Neumaier, Arnold. 1989.
\newblock {\em Distance-Regular Graphs}.
\newblock Springer-Verlag.

\bibitem[\protect\citename{Bundfuss \& D{\"u}r, }2009]{bundfuss2009adaptive}
Bundfuss, Stefan, \& D{\"u}r, Mirjam. 2009.
\newblock An adaptive linear approximation algorithm for copositive programs.
\newblock {\em SIAM Journal on Optimization}, {\bf 20}(1), 30--53.

\bibitem[\protect\citename{Cadima \& Jolliffe, }1995]{cadima1995loading}
Cadima, Jorge, \& Jolliffe, Ian~T. 1995.
\newblock Loading and correlations in the interpretation of principle compenents.
\newblock {\em Journal of applied Statistics}, {\bf 22}(2), 203--214.

\bibitem[\protect\citename{Camacho {\em et~al.}, }2020]{camacho2020all}
Camacho, J, Smilde, Age~K, Saccenti, E, \& Westerhuis, Johan~A. 2020.
\newblock All sparse PCA models are wrong, but some are useful. Part I: computation of scores, residuals and explained variance.
\newblock {\em Chemometrics and Intelligent Laboratory Systems}, {\bf 196}, 103907.

\bibitem[\protect\citename{Chan {\em et~al.}, }2016]{chan2016approximate}
Chan, Siu~On, Lee, James~R, Raghavendra, Prasad, \& Steurer, David. 2016.
\newblock Approximate constraint satisfaction requires large LP relaxations.
\newblock {\em Journal of the ACM (JACM)}, {\bf 63}(4), 1--22.

\bibitem[\protect\citename{Charikar {\em et~al.}, }2009]{charikar2009integrality}
Charikar, Moses, Makarychev, Konstantin, \& Makarychev, Yury. 2009.
\newblock Integrality gaps for Sherali-Adams relaxations.
\newblock {\em Pages  283--292 of:} {\em Proceedings of the forty-first annual ACM symposium on Theory of computing}.

\bibitem[\protect\citename{Chung \& Radcliffe, }2011]{chung2011spectra}
Chung, Fan, \& Radcliffe, Mary. 2011.
\newblock On the spectra of general random graphs.
\newblock {\em the electronic journal of combinatorics},  P215--P215.

\bibitem[\protect\citename{Coey {\em et~al.}, }2020]{coey2020outer}
Coey, Chris, Lubin, Miles, \& Vielma, Juan~Pablo. 2020.
\newblock Outer approximation with conic certificates for mixed-integer convex problems.
\newblock {\em Mathematical Programming Computation}, {\bf 12}(2), 249--293.

\bibitem[\protect\citename{Conrad, }n.d.]{conradsimultaneous}
Conrad, Keith.
\newblock Simultaneous commutativity of operators.
\newblock {\em University of Connecticut}.

\bibitem[\protect\citename{d'Aspremont {\em et~al.}, }2004]{d2004direct}
d'Aspremont, Alexandre, Ghaoui, Laurent, Jordan, Michael, \& Lanckriet, Gert. 2004.
\newblock A direct formulation for sparse PCA using semidefinite programming.
\newblock {\em Advances in neural information processing systems}, {\bf 17}.

\bibitem[\protect\citename{de~la Vega \& Kenyon-Mathieu, }2007]{de2007linear}
de~la Vega, Wenceslas~Fernandez, \& Kenyon-Mathieu, Claire. 2007.
\newblock Linear programming relaxations of maxcut.
\newblock {\em Pages  53--61 of:} {\em Proceedings of the eighteenth annual ACM-SIAM symposium on Discrete algorithms}.

\bibitem[\protect\citename{Delorme \& Poljak, }1993]{delorme1993laplacian}
Delorme, Charles, \& Poljak, Svatopluk. 1993.
\newblock Laplacian eigenvalues and the maximum cut problem.
\newblock {\em Mathematical Programming}, {\bf 62}(1-3), 557--574.

\bibitem[\protect\citename{Dey {\em et~al.}, }2021]{dey2021cutting}
Dey, Santanu~S, Kazachkov, Aleksandr~M, Lodi, Andrea, \& Munoz, Gonzalo. 2021.
\newblock Cutting plane generation through sparse principal component analysis.
\newblock {\em URL http://www. optimization-online. org/DB\_HTML/2021/02/8259. html}.

\bibitem[\protect\citename{Deza \& Laurent, }2009]{deza2009geometry}
Deza, Michel~Marie, \& Laurent, Monique. 2009.
\newblock {\em Geometry of cuts and metrics}.
\newblock  Vol. 15.
\newblock Springer.

\bibitem[\protect\citename{Fawzi, }2021]{fawzi2021polyhedral}
Fawzi, Hamza. 2021.
\newblock On polyhedral approximations of the positive semidefinite cone.
\newblock {\em Mathematics of Operations Research}, {\bf 46}(4), 1479--1489.

\bibitem[\protect\citename{Feige \& Ofek, }2005]{feige2005spectral}
Feige, Uriel, \& Ofek, Eran. 2005.
\newblock Spectral techniques applied to sparse random graphs.
\newblock {\em Random Structures \& Algorithms}, {\bf 27}(2), 251--275.

\bibitem[\protect\citename{Friedman {\em et~al.}, }1989]{friedman1989second}
Friedman, Joel, Kahn, Jeff, \& Szemeredi, Endre. 1989.
\newblock On the second eigenvalue of random regular graphs.
\newblock {\em Pages  587--598 of:} {\em Proceedings of the twenty-first annual ACM symposium on Theory of computing}.

\bibitem[\protect\citename{Frieze \& Karo{\'n}ski, }2016]{frieze2016introduction}
Frieze, Alan, \& Karo{\'n}ski, Micha{\l}. 2016.
\newblock {\em Introduction to random graphs}.
\newblock Cambridge University Press.

\bibitem[\protect\citename{Gally \& Pfetsch, }2016]{gally2016computing}
Gally, Tristan, \& Pfetsch, MARC~E. 2016.
\newblock Computing restricted isometry constants via mixed-integer semidefinite programming.
\newblock {\em preprint, submitted}.

\bibitem[\protect\citename{Gally {\em et~al.}, }2018]{gally2018framework}
Gally, Tristan, Pfetsch, Marc~E, \& Ulbrich, Stefan. 2018.
\newblock A framework for solving mixed-integer semidefinite programs.
\newblock {\em Optimization Methods and Software}, {\bf 33}(3), 594--632.

\bibitem[\protect\citename{Garstka {\em et~al.}, }2021]{garstka2021cosmo}
Garstka, Michael, Cannon, Mark, \& Goulart, Paul. 2021.
\newblock COSMO: A conic operator splitting method for convex conic problems.
\newblock {\em Journal of Optimization Theory and Applications}, {\bf 190}(3), 779--810.

\bibitem[\protect\citename{Goemans \& Williamson, }1995]{goemans1995improved}
Goemans, Michel~X, \& Williamson, David~P. 1995.
\newblock Improved approximation algorithms for maximum cut and satisfiability problems using semidefinite programming.
\newblock {\em Journal of the ACM (JACM)}, {\bf 42}(6), 1115--1145.

\bibitem[\protect\citename{Golub \& Van~Loan, }2013]{golub2013matrix}
Golub, Gene~H, \& Van~Loan, Charles~F. 2013.
\newblock {\em Matrix computations}.
\newblock JHU press.

\bibitem[\protect\citename{Helmberg \& Rendl, }2000]{helmberg2000spectral}
Helmberg, Christoph, \& Rendl, Franz. 2000.
\newblock A spectral bundle method for semidefinite programming.
\newblock {\em SIAM Journal on Optimization}, {\bf 10}(3), 673--696.

\bibitem[\protect\citename{Helmberg {\em et~al.}, }2000]{helmberg2000semidefinite}
Helmberg, Christoph, Rendl, Franz, \& Weismantel, Robert. 2000.
\newblock A semidefinite programming approach to the quadratic knapsack problem.
\newblock {\em Journal of combinatorial optimization}, {\bf 4}, 197--215.

\bibitem[\protect\citename{Henrion \& Garulli, }2005]{henrion2005positive}
Henrion, Didier, \& Garulli, Andrea. 2005.
\newblock {\em Positive polynomials in control}.
\newblock  Vol. 312.
\newblock Springer Science \& Business Media.

\bibitem[\protect\citename{Hojny \& Pfetsch, }2023]{hojny2023handling}
Hojny, Christopher, \& Pfetsch, Marc~E. 2023.
\newblock Handling symmetries in mixed-integer semidefinite programs.
\newblock {\em Pages  69--78 of:} {\em International Conference on Integration of Constraint Programming, Artificial Intelligence, and Operations Research}.
\newblock Springer.

\bibitem[\protect\citename{Horn \& Johnson, }2012]{horn2012matrix}
Horn, Roger~A, \& Johnson, Charles~R. 2012.
\newblock {\em Matrix analysis}.
\newblock Cambridge university press.

\bibitem[\protect\citename{Jeffers, }1967]{jeffers1967two}
Jeffers, John~NR. 1967.
\newblock Two case studies in the application of principal component analysis.
\newblock {\em Journal of the Royal Statistical Society: Series C (Applied Statistics)}, {\bf 16}(3), 225--236.

\bibitem[\protect\citename{Kelley, }1960]{kelley1960cutting}
Kelley, Jr, James~E. 1960.
\newblock The cutting-plane method for solving convex programs.
\newblock {\em Journal of the society for Industrial and Applied Mathematics}, {\bf 8}(4), 703--712.

\bibitem[\protect\citename{Kothari {\em et~al.}, }2021]{kothari2021approximating}
Kothari, Pravesh~K, Meka, Raghu, \& Raghavendra, Prasad. 2021.
\newblock Approximating rectangles by juntas and weakly exponential lower bounds for LP relaxations of CSPs.
\newblock {\em SIAM Journal on Computing}, {\bf 0}(0), STOC17--305.

\bibitem[\protect\citename{Krishnan \& Mitchell, }2006]{krishnan2006unifying}
Krishnan, Kartik, \& Mitchell, John~E. 2006.
\newblock A unifying framework for several cutting plane methods for semidefinite programming.
\newblock {\em Optimization methods and software}, {\bf 21}(1), 57--74.

\bibitem[\protect\citename{Lanckriet {\em et~al.}, }2004]{lanckriet2004learning}
Lanckriet, Gert~RG, Cristianini, Nello, Bartlett, Peter, Ghaoui, Laurent~El, \& Jordan, Michael~I. 2004.
\newblock Learning the kernel matrix with semidefinite programming.
\newblock {\em Journal of Machine learning research}, {\bf 5}(Jan), 27--72.

\bibitem[\protect\citename{Lasserre, }2001]{lasserre2001global}
Lasserre, Jean~B. 2001.
\newblock Global optimization with polynomials and the problem of moments.
\newblock {\em SIAM Journal on optimization}, {\bf 11}(3), 796--817.

\bibitem[\protect\citename{Locatelli, }2016]{locatelli2016exactness}
Locatelli, Marco. 2016.
\newblock Exactness conditions for an SDP relaxation of the extended trust region problem.
\newblock {\em Optimization Letters}, {\bf 10}(6), 1141--1151.

\bibitem[\protect\citename{Lov{\'a}sz, }1979]{Lovasz1979shannon}
Lov{\'a}sz, L{\'a}szl{\'o}. 1979.
\newblock On the Shannon capacity of a graph.
\newblock {\em IEEE Transactions on Information theory}, {\bf 25}(1), 1--7.

\bibitem[\protect\citename{Lubin {\em et~al.}, }2016]{lubin2016extended}
Lubin, Miles, Yamangil, Emre, Bent, Russell, \& Vielma, Juan~Pablo. 2016.
\newblock Extended formulations in mixed-integer convex programming.
\newblock {\em Pages  102--113 of:} {\em Integer Programming and Combinatorial Optimization: 18th International Conference, IPCO 2016, Li{\`e}ge, Belgium, June 1-3, 2016, Proceedings 18}.
\newblock Springer.

\bibitem[\protect\citename{Lubotzky {\em et~al.}, }1988]{lubotzky1988ramanujan}
Lubotzky, Alexander, Phillips, Ralph, \& Sarnak, Peter. 1988.
\newblock Ramanujan graphs.
\newblock {\em Combinatorica}, {\bf 8}(3), 261--277.

\bibitem[\protect\citename{Majumdar {\em et~al.}, }2019]{majumdar2019survey}
Majumdar, Anirudha, Hall, Georgina, \& Ahmadi, Amir~Ali. 2019.
\newblock A survey of recent scalability improvements for semidefinite programming with applications in machine learning, control, and robotics.
\newblock {\em arXiv preprint arXiv:1908.05209}.

\bibitem[\protect\citename{Majumdar {\em et~al.}, }2020]{majumdar2020recent}
Majumdar, Anirudha, Hall, Georgina, \& Ahmadi, Amir~Ali. 2020.
\newblock Recent scalability improvements for semidefinite programming with applications in machine learning, control, and robotics.
\newblock {\em Annual Review of Control, Robotics, and Autonomous Systems}, {\bf 3}, 331--360.

\bibitem[\protect\citename{Mirka \& Williamson, }2022]{mirka2022experimental}
Mirka, Renee, \& Williamson, David~P. 2022.
\newblock An Experimental Evaluation of Semidefinite Programming and Spectral Algorithms for Max Cut.
\newblock {\em In:} {\em 20th International Symposium on Experimental Algorithms (SEA 2022)}.
\newblock Schloss Dagstuhl-Leibniz-Zentrum f{\"u}r Informatik.

\bibitem[\protect\citename{Mohar \& Poljak, }1990]{mohar1990eigenvalues}
Mohar, Bojan, \& Poljak, Svatopluk. 1990.
\newblock Eigenvalues and the max-cut problem.
\newblock {\em Czechoslovak Mathematical Journal}, {\bf 40}(2), 343--352.

\bibitem[\protect\citename{Nemirovski, }2004]{nemirovski2004interior}
Nemirovski, Arkadi. 2004.
\newblock Interior point polynomial time methods in convex programming.
\newblock {\em Lecture notes}, {\bf 42}(16), 3215--3224.

\bibitem[\protect\citename{Nesterov {\em et~al.}, }1997]{nesterov1997semidefinite}
Nesterov, Yurii, {\em et~al.} 1997.
\newblock {\em Semidefinite relaxation and nonconvex quadratic optimization}.
\newblock Tech. rept. Universit{\'e} catholique de Louvain, Center for Operations Research and~….

\bibitem[\protect\citename{O'Donnell \& Schramm, }2018]{o2018sherali}
O'Donnell, Ryan, \& Schramm, Tselil. 2018.
\newblock Sherali--Adams Strikes Back.
\newblock {\em arXiv preprint arXiv:1812.09967}.

\bibitem[\protect\citename{O’donoghue {\em et~al.}, }2016]{o2016conic}
O’donoghue, Brendan, Chu, Eric, Parikh, Neal, \& Boyd, Stephen. 2016.
\newblock Conic optimization via operator splitting and homogeneous self-dual embedding.
\newblock {\em Journal of Optimization Theory and Applications}, {\bf 169}, 1042--1068.

\bibitem[\protect\citename{Parrilo, }2000]{parrilo2000structured}
Parrilo, Pablo~A. 2000.
\newblock {\em Structured semidefinite programs and semialgebraic geometry methods in robustness and optimization}.
\newblock California Institute of Technology.

\bibitem[\protect\citename{Parrilo, }2003]{parrilo2003semidefinite}
Parrilo, Pablo~A. 2003.
\newblock Semidefinite programming relaxations for semialgebraic problems.
\newblock {\em Mathematical programming}, {\bf 96}(2), 293--320.

\bibitem[\protect\citename{Pisinger, }2007]{pisinger2007quadratic}
Pisinger, David. 2007.
\newblock The quadratic knapsack problem—a survey.
\newblock {\em Discrete applied mathematics}, {\bf 155}(5), 623--648.

\bibitem[\protect\citename{Poljak \& Rendl, }1995]{poljak1995nonpolyhedral}
Poljak, Svatopluk, \& Rendl, Franz. 1995.
\newblock Nonpolyhedral relaxations of graph-bisection problems.
\newblock {\em SIAM Journal on Optimization}, {\bf 5}(3), 467--487.

\bibitem[\protect\citename{Poljak \& Tuza, }1994]{poljak1994expected}
Poljak, Svatopluk, \& Tuza, Zsolt. 1994.
\newblock The expected relative error of the polyhedral approximation of the max-cut problem.
\newblock {\em Operations Research Letters}, {\bf 16}(4), 191--198.

\bibitem[\protect\citename{Qualizza {\em et~al.}, }2012]{qualizza2012linear}
Qualizza, Andrea, Belotti, Pietro, \& Margot, Fran{\c{c}}ois. 2012.
\newblock Linear programming relaxations of quadratically constrained quadratic programs.
\newblock {\em Pages  407--426 of:} {\em Mixed Integer Nonlinear Programming}.
\newblock Springer.

\bibitem[\protect\citename{Ramana, }1997]{ramana1997polyhedra}
Ramana, Motakuri~V. 1997.
\newblock Polyhedra, spectrahedra, and semidefinite programming.
\newblock {\em Topics in semidefinite and interior-point methods, Fields Institute Communications}, {\bf 18}, 27--38.

\bibitem[\protect\citename{Reinelt, }1991]{reinelt1991tsplib}
Reinelt, Gerhard. 1991.
\newblock TSPLIB—A traveling salesman problem library.
\newblock {\em ORSA journal on computing}, {\bf 3}(4), 376--384.

\bibitem[\protect\citename{Rossi \& Ahmed, }2015]{nr-aaai15}
Rossi, Ryan~A., \& Ahmed, Nesreen~K. 2015.
\newblock The Network Data Repository with Interactive Graph Analytics and Visualization.
\newblock {\em In:} {\em AAAI}.

\bibitem[\protect\citename{Sherali \& Fraticelli, }2002]{sherali2002enhancing}
Sherali, Hanif~D, \& Fraticelli, Barbara~MP. 2002.
\newblock Enhancing RLT relaxations via a new class of semidefinite cuts.
\newblock {\em Journal of Global Optimization}, {\bf 22}(1), 233--261.

\bibitem[\protect\citename{Shor, }1990]{shor1990dual}
Shor, Naum~Zuselevich. 1990.
\newblock Dual quadratic estimates in polynomial and Boolean programming.
\newblock {\em Annals of Operations Research}, {\bf 25}(1), 163--168.

\bibitem[\protect\citename{Sivaramakrishnan, }2002]{sivaramakrishnan2002linear}
Sivaramakrishnan, Kartik~Krishnan. 2002.
\newblock {\em Linear programming approaches to semidefinite programming problems}.
\newblock Rensselaer Polytechnic Institute.

\bibitem[\protect\citename{Spielman, }2012]{spielman2012spectral}
Spielman, Daniel. 2012.
\newblock Spectral graph theory.
\newblock {\em Combinatorial scientific computing}, {\bf 18}.

\bibitem[\protect\citename{Tikhomirov \& Youssef, }2019]{tikhomirov2019spectral}
Tikhomirov, Konstantin, \& Youssef, Pierre. 2019.
\newblock The spectral gap of dense random regular graphs.
\newblock {\em The Annals of Probability}, {\bf 47}(1), 362--419.

\bibitem[\protect\citename{Trevisan, }2012]{trevisan2012max}
Trevisan, Luca. 2012.
\newblock Max cut and the smallest eigenvalue.
\newblock {\em SIAM Journal on Computing}, {\bf 41}(6), 1769--1786.

\bibitem[\protect\citename{Vandenberghe \& Boyd, }1996]{vandenberghe1996semidefinite}
Vandenberghe, Lieven, \& Boyd, Stephen. 1996.
\newblock Semidefinite programming.
\newblock {\em SIAM review}, {\bf 38}(1), 49--95.

\bibitem[\protect\citename{Vinzant, }2014]{vinzant2014spectrahedron}
Vinzant, Cynthia. 2014.
\newblock What is... a spectrahedron.
\newblock {\em Notices Amer. Math. Soc}, {\bf 61}(5), 492--494.

\bibitem[\protect\citename{Wang \& K{\i}l{\i}n{\c{c}}-Karzan, }2022]{wang2022tightness}
Wang, Alex~L, \& K{\i}l{\i}n{\c{c}}-Karzan, Fatma. 2022.
\newblock On the tightness of SDP relaxations of QCQPs.
\newblock {\em Mathematical Programming}, {\bf 193}(1), 33--73.

\bibitem[\protect\citename{Wang {\em et~al.}, }2023]{wang2023decomposition}
Wang, Yifei, Deng, Kangkang, Liu, Haoyang, \& Wen, Zaiwen. 2023.
\newblock A decomposition augmented lagrangian method for low-rank semidefinite programming.
\newblock {\em SIAM Journal on Optimization}, {\bf 33}(3), 1361--1390.

\bibitem[\protect\citename{Wang {\em et~al.}, }2021]{wang2021polyhedral}
Wang, Yuzhu, Tanaka, Akihiro, \& Yoshise, Akiko. 2021.
\newblock Polyhedral approximations of the semidefinite cone and their application.
\newblock {\em Computational Optimization and Applications}, {\bf 78}(3), 893--913.

\bibitem[\protect\citename{Wen \& Yin, }2013]{wen2013feasible}
Wen, Zaiwen, \& Yin, Wotao. 2013.
\newblock A feasible method for optimization with orthogonality constraints.
\newblock {\em Mathematical Programming}, {\bf 142}(1-2), 397--434.

\bibitem[\protect\citename{Yang {\em et~al.}, }2023]{yang2023inexact}
Yang, Heng, Liang, Ling, Carlone, Luca, \& Toh, Kim-Chuan. 2023.
\newblock An inexact projected gradient method with rounding and lifting by nonlinear programming for solving rank-one semidefinite relaxation of polynomial optimization.
\newblock {\em Mathematical Programming}, {\bf 201}(1-2), 409--472.

\bibitem[\protect\citename{Yang {\em et~al.}, }2015]{yang2015sdpnal+}
Yang, Liuqin, Sun, Defeng, \& Toh, Kim-Chuan. 2015.
\newblock SDPNAL+: a majorized semismooth Newton-CG augmented Lagrangian method for semidefinite programming with nonnegative constraints.
\newblock {\em Mathematical Programming Computation}, {\bf 7}(3), 331--366.

\bibitem[\protect\citename{Yonekura \& Kanno, }2010]{yonekura2010global}
Yonekura, Kazuo, \& Kanno, Yoshihiro. 2010.
\newblock Global optimization of robust truss topology via mixed integer semidefinite programming.
\newblock {\em Optimization and Engineering}, {\bf 11}, 355--379.

\bibitem[\protect\citename{Zhao {\em et~al.}, }2010]{zhao2010newton}
Zhao, Xin-Yuan, Sun, Defeng, \& Toh, Kim-Chuan. 2010.
\newblock A Newton-CG augmented Lagrangian method for semidefinite programming.
\newblock {\em SIAM Journal on Optimization}, {\bf 20}(4), 1737--1765.

\bibitem[\protect\citename{Zou {\em et~al.}, }2006]{zou2006sparse}
Zou, Hui, Hastie, Trevor, \& Tibshirani, Robert. 2006.
\newblock Sparse principal component analysis.
\newblock {\em Journal of computational and graphical statistics}, {\bf 15}(2), 265--286.

\end{thebibliography}

%
%

\end{document}